\documentclass[12pt, reqno]{amsart}
\usepackage{amsrefs}
\usepackage{mathrsfs}
\usepackage{amsfonts}
\usepackage[centertags]{amsmath}
\usepackage{amssymb}
\usepackage{amsthm}
\usepackage{graphicx}
\usepackage{caption}
\usepackage{enumerate}
\usepackage{geometry}
\usepackage{appendix}
\usepackage{color}
\usepackage[all]{xy}
\usepackage{float}
\usepackage{amsmath,amsfonts,amssymb,graphicx}
\usepackage{longtable}
\usepackage{resizegather}
\usepackage{framed}
\usepackage{tikz}
\usepackage{multirow}

\newcommand\scalemath[2]{\scalebox{#1}{\mbox{\ensuremath{\displaystyle #2}}}}
\newtheorem{theorem}{Theorem}[section]
\newtheorem{corollary}[theorem]{Corollary}
\newtheorem{lemma}[theorem]{Lemma}
\newtheorem{proposition}[theorem]{Proposition}
\newtheorem{definition}[theorem]{Definition}
\newtheorem{remark}[theorem]{Remark}
\newtheorem{example}[theorem]{Example}

\newtheorem{conjecture}[theorem]{Conjecture}

\numberwithin{equation}{section}

\DeclareMathOperator*{\diag}{diag}

\DeclareMathOperator*{\seq}{Seq}
\DeclareMathOperator*{\res}{Res}

\DeclareMathOperator*{\pr}{pr}
\DeclareMathOperator*{\sgn}{sgn}

\usetikzlibrary{arrows,automata,positioning}
\newcommand{\midarrow}{\tikz \draw[thin,-angle 45] (0,0) -- +(.25,0);}

\begin{document}

\title{Cluster algebras and snake modules}
\author{Bing Duan, Jian-Rong Li, \and Yan-Feng Luo$^\dag$}\thanks{$\dag$ Corresponding author}

\address{Bing Duan: School of Mathematics and Statistics, Lanzhou University, Lanzhou 730000, P. R. China.}
\email{duan890818@163.com}

\address{Jian-Rong Li, Dept. of Mathematics, Weizmann Institute of Science, Rehovot 7610001, Israel, and  Einstein Institute of Mathematics, The Hebrew University of Jerusalem, Jerusalem 9190401, Israel, and School of Mathematics and Statistics, Lanzhou University, Lanzhou 730000, P. R. China.}
\email{lijr07@gmail.com}

\address{Yan-Feng Luo: School of Mathematics and Statistics, Lanzhou University, Lanzhou 730000, P. R. China.}
\email{luoyf@lzu.edu.cn}
\date{}

\maketitle
\begin{abstract}
Snake modules introduced by Mukhin and Young form a family of modules of quantum affine algebras. The aim of this paper is to prove that the Hernandez-Leclerc conjecture about monoidal categorifications of cluster algebras is true for prime snake modules of types $A_{n}$ and $B_{n}$. We prove that prime snake modules are real. We introduce $S$-systems consisting of equations satisfied by the $q$-characters of prime snake modules of types $A_{n}$ and $B_{n}$. Moreover, we show that every equation in the $S$-system of type $A_n$ (respectively, $B_n$) corresponds to a mutation in the cluster algebra $\mathscr{A}$ (respectively, $\mathscr{A}'$) constructed by Hernandez and Leclerc and every prime snake module of type $A_n$ (respectively, $B_n$) corresponds to some cluster variable in $\mathscr{A}$ (respectively, $\mathscr{A}'$). In particular, this proves that the Hernandez-Leclerc conjecture is true for all prime snake modules of types $A_{n}$ and $B_{n}$.

\hspace{0.15cm}


\hspace{0.15cm}

\end{abstract}

\section{Introduction}
Let $\mathfrak{g}$ be a simple Lie algebra over the field of complex numbers and $U_q \widehat{\mathfrak{g}}$ the corresponding quantum affine algebra. Snake modules introduced by Mukhin and Young in \cites{MY12a,MY12b} are modules of quantum affine algebras. The family of snake modules contains all minimal affinizations which were introduced by Chari in \cite{C95}.

A simple $U_q \widehat{\mathfrak{g}}$-module $M$ is called \textit{real} if $M \otimes M$ is simple, see \cite{Le03}. A simple $U_q \widehat{\mathfrak{g}}$-module $M$ is called \textit{prime} if either $M$ is trivial or there are no non-trivial $U_q \widehat{\mathfrak{g}}$-modules $M_{1}$, $M_{2}$ such that $M= M_1\otimes M_2$, see \cite{CP97}.

Chari and Pressley classified all prime $U_{q} \widehat{\mathfrak{sl}}_{2} $-modules in \cite{CP91}. Some prime $U_{q} \widehat{\mathfrak{g}}$-modules including minimal affinizations were classified in \cite{CMY13} by considering certain homological properties. In \cite{MY12b}, Mukhin and Young classified all prime snake modules of types $A_{n}$ and $B_{n}$ and proved that snake modules of types $A_{n}$ and $B_{n}$ can be uniquely (up to permutation) decomposed into a tensor of prime snake modules. We show that all prime snake modules of types $A_{n}$ and $B_{n}$ are real (Theorem \ref{real snake modules}).

The theory of cluster algebras were introduced by Fomin and Zelevinsky in \cite{FZ02}. It has many applications to mathematics and physics.

Let $\mathcal{C}$ be the category of all finite-dimensional $U_q \widehat{\mathfrak{g}}$-modules. In \cite{HL10}, Hernandez and Leclerc introduced a full subcategory $\mathcal{C}_{\ell}$ ($\ell \in \mathbb{Z}_{\geq 0}$) of $\mathcal{C}$. Let $I$ be the set of vertices of the Dynkin diagram of $\mathfrak{g}$ and let $I = I_0 \sqcup I_1$ be a partition of $I$ such that every edge connects a vertex of $I_0$ with a vertex of $I_1$. For $i \in I$, let $\xi_i = 0$ if $i \in I_0$ and $\xi_i=1$ if $i \in I_1$. Every object $V$ in $\mathcal{C}_{\ell}$ satisfies: for every composition factor $S$ of $V$ and every $i\in I$, the roots of the Drinfeld polynomial $\pi_{i,S}(u)$ belong to $\{q^{-2k-\xi_i}| 0 \leq k \leq \ell \}$.

In \cite{HL10}, Hernandez and Leclerc introduced the concept of monoidal categorifications of cluster algebras. They proposed the following conjecture, see \cite[Conjecture 13.2]{HL10}, \cite[Conjecture 5.2]{HL13}, \cite[Conjecture 9.1]{Le10}.

\begin{conjecture}[\cite{HL10, Le10, HL13}] \label{Hernandez-Leclerc conjecture}
Let $\ell \in \mathbb{Z}_{\geq 1}$. The Grothendieck ring of $\mathcal{C}_{\ell}$ has a cluster algebra structure (the cluster algebra is denoted by $\mathscr{A}_{\ell}$). There is a bijection between the cluster monomials (respectively, cluster variables) in $\mathscr{A}_{\ell}$ and real simple modules (respectively, prime real simple modules) in $\mathcal{C}_{\ell}$.
\end{conjecture}

In the case of types $A_n$ and $D_4$, $\ell = 1$, Conjecture \ref{Hernandez-Leclerc conjecture} was proved in \cite{HL10}. In the case of types $ADE$, $\ell=1$, Conjecture \ref{Hernandez-Leclerc conjecture} was proved in \cite{Nak11}. The work of \cite{Nak11} was generalized to all acyclic quivers by Kimura and Qin \cite{KQ14} and Lee \cite{Lee13}. In the case of type $A_3$, $\ell=2$, Conjecture \ref{Hernandez-Leclerc conjecture} was proved in \cite{YMLZ15}. It was proved in \cite{HL13} that Conjecture \ref{Hernandez-Leclerc conjecture} is true for Kirillov--Reshetikhin modules in all types. It is shown that Conjecture \ref{Hernandez-Leclerc conjecture} is true for all minimal affinizations of types $G_2$, $A_n$ and $B_n$ in \cite{LQ17} and \cite{ZDLL16}.

In \cite{Q15}, Qin proved half of Conjecture \ref{Hernandez-Leclerc conjecture}: all cluster monomials correspond to simple modules. The other half of Conjecture \ref{Hernandez-Leclerc conjecture} is still unknown in general. Given a simple $U_{q} \widehat{\mathfrak{g}}$-module which is both prime and real. In general, it is not known how to obtain this module using a sequence of mutations starting from the initial seed constructed in \cite{HL10, HL13}.

Moreover, the classification of prime and real modules are not known in general. Recently, Lapid and Minguez \cite{LM16} classified real modules in a family of representations of the general linear group over a non-archimedean local field. Their results can be translated to the language of representation theory of the quantum affine algebra $U_q \widehat{\mathfrak{g}}$ ($\mathfrak{g}$ is of type $A_n$ and $q$ is not a root of unity) using the quantum Schur-Weyl duality \cite{CP96}, \cite{LM16}.

The aim of this paper is to prove that the other half of Conjecture \ref{Hernandez-Leclerc conjecture} is true for all prime snake modules of types $A_{n}$ and $B_{n}$. More precisely, we prove that every prime snake module is a cluster variable in some cluster algebra introduced in \cite{HL13}. To this aim, we introduce two systems of equations satisfied by the $q$-characters of prime snake modules of types $A_{n}$ and $B_{n}$. We call these systems the \textit{$S$-systems} of types $A_{n}$ and $B_{n}$ respectively. The equations in the $S$-systems of types $A_n$ and $B_n$ are of the form
\begin{align}
[\mathcal{S}_1]  [\mathcal{S}_2] = [\mathcal{S}_3] [\mathcal{S}_4] + [\mathcal{S}_5] [\mathcal{S}_6], \label{general form of an equation in S-system 1}
\end{align}
where $\mathcal{S}_i$ ($i \in \{1,2,\ldots,6\}$) is a prime snake module and $[\mathcal{S}_i]$ is the equivalence class of $\mathcal{S}_i$ in the Grothendieck ring of $\mathcal{C}$.
Moreover, $\mathcal{S}_3 \otimes \mathcal{S}_4$ and $\mathcal{S}_5 \otimes \mathcal{S}_6$ are simple (Theorem \ref{equation simple}). By Equation (\ref{general form of an equation in S-system 1}), $\mathcal{S}_1 \otimes \mathcal{S}_2$ is not simple. Therefore, some tensor products of prime snake modules are simple and some tensor products of prime snake modules are not simple.

Let $\mathscr{A}$ (respectively, $\mathscr{A}'$) be the cluster algebra for the quantum affine algebra of type $A_n$ (respectively, $B_n$) introduced in \cite{HL13}. We show that the equations in the $S$-system of type $A_n$ (respectively, $B_n$) correspond to mutations in $\mathscr{A}$ (respectively, $\mathscr{A}'$) and prime snake modules of type $A_n$ (respectively, $B_n$) correspond to some cluster variables in $\mathscr{A}$ (respectively, $\mathscr{A}'$). In particular, this proves that the Hernandez-Leclerc conjecture is true for all prime snake modules of types $A_{n}$ and $B_{n}$.

The procedure of proving that prime snake modules of type $A_n$ (respectively, $B_n$) correspond to some cluster variables in $\mathscr{A}$ (respectively, $\mathscr{A}'$) is as follows. For a prime snake module $L(S)$ with highest $l$-weight monomial $S$, we define a set (Section \ref{Fundamental segments and distinguished factors})
\[
\mathcal{FS}(S)=\{M_1, M_2, \ldots, M_q\},
\]
where every $M_i$ is the highest weight monomial of a minimal affinization or a certain simple $U_q \widehat{\mathfrak{g}}$-module. We construct a mutation sequence $\text{Seq}_1, \text{Seq}_2, \ldots, \text{Seq}_q$ for $L(S)$ (Section \ref{mutation sequences for prime snake modules}), where $\text{Seq}_i$ is the mutation sequence for the simple $U_q \widehat{\mathfrak{g}}$-module $L(M_i)$ with highest $l$-weight monomial $M_i$. Therefore, prime snake modules of type $A_n$ (respectively, $B_n$) correspond to some cluster variables in $\mathscr{A}$ (respectively, $\mathscr{A}'$).

When $M_i$ is the highest weight monomial of a minimal affinization, the mutation sequence $\text{Seq}_i$ is similar to the mutation sequence for a minimal affinization in \cite{ZDLL16}. A minimal affinization is a $U_q \widehat{\mathfrak{g}}$-module with a highest weight of the form (see Section \ref{New notation of snake modules}):
\begin{align*}
S^{(t)}_{k_{1}^{(i_{1})}, k_{2}^{(i_{2})}, \ldots, k_{m-1}^{(i_{m-1})}, k_{m}^{(i_{m})}},
\end{align*}
where $i_1 < i_2 < \ldots < i_{m-1} < i_{m}$, or $i_1 > i_2 > \ldots > i_{m-1} > i_{m}$. In \cite{ZDLL16}, the mutation sequences for minimal affinizations with $i_1 < i_2 < \ldots < i_{m-1} < i_{m}$ are in a cluster algebra $\mathscr{A}$ and the mutation sequences for minimal affinizations with $i_1 > i_2 > \ldots > i_{m-1} > i_{m}$ are in a cluster algebra $\widetilde{\mathscr{A}}$ which is dual to $\mathscr{A}$. In this paper, on the basis of \cite{ZDLL16}, we modify the mutation sequences in \cite{ZDLL16} so that the mutation sequences for all minimal affinizations of type $A_n$ (respectively, $B_n$) are in the same cluster algebra.

The paper is organized as follows. In Section \ref{Preliminaries}, we give some background information about cluster algebras and finite-dimensional representations of quantum affine algebras. In Section \ref{Snake modules in types $A_{n}$ and $B_{n}$}, we recall the definitions of snake modules and path description of $q$-characters for snake modules of types $A_{n}$ and $B_{n}$. Moreover, we show that all prime snake modules of types $A_n$ and $B_n$ are real (Theorem \ref{real snake modules}). In Section \ref{S-systems of types $A_{n}$, $B_{n}$}, we describe the $S$-systems of types $A_{n}$ and $B_{n}$. In Section \ref{relation between M-systems and cluster algebras}, we show that the Hernandez-Leclerc conjecture is true for all prime snake modules of types $A_{n}$ and $B_{n}$. In Section \ref{examples of mutation sequences for some snake modules}, we give some examples of mutation sequences for some prime snake modules. In Sections \ref{prove theorem of Section 4}, \ref{proof of S-systems}, and \ref{proof irreducible}, we prove Theorem \ref{real snake modules}, Theorem \ref{S-systems} and Theorem \ref{equation simple} respectively.

\section{Preliminaries}\label{Preliminaries}

\subsection{Cluster algebras}
Cluster algebras were invented by Fomin and Zelevinsky in \cite{FZ02}. Let $\mathbb{Q}$ be the rational field and $\mathcal{F}=\mathbb{Q}(x_{1}, x_{2}, \ldots, x_{n})$ the field of rational functions. A seed in $\mathcal{F}$ is a pair $\Sigma=({\bf y}, Q)$, where ${\bf y} = (y_{1}, y_{2}, \ldots, y_{n})$ is a free generating set of $\mathcal{F}$, and $Q$ is a quiver with vertices labeled by $1, 2, \ldots, n$. Assume that $Q$ has neither loops nor $2$-cycles. For $k=1, 2, \ldots, n$, one defines a mutation $\mu_k$ by $\mu_k({\bf y}, Q) = ({\bf y}', Q')$. Here ${\bf y}' = (y_1', \ldots, y_n')$, $y_{i}'=y_{i}$, for $i\neq k$, and
\begin{equation}
y_{k}'=\frac{\prod_{i\rightarrow k} y_{i}+\prod_{k\rightarrow j} y_{j}}{y_{k}}, \label{exchange relation}
\end{equation}
where the first (respectively, second) product in the right-hand side is over all arrows of $Q$ with target (respectively, source) $k$, and $Q'$ is obtained from $Q$ by
\begin{enumerate}
\item[(i)] adding a new arrow $i\rightarrow j$ for every existing pair of arrows $i\rightarrow k$ and $k\rightarrow j$;

\item[(ii)] reversing the orientation of every arrow with target or source equal to $k$;

\item[(iii)] erasing every pair of opposite arrows possible created by (i).
\end{enumerate}
The mutation class $\mathcal{C}(\Sigma)$ is the set of all seeds obtained from $\Sigma$ by a finite sequence of mutations. If $\Sigma'=((y_{1}', y_{2}', \ldots, y_{n}'), Q')$ is a seed in $\mathcal{C}(\Sigma)$, then the subset $\{y_{1}', y_{2}', \ldots, y_{n}'\}$ is called a \textit{cluster}, and its elements are called \textit{cluster variables}. The \textit{cluster algebra} $\mathcal{A}_{\Sigma}$ is the subring of $\mathcal{F}$ generated by all cluster variables. \textit{Cluster monomials} are monomials in the cluster variables supported on a single cluster.

In this paper, the initial seed in the cluster algebra we use is of the form $\Sigma=({\bf y}, Q)$, where ${\bf y}$ is an infinite set and $Q$ is an infinite quiver.

\begin{definition}[{\cite[Definition 3.1]{GG14}}] \label{definition of cluster algebras of infinite rank}
Let $Q$ be a quiver without loops or $2$-cycles and with a countably infinite number of vertices labeled by all integers $i \in \mathbb{Z}$. Furthermore, for each vertex $i$ of $Q$ let the number of arrows incident with $i$ be finite. Let ${\bf y} = \{y_i \mid i \in \mathbb{Z}\}$. An infinite initial seed is the pair $({\bf y}, Q)$. By finite sequences of mutations at vertices of $Q$ and simultaneous mutations of the set ${\bf y}$ using the exchange relation (\ref{exchange relation}), one obtains a family of infinite seeds. The sets of variables in these seeds are called the infinite clusters and their elements are called the cluster variables. The cluster algebra of infinite rank of type $Q$ is the subalgebra of $\mathbb{Q}({\bf y})$ generated by the cluster variables.
\end{definition}

Two quivers $Q_1$ and $Q_2$ related by a sequence of mutations are called \textit{mutation equivalent}, and we write as $Q_1 \sim Q_2$.

\subsection{Quantum affine algebras}\label{definition of quantum affine algebras}
Let $\mathfrak{g}$ be a simple Lie algebra and $I=\{1, \ldots, n\}$ the indices of the Dynkin diagram of $\mathfrak{g}$ (we use the same labeling of the vertices of the Dynkin diagram of $\mathfrak{g}$ as the one used in \cite{Car05}). Let $C=(c_{ij})_{i,j\in I}$ be the Cartan matrix of $\mathfrak{g}$, where $c_{ij}=\frac{2 ( \alpha_i, \alpha_j ) }{( \alpha_i, \alpha_i )}$. There is a matrix $D=\diag(d_{i}\mid i\in I)$ with entries in $\mathbb{Z}_{>0}$ such that $B=DC=(b_{ij})_{i,j\in I}$ is symmetric. We have $D=\diag(d_{i}\mid i\in I)$, where $d_i=1$, $i\in I$, for type $A_n$ and $d_i=2$, $i=1,\ldots, n-1$, $d_n=1$, for type $B_n$. Let $t=\max \{d_{i}\mid i\in I\}$. Then $t=1$ for type $A_{n}$ and $t=2$ for type $B_{n}$.

Let $q_i=q^{d_i}$, $i\in I$. Let $Q$ (respectively, $Q^+$) and $P$ (respectively, $P^+$) denote the $\mathbb{Z}$-span (respectively, $\mathbb{Z}_{\geq 0}$-span) of the simple roots and fundamental weights respectively. Let $\leq$ be the partial order on $P$ in which $\lambda \leq \lambda'$ if and only if $\lambda' - \lambda \in Q^+$.

Quantum groups were introduced independently by Jimbo \cite{Jim85} and Drinfeld \cite{Dri87}. Quantum affine algebras form a family of infinite-dimensional quantum groups. Let $\widehat{\mathfrak{g}}$ denote the untwisted affine algebra corresponding to $\mathfrak{g}$. In this paper, we fix a $q\in \mathbb{C}^{\times}$, not a root of unity. The quantum affine algebra $U_q\widehat{\mathfrak{g}}$ in Drinfeld's new realization, see \cite{Dri88}, is generated by $x_{i, n}^{\pm}$ ($i\in I, n\in \mathbb{Z}$), $k_i^{\pm 1}$ $(i\in I)$, $h_{i, n}$ ($i\in I, n\in \mathbb{Z}\backslash \{0\}$) and central elements $c^{\pm 1/2}$, subject to certain relations.

The algebra $U_q\mathfrak{g}$ is isomorphic to a subalgebra of $U_q\widehat{\mathfrak{g}}$. Therefore, $U_q\widehat{\mathfrak{g}}$-modules restrict to $U_q\mathfrak{g}$-modules.

\subsection{Finite-dimensional $U_q \widehat{\mathfrak{g}}$-modules and their $q$-characters} We recall some known results on finite-dimensional $U_q\widehat{\mathfrak{g}}$-modules and their $q$-characters, see \cite{CP94, CP95a, FR98, MY12a} for details.

Let $\mathcal{P}$ be the free abelian multiplicative group of monomials in infinitely many formal variables $(Y_{i, a})_{i\in I, a\in \mathbb{C}^{\times}}$. Then $\mathbb{Z}\mathcal{P} = \mathbb{Z}[Y_{i, a}^{\pm 1}]_{i\in I, a\in \mathbb{C}^{\times}}$. For each $j\in I$, a monomial $m=\prod_{i\in I, a\in \mathbb{C}^{\times}} Y_{i, a}^{u_{i, a}}$, where $u_{i, a}$ are some integers, is said to be \textit{$j$-dominant} (respectively, \textit{$j$-anti-dominant}) if $u_{j, a} \geq 0$ (respectively, $u_{j, a} \leq 0$) for all $a\in \mathbb{C}^{\times}$. A monomial is called \textit{dominant} (respectively, \textit{anti-dominant}) if it is $j$-dominant (respectively, $j$-anti-dominant) for all $j\in I$. Let $\mathcal{P}^+ \subset \mathcal{P}$ denote the set of all dominant monomials and for $i\in I$, let $\mathcal{P}^+_{i}\subset \mathcal{P}$ denote the set of all $i$-dominant monomials.

Every finite-dimensional simple $U_{q}\mathfrak{\widehat{g}}$-module is parametrized by a dominant monomial in $\mathcal{P}^+$ \cite{CP94, CP95a}. That is, for a dominant monomial $m=\prod_{i\in I, a \in \mathbb{C}^{\times}}Y_{i,a}^{u_{i,a}}$, there is a corresponding simple $U_{q}\mathfrak{\widehat{g}}$-module $L(m)$.

The $q$-character of a $U_q\widehat{\mathfrak{g}}$-module $V$ is given by
\begin{align*}
\chi_q(V) = \sum_{m\in \mathcal{P}} \dim(V_{m}) m \in \mathbb{Z}\mathcal{P},
\end{align*}
where $V_m$ is the $l$-weight space with $l$-weight $m$ \cite{FR98}. We use $\mathscr{M}(V)$ to denote the set of all monomials in $\chi_q(V)$ for a finite-dimensional $U_q\widehat{\mathfrak{g}}$-module $V$. For $m_+ \in \mathcal{P}^+$, we use $\chi_q(m_+)$ to denote $\chi_q(L(m_+))$. We also write $m \in \chi_q(m_+)$ if $m \in \mathscr{M}(\chi_q(m_+))$.

The following lemma is well-known.
\begin{lemma}\label{well-known}
Let $m_1, m_2$ be two monomials. Then $L(m_1m_2)$ is a sub-quotient of $L(m_1) \otimes L(m_2)$. In particular, $\mathscr{M}(L(m_1m_2)) \subseteq \mathscr{M}(L(m_1))\mathscr{M}(L(m_2))$.
\end{lemma}

A finite-dimensional $U_q\widehat{\mathfrak{g}}$-module $V$ is said to be \textit{special} if $\mathscr{M}(V)$ contains exactly one dominant monomial. It is \textit{anti-special} if $\mathscr{M}(V)$ contains exactly one anti-dominant monomial. It is \textit{thin} if no $l$-weight space of $V$ has dimension greater than 1. Clearly, a special or anti-special module must be simple.

The elements $A_{i, a} \in \mathcal{P}, i\in I, a\in \mathbb{C}^{\times}$, are defined by
\begin{gather}
\begin{align*}
A_{i, a} = Y_{i, aq_{i}}Y_{i, aq_{i}^{-1}} \left(\prod_{j:c_{ji}=-1}Y_{j, a}^{-1}\right) \left(\prod_{j:c_{ji}=-2}Y_{j, aq}^{-1}Y_{j, aq^{-1}}^{-1}\right) \left(\prod_{j:c_{ji}=-3}Y_{j, aq^{2}}^{-1}Y_{j, a}^{-1}Y_{j, aq^{-2}}^{-1}\right),
\end{align*}
\end{gather}
see \cite{FR98}. Let $\mathcal{Q}$ be the subgroup of $\mathcal{P}$ generated by $A_{i, a}, i\in I, a\in \mathbb{C}^{\times}$. Let $\mathcal{Q}^{\pm}$ be the monoids generated by $A_{i, a}^{\pm 1}, i\in I, a\in \mathbb{C}^{\times}$. There is a partial order $\leq$ on $\mathcal{P}$ in which
\begin{align*}
m\leq m' \text{ if and only if } m'm^{-1}\in \mathcal{Q}^{+}. 
\end{align*}
For all $m_+ \in \mathcal{P}^+$, we have $\mathscr{M}(L(m_+)) \subset m_+\mathcal{Q}^{-}$, see \cite{FM01}.

The concept of \textit{right negative} was introduced in Section 6 of \cite{FM01}.
\begin{definition}
A monomial $m$ is called \textit{right negative} if for all $a \in \mathbb{C}^{\times}$, for $L= \max\{l\in \mathbb{Z}\mid u_{i,aq^{l}}(m)\neq 0 \text{ for some i $\in I$}\}$ we have $u_{j,aq^{L}}(m)\leq 0$ for $j\in I$.
\end{definition}
For $i\in I, a\in \mathbb{C}^{\times}$, $A_{i,a}^{-1}$ is right-negative. A product of right-negative monomials is right-negative. If $m$ is right-negative and $m'\leq m$, then $m'$ is right-negative, see \cite{FM01, Her06}. All monomials in the $q$-character of a Kirillov-Reshetikhin module is right-negative except the highest $l$-weight monomial, see {\cite[Lemma 4.4]{Her06}}.

We need the following result from \cite{FM01}, \cite{HL10}.
\begin{proposition}[{\cite[Proposition 5.3]{HL10}}] \label{dominant monomials determine q-characters}
Let $V, W$ be two $U_q\widehat{\mathfrak{g}}$-modules. If $\chi_q(V)$ and $\chi_q(W)$ have the same dominant monomials with the same multiplicities, then $\chi_q(V) = \chi_q(W)$.
\end{proposition}

\subsection{$q$-characters of $U_q\widehat{\mathfrak{sl}}_2$-modules and the Frenkel-Mukhin algorithm}
We recall the results of the $q$-characters of $U_q \widehat{\mathfrak{sl}}_2$-modules which are well-understood, see \cite{CP91, FR98}.

Let $W_{k}^{(a)}$ be the simple $U_q \widehat{\mathfrak{sl}}_2$-module with highest weight monomial
\begin{align*}
X_{k}^{(a)}=\prod_{i=0}^{k-1} Y_{aq^{k-2i-1}},
\end{align*}
where $Y_a=Y_{1, a}$. Then the $q$-character of $W_{k}^{(a)}$ is given by
\begin{align}\label{q-characters of Uqsl_2 KR module}
\chi_q(W_{k}^{(a)})=X_{k}^{(a)} \sum_{i=0}^{k} \prod_{j=0}^{i-1} A_{aq^{k-2j}}^{-1},
\end{align}
where $A_a=Y_{aq^{-1}}Y_{aq}$.

For $a \in \mathbb{C}^{\times}, k\in \mathbb{Z}_{\geq 1}$, the set $\Sigma_{k}^{(a)} =\{aq^{k-2i-1}\}_{i=0, \ldots, k-1}$ is called a \textit{$q$-string}. Two $q$-strings $\Sigma_{k}^{(a)}$ and $\Sigma_{k'}^{(a')}$ are said to be in \textit{general position} if the union $\Sigma_{k}^{(a)} \cup \Sigma_{k'}^{(a')}$ is not a $q$-string or $\Sigma_{k}^{(a)} \subset \Sigma_{k'}^{(a')}$ or $\Sigma_{k'}^{(a')} \subset \Sigma_{k}^{(a)}$.

Denote by $L(m_+)$ the simple $U_q \widehat{\mathfrak{sl}}_2$-module with highest weight monomial $m_+$. Let $m_{+} \neq 1$ and $m_{+} \in \mathbb{Z}[Y_a]_{a\in \mathbb{C}^{\times}}$ be a dominant monomial. Then $m_+$ can be uniquely (up to permutation) written in the form
\begin{align*}
m_+=\prod_{i=1}^{s} \left( \prod_{b\in \Sigma_{k_i}^{(a_i)}} Y_{b} \right),
\end{align*}
where $s$ is an integer, $\Sigma_{k_i}^{(a_i)}, i=1, \ldots, s$, are strings which are pairwise in general position and
\begin{align}\label{q-characters of Uqsl_2 module}
L(m_+)=\bigotimes_{i=1}^s W_{k_i}^{(a_i)}, \qquad \chi_q(L(m_+))=\prod_{i=1}^s \chi_q(W_{k_i}^{(a_i)}).
\end{align}

Let $i \in I$. We call $n=Y_{i, a}Y_{i, aq_i^{2}} \cdots Y_{i, aq_i^{2k-2}}$ a $q_i$-string in a monomial $m$ if $n$ is a factor of $m$. We say that two $q_i$-strings $n_1$ and $n_2$ are in general position if $n_1 n_2$ is not a $q_i$-string or $n_1$ is a factor of $n_2$ or $n_2$ is a factor of $n_1$.

For $j\in I$, let
\begin{align*}
\beta_j : \mathbb{Z}[Y_{i, a}^{\pm 1}]_{i\in I; a\in \mathbb{C}^{\times}} \to \mathbb{Z}[Y_{a}^{\pm 1}]_{a\in \mathbb{C}^{\times}}
\end{align*}
be the ring homomorphism such that for all $a\in \mathbb{C}^{\times}$, $\beta_j(Y_{k, a})=1$ for $k\neq j$ and $\beta_j(Y_{j, a})=Y_{a}$.

Let $V$ be a $U_q \widehat{\mathfrak{g}}$-module. Then $\beta_i(\chi_q(V))$, $i\in I$, is the \textit{$q$-character} of $V$ considered as a $U_{q_i} \widehat{\mathfrak{sl}_2} $-module.

The Frenkel-Mukhin algorithm was introduced to compute the $q$-characters of $U_q \widehat{\mathfrak{g}}$-modules in Section $5$ of \cite{FM01}. The algorithm is based on the $q$-characters of $U_{q_i} \widehat{\mathfrak{sl}_2}$-modules. In some cases, the Frenkel-Mukhin algorithm does not return all terms in the $q$-character of a module. There are some counterexamples given in \cite{NN11}. However, the Frenkel-Mukhin algorithm produces the correct $q$-characters of modules in many cases. In particular, if a module $L(m_+)$ is special, then the Frenkel-Mukhin algorithm applied to $m_+$ produces the correct $q$-character $\chi_q(L(m_+))$, see Theorem 5.9 of \cite{FM01}.

Let $\mathcal{P}^+_{i}\subset \mathcal{P}$, $i\in I$, denote the set of all $i$-dominant monomials. We will need the following proposition.

\begin{proposition}[{\cite[Proposition 3.1]{Her05}; \cite[Proposition 5.9]{HL10}}] \label{decomposition of a q-character}
Let $V$ be a $U_q \widehat{\mathfrak{g}}$-module and fix $i \in I$. Then there is a unique
decomposition of $\chi_q(V)$ as a finite sum
\begin{align}
\chi_q(V) = \sum_{m \in \mathcal{P}^+_{i}} \lambda_m \varphi_i(m),
\end{align}
and the $\lambda_m$ are non-negative integers.
\end{proposition}
Here $\varphi_i(m)$ ($m\in \mathcal{P}^+_{i}$) is a polynomial defined as follows, see Section $5.2.1$ of \cite{HL10}. Let $m\in \mathcal{P}^+_{i}$ be an $i$-dominant monomial. Let $\overline{m}$ be the monomial obtained from $m$ by replacing $Y_{j,a}$ with $Y_a$ if $j = i$ and by $1$ if $j \neq i$. Then the $q$-character $\chi_q(L(\overline{m}))$ of the $U_q \widehat{\mathfrak{sl}}_2$-module $L(\overline{m})$ is given by (\ref{q-characters of Uqsl_2 KR module}), (\ref{q-characters of Uqsl_2 module}). Write $\chi_q(L(\overline{m})) = \overline{m}(1 + \sum_{p} \overline{M}_p)$, where the $\overline{M}_p$ are monomials in the variables $A^{-1}_{a}$ $(a \in \mathbb{C}^{\times})$. Then one sets $\varphi_i(m) := m(1 + \sum_p M_p)$ where each $M_p$ is obtained from the corresponding $\overline{M}_p$ by replacing each variable $A^{-1}_{a}$ by $A^{-1}_{i,a}$.

The following corollary follows from Proposition \ref{decomposition of a q-character}, see \cite{HL10}.
\begin{corollary}[{\cite{HL10}}] \label{Uqsl_2 arguments}
Let $m \in \mathcal{P}^{+}$ and $mM$ a monomial of $\chi_q(L(m))$, where $M$ is a monomial in the variables $A_{j,a}^{-1}$, $j \in I$. If $M$ contains no variable $A^{-1}_{i,a}$, then $mM$ is an $i$-dominant monomial and $\varphi_i(mM)$ is contained in $\chi_q(L(m))$. In particular, $\varphi_i(m)$ is contained in $\chi_q(L(m))$.
\end{corollary}

By the Frenkel-Mukhin algorithm \cite{FM01} and the formulas (\ref{q-characters of Uqsl_2 KR module}), (\ref{q-characters of Uqsl_2 module}), we have the following result which is used frequently in our proof.

\begin{lemma} \label{need to expand from right most factors}
Let $L(m_+)$ be a special module, where $m_+$ is its dominant monomial. Then every monomial in $\chi_q(m_+)$ is a monomial in some $\varphi_i(m)$, where $i \in I$ and $m$ is an $i$-dominant monomial in $\chi_q(m_+)$. The $l$-weights of the monomials in $\varphi_i(m)$ are less or equal to the $l$-weight of $m$.

From now on, we fix an $a\in \mathbb{C}^{\times}$ and for convenient we write $i_k=Y_{i,aq^k}$, $A_{i,k}=A_{i,aq^k}$ for $i \in I$, $k \in \mathbb{Z}$. 

Suppose that $\beta_i(m)={\bf i}_1 \cdots {\bf i}_p$, where ${\bf i}_j = i_{s_j} i_{s_j+2r_i} \cdots i_{s_j+2k_jr_{i}-2r_{i}}$, $s_j \in \mathbb{Z}$, $k_j \in \mathbb{Z}_{\geq 1}$, $j=1, \ldots, p$, are $q_i$-strings which are pairwise in general position and ${\bf i}_{j'} \not\subset {\bf i}_{j''}$, $j' \neq j''$, $j', j'' \in \{1, \ldots, p\}$. Then $m A_{i, s_j+2hr_i-r_i}^{-1}$ ($h \in \{1,\ldots,k_j-1\}$) is not a monomial in $\chi_q(m_+)$.
\end{lemma}

For example, in type $A_3$, $2_{-5}2_{-3}1_0 A^{-1}_{2,-4}=3_{-4}1_{-4}1_0$ is not in $\chi_q(2_{-5}2_{-3}1_0)$.

In type $B_3$, $1_{-22}2_{-16}2_{-12}3_{-7}3_{-5}2_0 A^{-1}_{2,-14}= 1_{-22}1_{-14}3_{-15}3_{-13}3_{-7}3_{-5}2_0$ is not in $\chi_q(1_{-22}2_{-16}2_{-12}3_{-7}3_{-5}2_0)$.

\section{Snake modules of types $A_{n}$ and $B_{n}$}\label{Snake modules in types $A_{n}$ and $B_{n}$}
In this section, we recall the definition of snake modules which were introduced by Mukhin and Young in \cite{MY12a, MY12b}. In the following, we assume that $\mathfrak{g}$ is of type $A_{n}$ or $B_{n}$.

\subsection{Snake positions and minimal snake positions}\label{Snake positions and minimal snake positions}
We recall the definitions of snake positions and minimal snake positions introduced in Section $4$ of \cite{MY12a} and Section $3$ of \cite{MY12b}.
A subset $\mathcal{X} \subset I \times \mathbb{Z}$ and an injective mapping $\iota: \mathcal{X} \to \mathbb{Z}\times \mathbb{Z}$ are defined as follows.
\begin{gather}
\begin{align*}
&\text{Type $A_{n}$}: \text{ Let } \mathcal{X}:=\{(i,k)\in I \times \mathbb{Z}: i-k \equiv 0 {\hskip -0.7em}\pmod 2\} \text{ and } \iota(i,k)=(i,k).\\
&\text{Type $B_{n}$}: \text{ Let } \mathcal{X}:=\{(n,2k): k\in \mathbb{Z} \}\sqcup \{(i,k)\in I \times \mathbb{Z}: i<n \text{ and } k \equiv 1{\hskip -0.7em} \pmod 2\} \text{ and }\\
&\qquad \qquad \quad \iota(i,k) =
\begin{cases}
     (2i, k),  & \text{if } i<n \text{ and } 2n+k-2i \equiv 1{\hskip -0.7em}\pmod 4,  \\
     (4n-2-2i,k), & \text{if } i<n \text{ and } 2n+k-2i \equiv 3{\hskip -0.7em} \pmod 4, \\
     (2n-1, k), & \text{if } i=n.
\end{cases}
\end{align*}
\end{gather}
For two sets $A$ and $B$, we define a mapping $\pr_1:A \times B \to A$ given by $\pr_1(a,b) = a$.

Let $(i,k) \in \mathcal{X}$. A point $(i',k')$ is said to be in \textit{snake position} with respect to $(i,k)$ if
\begin{gather}
\begin{align*}
&\text{Type $A_{n}$}: \; k'-k \geq |i'-i|+2 \ \text{and} \  k'-k \equiv |i'-i|{\hskip -0.7em} \pmod 2.  \\
&\text{Type $B_{n}$}: \; i=i'=n:\ k'-k \geq  2 \ \text{and} \ k'-k \equiv  2 {\hskip -0.7em}\pmod 4,  \\
&\qquad \qquad \quad i \neq i'=n \text{ or } i'\neq i=n: \ k'-k \geq 2|i'-i|+3  \ \text{and} \  k'-k \equiv  2|i'-i|-1 {\hskip -0.7em}\pmod 4, \\
&\qquad \qquad \quad i < n \text{ and }i' < n:\ k'-k \geq 2|i'-i|+4 \  \text{and} \ k'-k \equiv  2|i'-i|{\hskip -0.7em} \pmod 4.
\end{align*}
\end{gather}
The point $(i',k')$ is in \textit{minimal} snake position to $(i,k)$ if $k'-k$ is equal to the given lower bound.

\begin{remark}
The above condition for type $A_{n}$ is slightly different from the condition for type $A_n$ in Section $4.2$ of \cite{MY12a} and Section $3.2$ of \cite{MY12b}.
\end{remark}

\subsection{Prime snake positions}\label{prime snake positions}
Let $(i,k) \in \mathcal{X}$. A point $(i',k') \in \mathcal{X}$ is said to be in \textit{prime snake position} with respect to $(i,k)$ if
\begin{gather}
\begin{align*}
&\text{Type $A_{n}$}: \; \min \{ 2n+2-i-i', i+i' \} \geq k'-k \geq |i'-i|+2 \ \text{and} \  k'-k \equiv |i'-i|{\hskip -0.7em} \pmod 2.  \\
&\text{Type $B_{n}$}:\; i=i'=n:\ 4n-2 \geq  k'-k \geq  2 \ \text{and} \ k'-k \equiv  2 {\hskip -0.7em}\pmod 4,  \\
&{\hskip 5em} i \neq i'=n \text{ or } i'\neq i=n: \  2i'+2i-1 \geq  k'-k \geq 2|i'-i|+3  \ \text{and} \  k'-k \equiv  2|i'-i|-1 {\hskip -0.7em}\pmod 4, \\
&{\hskip 5em} i < n \text{ and }i' < n: \  2i'+2i \geq  k'-k \geq 2|i'-i|+4 \  \text{and} \ k'-k \equiv  2|i'-i| {\hskip -0.7em}\pmod 4.
\end{align*}
\end{gather}
\begin{remark}
The above condition for type $A_{n}$ is slightly different from the condition for type $A_n$ in Section $3.3$ of \cite{MY12b}.
\end{remark}

\subsection{Snakes and snake modules}\label{Snakes and snake modules}
A finite sequence $(i_{t},k_{t})$, $1 \leq t \leq T$, $T \in \mathbb{Z}_{\geq 0}$, of points in $\mathcal{X}$ is called a \textit{snake} if for all $2 \leq t \leq T$, the point $(i_{t},k_{t})$ is in snake position with respect to $(i_{t-1},k_{t-1})$ \cite{MY12a, MY12b}. It is called a \textit{minimal} (respectively, \textit{prime}) snake if all successive points are in minimal (respectively, prime) snake position \cite{MY12a, MY12b}.

The simple module $L(m)$ is called a \textit{snake module} (respectively, a \textit{minimal snake module}) if $m=\prod_{t=1}^{T} (i_{t})_{k_{t}}$ for some snake $(i_{t},k_{t})_{1 \leq t \leq T}$ (respectively, for some minimal snake $(i_{t},k_{t})_{1 \leq t \leq T})$ \cite{MY12a, MY12b}. In this case, we say $(i_{t},k_{t})_{1 \leq t \leq T}$ is the snake of $L(m)$.

\begin{theorem}[{\cite[Proposition 3.1]{MY12b}}]\label{prime snake modules}
A snake module is prime if and only if its snake is prime. Every snake module can be uniquely (up to permutation) decomposed into a tensor of prime snake modules.
\end{theorem}

Now we are ready for our main result in this section.

\begin{theorem}\label{real snake modules}
Prime snake modules are real.
\end{theorem}

We will prove Theorem \ref{real snake modules} in Section \ref{prove theorem of Section 4}.

\begin{remark}
The fact that prime snake modules of type $A_n$  are real is also proved in \cite{LM16}.
\end{remark}

\vskip 0.05in

Throughout this paper, when we write the highest $l$-weight monomial 
\[
m=(i_1)_{k_{1}} (i_2)_{k_2} \cdots (i_T)_{k_T}
\]
of a snake module $L(m)$, we always assume that $k_t$, $1 \leq t \leq T$, are in increasing order.

\subsection{Path description of $q$-characters for snake modules of types $A_{n}$ and $B_{n}$}\label{Path description of q-characters}
We review the path description of $q$-characters for snake modules of types $A_{n}$ and $B_{n}$, see Section $5$ of \cite{MY12a} and Section $6$ of \cite{MY12b} for details.

A \textit{path} is a finite sequence of points in the plane $\mathbb{R}^{2}$. We write $(j,\ell) \in p$ if $(j,\ell)$ is a point of the path $p$.

The following is the case of type $A_{n}$. For all $(i,k)\in \mathcal{X}$, let
\begin{align*}
\mathscr{P}_{i,k}=\{ & ((0,y_{0}),(1,y_{1}),\ldots,(n+1,y_{n+1})):  y_{0}=i+k, \\
&y_{n+1}=n+1-i+k, \text{ and } y_{i+1}-y_{i}\in \{1,-1\}, \  0\leq i\leq n\}.
\end{align*}

The sets $C_{p}^{\pm}$ of upper and lower corners of a path $p=((r,y_{r}))_{0\leq r \leq n+1}\in \mathscr{P}_{i,k}$ are defined as follows:
\begin{align*}
C^{+}_{p}=\{(r,y_{r})\in p: r\in I, \ y_{r-1}=y_{r}+1=y_{r+1}\},\\
C^{-}_{p}=\{(r,y_{r})\in p: r\in I, \ y_{r-1}=y_{r}-1=y_{r+1}\}.
\end{align*}

\begin{figure}[H]
\resizebox{1.0\width}{1.0\height}{
\begin{minipage}[b]{0.4\linewidth}
\centerline{
\begin{tikzpicture}
\draw[step=.5cm,gray,thin] (-0.5,5.5) grid (2,8) (-0.5,5.5)--(2,5.5);
\draw[fill] (0,8.3) circle (2pt) -- (0.5,8.3) circle (2pt) --(1,8.3) circle (2pt) --(1.5,8.3) circle (2pt);
\draw[thick] (-0.5,7)--(0,7.5)--(0.5,8)--(1,7.5)--(1.5,7)--(2,6.5);
\draw[thick] (-0.5,7)--(0,6.5)--(0.5,6)--(1,5.5)--(1.5,6)--(2,6.5);
\node [above] at (0,8.3) {$1$};
\node [above] at (0.5,8.3) {$2$};
\node [above] at (1,8.3) {$3$};
\node [above] at (1.5,8.3) {$4$};
\node [left] at (-0.5,8) {$0$};
\node [left] at (-0.5,7.5) {$1$};
\node [left] at (-0.5,7) {$2$};
\node [left] at (-0.5,6.5) {$3$};
\node [left] at (-0.5,6) {$4$};
\node [left] at (-0.5,5.5) {$5$};
\end{tikzpicture}}
\end{minipage}
\begin{minipage}[b]{0.4\linewidth}
\centerline{
\begin{tikzpicture}
\draw[step=.5cm,gray,thin] (-0.5,5.5) grid (2,8) (-0.5,5.5)--(2,5.5);
\draw[fill] (0,8.3) circle (2pt) -- (0.5,8.3) circle (2pt) --(1,8.3) circle (2pt) --(1.5,8.3) circle (2pt);
\begin{scope}[thick, every node/.style={sloped,allow upside down}]
\draw (-0.5,7)--node {\midarrow}(0,6.5);
\draw (0,6.5)--node {\midarrow}(0.5,6);
\draw (0.5,6)--node {\midarrow}(1,5.5);
\draw (1,5.5)--node {\midarrow}(1.5,6);
\draw (1.5,6)--node {\midarrow}(2,6.5);
\draw (-0.5,7)--node {\midarrow}(0,7.5);
\draw (0,7.5)--node {\midarrow}(0.5,8);
\draw (0.5,8)--node {\midarrow}(1,7.5);
\draw (1,7.5)--node {\midarrow}(1.5,7);
\draw (1.5,7)--node {\midarrow}(2,6.5);
\draw (0,7.5)--node {\midarrow}(0.5,7);
\draw (0.5,7)--node {\midarrow}(1,6.5);
\draw (1,6.5)--node {\midarrow}(1.5,6);
\draw (0,6.5)--node {\midarrow}(0.5,7);
\draw (0.5,7)--node {\midarrow}(1,7.5);
\draw (0.5,6)--node {\midarrow}(1,6.5);
\draw (1,6.5)--node {\midarrow}(1.5,7);
\end{scope}
\node [above] at (0,8.3) {$1$};
\node [above] at (0.5,8.3) {$2$};
\node [above] at (1,8.3) {$3$};
\node [above] at (1.5,8.3) {$4$};
\node [left] at (-0.5,8) {$0$};
\node [left] at (-0.5,7.5) {$1$};
\node [left] at (-0.5,7) {$2$};
\node [left] at (-0.5,6.5) {$3$};
\node [left] at (-0.5,6) {$4$};
\node [left] at (-0.5,5.5) {$5$};
\end{tikzpicture}}
\end{minipage}}
\caption{In type $A_4$: left, $\mathscr{P}_{(2,0)}$ corresponds to a rectangle; right, the paths corresponding to monomials of $\chi_q(L(2_0))$.}\label{rectangle 1}
\end{figure}
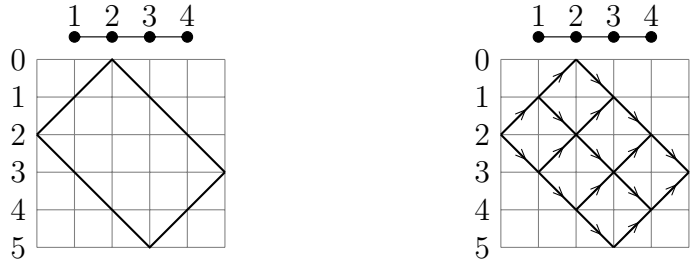

Each $\mathscr{P}_{(i,k)}$, $i\in I$, $k\in \mathbb{Z}$, corresponds to a rectangle, see Figure \ref{rectangle 1}. The monomials appearing in $\chi_q(L(i_k))$ correspond to paths in the rectangle, see Section $5$ of \cite{MY12a} and Section $6$ of \cite{MY12b} for details.

The following is the case of type $B_{n}$. Fix an $\varepsilon$ such that $0<\varepsilon <1/2$, for all $\ell \in 2\mathbb{Z}$, the set $\mathscr{P}_{n,\ell}$ is defined as follows.

For all $\ell \equiv 2 \pmod 4$,
\[
\begin{split}
\mathscr{P}_{n,\ell}= \{ & ((0,y_{0}), (2,y_{1}), \ldots, (2n-4,y_{n-2}), (2n-2,y_{n-1}),(2n-1,y_{n})): y_{0}=\ell+2n-1,\\
& y_{i+1}-y_{i}\in \{2,-2\}, \  0\leq i\leq n-2,\ \text{and}\ y_{n}-y_{n-1}\in \{1+\epsilon, -1-\epsilon\}\}.
\end{split}
\]

For all $\ell \equiv 0 \pmod 4$,
\[
\begin{split}
\mathscr{P}_{n,\ell}=& \{((4n-2,y_{0}),(4n-4,y_{1}),\ldots,(2n+2,y_{n-2}), (2n,y_{n-1}),(2n-1,y_{n})): y_{0}=\ell+2n-1,\\
& \ y_{i+1}-y_{i}\in \{2,-2\}, \   0\leq i\leq n-2, \ \text{and}\  y_{n}-y_{n-1}\in \{1+\epsilon, -1-\epsilon\}\}.
\end{split}
\]

For all $(i,k)\in \mathcal{X}$, $i< n$, let
\begin{align*}
\mathscr{P}_{i,k}=&\{(a_{0},a_{1},\ldots,a_{n},\overline{a}_{n},\ldots,\overline{a}_{1},\overline{a}_{0}):  (a_{0},a_{1},\ldots,a_n)\in \mathscr{P}_{n,k-(2n-2i-1)},\\
&(\overline{a}_{0},\overline{a}_{1},\ldots,\overline{a}_{n})\in \mathscr{P}_{n,k+(2n-2i-1)}, \text { and } a_{n}-\overline{a}_{n}= (0,y) \text{ where } y>0\}.
\end{align*}
The sets $C^{\pm}_{p}$ of upper and lower corners of a path $p=((j_{r},\ell_{r}))_{0\leq r \leq |p|-1}\in \mathscr{P}_{i,k}$,
where $|p|$ is the number of points in the path $p$, are defined as follows:
\begin{align*}
C^{+}_{p}=\ & \iota^{-1}\{(j_{r},\ell_{r})\in p: j_{r}\not \in \{0,2n-1,4n-2\}, \ \ell_{r-1}> \ell_{r},\ \ell_{r+1}>\ell_{r}\} \\
\ & \sqcup \{(n,\ell)\in \mathcal{X}: (2n-1,\ell-\epsilon)\in p \text{ and }(2n-1,\ell+\epsilon)\not\in p\},
\end{align*}
\begin{align*}
C^{-}_{p}=\ & \iota^{-1}\{(j_{r},\ell_{r})\in p: j_{r}\not \in \{0,2n-1,4n-2\}, \ \ell_{r-1}< \ell_{r},\ \ell_{r+1}<\ell_{r}\} \\
\ & \sqcup \{(n,\ell)\in \mathcal{X}: (2n-1,\ell-\epsilon)\not\in p \text{ and } (2n-1,\ell+\epsilon)\in p\}.
\end{align*}

Each $\mathscr{P}_{(i,k)}$, $i\in I$, $k\in \mathbb{Z}$, corresponds to a rectangle ($i\neq n$) or a triangle ($i=n$), see Figure \ref{(4,2)}--Figure \ref{(2,1)}. The monomials appearing in $\chi_q(L(i_k))$ correspond to paths in the rectangle or the triangle, see Section $5$ of \cite{MY12a} and Section $6$ of \cite{MY12b} for details.
\begin{figure}[H]
\resizebox{.5\width}{.5\height}{
\begin{minipage}[b]{0.7\linewidth}
\centerline{
\begin{tikzpicture}
\draw[step=.5cm,gray,thin] (-1,0) grid (6,8);
\draw[fill] (0,8.3) circle (2pt)--(1,8.3) circle (2pt);
\draw[fill] (1,8.3) circle (2pt)--(2,8.3) circle (2pt);
\draw[fill] (3,8.3) circle (2pt)--(4,8.3) circle (2pt);
\draw[fill] (4,8.3) circle (2pt)--(5,8.3) circle (2pt);
\draw[fill] (2.5,8.3) circle (2pt);
\node [above] at (0,8.3) {1};
\node [above] at (1,8.3) {2};
\node [above] at (2,8.3) {3};
\node [above] at (2.5,8.3) {4};
\node [above] at (3,8.3) {3};
\node [above] at (4,8.3) {2};
\node [above] at (5,8.3) {1};
\draw [double,->] (2.075,8.3)--(2.425,8.3);
\draw [double,->] (2.925,8.3)--(2.575,8.3);
\node [left] at (-1,8) {0};
\node [left] at (-1,7) {2};
\node [left] at (-1,6) {4};
\node [left] at (-1,5) {6};
\node [left] at (-1,4) {8};
\node [left] at (-1,3) {10};
\node [left] at (-1,2) {12};
\node [left] at (-1,1) {14};
\node [left] at (-1,0) {16};
\draw[thick] (-1,3.5)--(0.5,2)--(1.5,1)--(2,0.5)--(2.5,-0.1);
\draw[thick] (-1,3.5)--(0,4.5)--(1,5.5)--(2,6.5)--(2.5,7.1);
\draw[thick] (2.5,-0.1)--(2.5,7.1);
\end{tikzpicture}}
\end{minipage}
\begin{minipage}[b]{0.7\linewidth}
\centerline{
\begin{tikzpicture}
\draw[step=.5cm,gray,thin] (-1,0) grid (6,8);
\draw[fill] (0,8.3) circle (2pt)--(1,8.3) circle (2pt);
\draw[fill] (1,8.3) circle (2pt)--(2,8.3) circle (2pt);
\draw[fill] (3,8.3) circle (2pt)--(4,8.3) circle (2pt);
\draw[fill] (4,8.3) circle (2pt)--(5,8.3) circle (2pt);
\draw[fill] (2.5,8.3) circle (2pt);
\node [above] at (0,8.3) {1};
\node [above] at (1,8.3) {2};
\node [above] at (2,8.3) {3};
\node [above] at (2.5,8.3) {4};
\node [above] at (3,8.3) {3};
\node [above] at (4,8.3) {2};
\node [above] at (5,8.3) {1};
\draw [double,->] (2.075,8.3)--(2.425,8.3);
\draw [double,->] (2.925,8.3)--(2.575,8.3);
\node [left] at (-1,8) {0};
\node [left] at (-1,7) {2};
\node [left] at (-1,6) {4};
\node [left] at (-1,5) {6};
\node [left] at (-1,4) {8};
\node [left] at (-1,3) {10};
\node [left] at (-1,2) {12};
\node [left] at (-1,1) {14};
\node [left] at (-1,0) {16};
\begin{scope}[thick, every node/.style={sloped,allow upside down}]
\draw (-1,3.5)--node {\midarrow}(0.5,2);
\draw (0,2.5)--node {\midarrow}(1.5,1);
\draw (1.5,1)--node {\midarrow}(2,0.5);
\draw (2,0.5)--node {\midarrow}(2.5,-0.1);
\draw (2,0.5)--node {\midarrow}(2.5,1.1);
\draw (-1,3.5)--node {\midarrow}(0.5,5);
\draw (0,4.5)--node {\midarrow}(1.5,6);
\draw (1.5,6)--node {\midarrow}(2,6.5);
\draw (2,6.5)--node {\midarrow}(2.5,7.1);
\draw (2,6.5)--node {\midarrow}(2.5,5.9);
\draw (0,4.5)--node {\midarrow}(1.5,3);
\draw (1.5,3)--node {\midarrow}(2,2.5);
\draw (2,2.5)--node {\midarrow}(2.5,1.9);
\draw (2,2.5)--node {\midarrow}(2.5,3.1);
\draw (1,5.5)--(1.5,5) (1.5,5)--node {\midarrow}(2,4.5);
\draw (2,4.5)--node {\midarrow}(2.5,3.9);
\draw (2,4.5)--node {\midarrow}(2.5,5.1);
\draw (0,2.5)--node {\midarrow}(1.5,4);
\draw (1.5,4)--node {\midarrow}(2,4.5);
\draw (2,4.5)--node {\midarrow}(2.5,5.1);
\draw (1,1.5)--(1.5,2) (1.5,2)--node {\midarrow}(2,2.5);
\draw (2,2.5)--node {\midarrow}(2.5,3.1);
\end{scope}
\end{tikzpicture}}
\end{minipage}}
\caption{In type $B_4$: left, $\mathscr{P}_{(4,2)}$ corresponds to a triangle; right, the paths corresponding to monomials of $\chi_q(L(4_2))$.}\label{(4,2)}
\end{figure}
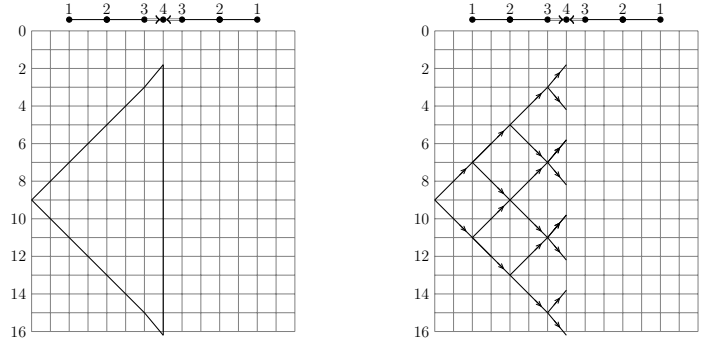

\begin{figure}[H]
\resizebox{.6\width}{.6\height}{
\begin{minipage}[b]{0.7\linewidth}
\centerline{
\begin{tikzpicture}
\draw[step=.5cm,gray,thin] (-1,0) grid (6,8);
\draw[fill] (0,8.3) circle (2pt)--(1,8.3) circle (2pt);
\draw[fill] (1,8.3) circle (2pt)--(2,8.3) circle (2pt);
\draw[fill] (3,8.3) circle (2pt)--(4,8.3) circle (2pt);
\draw[fill] (4,8.3) circle (2pt)--(5,8.3) circle (2pt);
\draw[fill] (2.5,8.3) circle (2pt);
\node [above] at (0,8.3) {1};
\node [above] at (1,8.3) {2};
\node [above] at (2,8.3) {3};
\node [above] at (2.5,8.3) {4};
\node [above] at (3,8.3) {3};
\node [above] at (4,8.3) {2};
\node [above] at (5,8.3) {1};
\draw [double,->] (2.075,8.3)--(2.425,8.3);
\draw [double,->] (2.925,8.3)--(2.575,8.3);
\node [left] at (-1,8) {0};
\node [left] at (-1,7) {2};
\node [left] at (-1,6) {4};
\node [left] at (-1,5) {6};
\node [left] at (-1,4) {8};
\node [left] at (-1,3) {10};
\node [left] at (-1,2) {12};
\node [left] at (-1,1) {14};
\node [left] at (-1,0) {16};
\draw[thick] (6,4.5)--(5,3.5)--(4,2.5)--(3,1.5)--(2.5,0.9);
\draw[thick] (6,4.5)--(5,5.5)--(4,6.5)--(3,7.5)--(2.5,8.1);
\draw[thick] (2.5,0.9)--(2.5,8.1);
\end{tikzpicture}}
\end{minipage}
\begin{minipage}[b]{0.7\linewidth}
\centerline{
\begin{tikzpicture}
\draw[step=.5cm,gray,thin] (-1,0) grid (6,8);
\draw[fill] (0,8.3) circle (2pt)--(1,8.3) circle (2pt);
\draw[fill] (1,8.3) circle (2pt)--(2,8.3) circle (2pt);
\draw[fill] (3,8.3) circle (2pt)--(4,8.3) circle (2pt);
\draw[fill] (4,8.3) circle (2pt)--(5,8.3) circle (2pt);
\draw[fill] (2.5,8.3) circle (2pt);
\node [above] at (0,8.3) {1};
\node [above] at (1,8.3) {2};
\node [above] at (2,8.3) {3};
\node [above] at (2.5,8.3) {4};
\node [above] at (3,8.3) {3};
\node [above] at (4,8.3) {2};
\node [above] at (5,8.3) {1};
\draw [double,->] (2.075,8.3)--(2.425,8.3);
\draw [double,->] (2.925,8.3)--(2.575,8.3);
\node [left] at (-1,8) {0};
\node [left] at (-1,7) {2};
\node [left] at (-1,6) {4};
\node [left] at (-1,5) {6};
\node [left] at (-1,4) {8};
\node [left] at (-1,3) {10};
\node [left] at (-1,2) {12};
\node [left] at (-1,1) {14};
\node [left] at (-1,0) {16};
\begin{scope}[thick, every node/.style={sloped,allow upside down}]
\draw (6,4.5)--node {\midarrow}(4.5,3);
\draw (5,3.5)--node {\midarrow}(3.5,2);
\draw (4,2.5)--(3.5,2) (3.5,2)--node {\midarrow}(3,1.5);
\draw (3,1.5)--node {\midarrow}(2.5,0.9);
\draw (3,1.5)--node {\midarrow}(2.5,2.1);
\draw (6,4.5)--node {\midarrow}(4.5,6);
\draw (5,5.5)--node {\midarrow}(3.5,7);
\draw (4,6.5)--(3.5,7) (3.5,7)--node {\midarrow}(3,7.5);
\draw (3,7.5)--node {\midarrow}(2.5,8.1);
\draw (3,7.5)--node {\midarrow}(2.5,6.9);
\draw (5,5.5)--node {\midarrow}(3.5,4);
\draw (4,4.5)--(3,3.5) (3,3.5)--node {\midarrow}(2.5,2.9);
\draw (3,3.5)--node {\midarrow}(2.5,4.1);
\draw (5,3.5)--node {\midarrow}(3.5,2);
\draw (4,2.5)--(3.5,3) (3.5,3)--node {\midarrow}(3,3.5);
\draw (5,3.5)--node {\midarrow}(3.5,5) (3.5,5)--node {\midarrow}(3,5.5);
\draw (4,6.5)--(3,5.5) (3,5.5)--node {\midarrow}(2.5,4.9);
\draw (3,5.5)--node {\midarrow}(2.5,6.1);
\end{scope}
\end{tikzpicture}}
\end{minipage}}
\caption{In type $B_4$: left, $\mathscr{P}_{(4,0)}$ corresponds to a triangle; right, the paths corresponding to monomials of $\chi_q(L(4_0))$.}\label{(4,0)}
\end{figure}
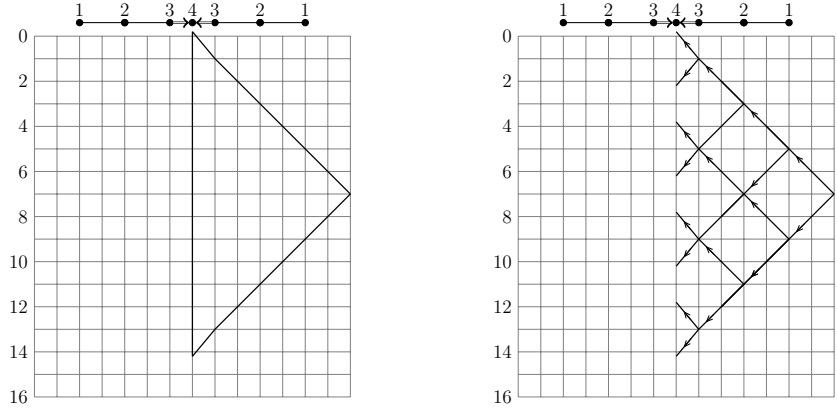

\begin{figure}[H]
\resizebox{.6\width}{.6\height}{
\begin{minipage}[b]{0.7\linewidth}
\centerline{
\begin{tikzpicture}
\draw[step=.5cm,gray,thin] (-1,0) grid (6,8);
\draw[fill] (0,8.3) circle (2pt)--(1,8.3) circle (2pt);
\draw[fill] (1,8.3) circle (2pt)--(2,8.3) circle (2pt);
\draw[fill] (3,8.3) circle (2pt)--(4,8.3) circle (2pt);
\draw[fill] (4,8.3) circle (2pt)--(5,8.3) circle (2pt);
\draw[fill] (2.5,8.3) circle (2pt);
\node [above] at (0,8.3) {1};
\node [above] at (1,8.3) {2};
\node [above] at (2,8.3) {3};
\node [above] at (2.5,8.3) {4};
\node [above] at (3,8.3) {3};
\node [above] at (4,8.3) {2};
\node [above] at (5,8.3) {1};
\draw [double,->] (2.075,8.3)--(2.425,8.3);
\draw [double,->] (2.925,8.3)--(2.575,8.3);
\node [left] at (-1,8) {0};
\node [left] at (-1,7) {2};
\node [left] at (-1,6) {4};
\node [left] at (-1,5) {6};
\node [left] at (-1,4) {8};
\node [left] at (-1,3) {10};
\node [left] at (-1,2) {12};
\node [left] at (-1,1) {14};
\node [left] at (-1,0) {16};
\draw[thick] (-1,5.5)--(1,7.5)--(2,6.5)--(2.5,5.9)--(2.5,6.1)--(6,2.5);
\draw[thick] (-1,5.5)--(2,2.5)--(2.5,1.9)--(2.5,2.1)--(4,0.5)--(6,2.5);
\end{tikzpicture}}
\end{minipage}
\begin{minipage}[b]{0.7\linewidth}
\centerline{
\begin{tikzpicture}
\draw[step=.5cm,gray,thin] (-1,0) grid (6,8);
\draw[fill] (0,8.3) circle (2pt)--(1,8.3) circle (2pt);
\draw[fill] (1,8.3) circle (2pt)--(2,8.3) circle (2pt);
\draw[fill] (3,8.3) circle (2pt)--(4,8.3) circle (2pt);
\draw[fill] (4,8.3) circle (2pt)--(5,8.3) circle (2pt);
\draw[fill] (2.5,8.3) circle (2pt);
\node [above] at (0,8.3) {1};
\node [above] at (1,8.3) {2};
\node [above] at (2,8.3) {3};
\node [above] at (2.5,8.3) {4};
\node [above] at (3,8.3) {3};
\node [above] at (4,8.3) {2};
\node [above] at (5,8.3) {1};
\draw [double,->] (2.075,8.3)--(2.425,8.3);
\draw [double,->] (2.925,8.3)--(2.575,8.3);
\node [left] at (-1,8) {0};
\node [left] at (-1,7) {2};
\node [left] at (-1,6) {4};
\node [left] at (-1,5) {6};
\node [left] at (-1,4) {8};
\node [left] at (-1,3) {10};
\node [left] at (-1,2) {12};
\node [left] at (-1,1) {14};
\node [left] at (-1,0) {16};
\draw[thick] (2.5,2.1)--(2.5,1.9) (2.5,6.1)--(2.5,5.9) (2.5,2.9)--(2.5,3.1) (1,3.5)--(2,2.5) (1,5.5)--(2,4.5) (2.5,3.9)--(2.5,4.1);
\begin{scope}[thick, every node/.style={sloped,allow upside down}]
\draw (1,7.5)--node {\midarrow}(2,6.5);
\draw (2,6.5)--node {\midarrow}(2.5,5.9);
\draw (2.5,6.1)--node {\midarrow}(4,4.5);
\draw (3.5,5)--node {\midarrow}(5,3.5);
\draw (4.5,4)--node {\midarrow}(6,2.5);
\draw (-1,5.5)--node {\midarrow}(0.5,7);
\draw (0.5,7)--node {\midarrow}(1,7.5);
\draw (-1,5.5)--node {\midarrow}(0.5,4);
\draw (0,4.5)--node {\midarrow}(1.5,3);
\draw (2,2.5)--node {\midarrow}(2.5,1.9);
\draw (2,2.5)--node {\midarrow}(2.5,3.1);
\draw (2.5,2.9)--node {\midarrow}(3,3.5);
\draw (3,3.5)--node {\midarrow}(4,4.5);
\draw (2.5,2.1)--node {\midarrow}(4,0.5);
\draw (4,0.5)--node {\midarrow}(5.5,2);
\draw (5,1.5)--node {\midarrow}(6,2.5);
\draw (0,6.5)--node {\midarrow}(1.5,5);
\draw (2,4.5)--node {\midarrow}(2.5,3.9);
\draw (2,4.5)--node {\midarrow}(2.5,5.1);
\draw (2.5,4.1)--node {\midarrow}(4,2.5);
\draw (3.5,3)--node {\midarrow}(5,1.5);
\draw (2.5,4.9)--node {\midarrow}(3,5.5);
\draw (3,1.5)--node {\midarrow}(4,2.5);
\draw (4,2.5)--node {\midarrow}(5,3.5);
\draw (0,4.5)--node {\midarrow} (1,5.5);
\draw (1,5.5)--node {\midarrow} (2,6.5);
\draw (1,3.5)--node {\midarrow} (2,4.5);
\draw (2.5,2)--node {\midarrow} (2.5,3);
\draw (2.5,3)--node {\midarrow} (2.5,4);
\draw (2.5,4)--node {\midarrow} (2.5,5);
\draw (2.5,5)--node {\midarrow} (2.5,6);
\end{scope}
\end{tikzpicture}}
\end{minipage}}
\caption{In type $B_4$: left, $\mathscr{P}_{(2,1)}$ corresponds to a rectangle; right, the paths corresponding to monomials of $\chi_q(L(2_1))$.}\label{(2,1)}
\end{figure}
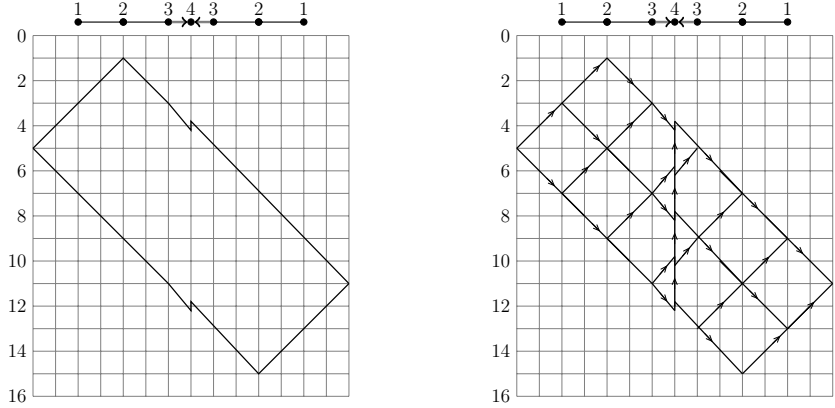

A mapping $m$ sending paths to monomials is defined by
\begin{align}
m: \bigsqcup_{(i,k)\in \mathcal{X}}{\hskip -0.5em}\mathscr{P}_{i,k} & \longrightarrow \mathbb{Z}[j^{\pm}_{\ell}]_{(j,\ell)\in \mathcal{X}} \nonumber \\
p& \longmapsto  m(p)=\prod_{(j,\ell)\in C^{+}_{p}}{\hskip -0.5em}j_\ell{\hskip -0.5em}\prod_{(j,\ell)\in C^{-}_{p}}{\hskip -0.5em}j^{-1}_\ell.
\end{align}

We always identify a path $p$ with the monomial $m(p)$.

Let $p,p'$ be paths. We say that $p$ is \textit{strictly above} $p'$ or $p'$ is \textit{strictly below} $p$ if
\begin{align*}
(x,y)\in p \text{ and } (x,z)\in p' \Longrightarrow y < z.
\end{align*}
We say that a $T$-tuple of paths $(p_{1},\ldots,p_{T})$ is \textit{non-overlapping} if $p_{s}$ is strictly above $p_{t}$ for all $s<t$.
For any snake $(i_{t},k_{t})\in \mathcal{X}$, $1\leq t\leq T$, $T\in \mathbb{Z}_{\geq1}$, let
\begin{align*}
\overline{\mathscr{P}}_{(i_{t},k_{t})_{1\leq t\leq T}}=\{(p_{1},\ldots,p_{T}): p_{t}\in \mathscr{P}_{i_{t},k_{t}}, \ 1\leq t\leq T, \ (p_{1},\ldots,p_{T})\text { is } \text {non-overlapping} \}.
\end{align*}
\begin{theorem}[{\cite[Theorem 6.1]{MY12a}; \cite[Theorem 6.5]{MY12b}}] \label{path description of q-characters}
Let $(i_{\ell},k_{\ell}) \in \mathcal{X}$, $1 \leq \ell \leq T$, be a snake of length $T \in \mathbb{Z}_{\geq 1}$. Then
\begin{align}
\chi_{q} (L(\prod_{\ell=1}^{T} (i_{\ell})_{k_\ell}))=\sum_{(p_{1},\ldots,p_{T}) \in \overline{\mathscr{P}}_{(i_{\ell},k_{\ell})_{1 \leq \ell \leq T}}} \prod_{\ell=1}^{T}m(p_{\ell}).
\end{align}
The module $L(\prod_{\ell=1}^{T} (i_{\ell})_{k_\ell})$ is thin, special and anti-special.
\end{theorem}

In view of Theorem \ref{path description of q-characters}, the $q$-characters of snake modules of types $A_n$ and $B_n$ with length $T$ are given by a set of $T$-tuples of non-overlapping paths, the path in each $T$-tuple is non-overlapping. This property is called the \textit{non-overlapping property}.

We also need the following notations in this paper. For all $(i,k)\in \mathcal{X}$, let $p^{+}_{i,k}$ be the highest path which is the unique path in $\mathscr{P}_{i,k}$ with no lower corners and $p^{-}_{i,k}$ the lowest path which is the unique path in $\mathscr{P}_{i,k}$ with no upper corners.

\section{$S$-systems of types $A_{n}$ and $B_{n}$ }\label{S-systems of types $A_{n}$, $B_{n}$}
In this section, we introduce a closed system of equations which contains just all prime snake modules of type $A_n$ (respectively, $B_n$).

\subsection{Another notation of snake modules}\label{New notation of snake modules}
In order to introduce the $S$-systems, we need to use another notation of snake modules. We fix an $a\in \mathbb{C}^{\times}$ and denote $i_s = Y_{i, aq^s}$, where $i \in I$, $s \in \mathbb{Z}$. For $i,j \in I$, let $\varepsilon_{i,j}= -\delta_{in} - \delta_{jn}$, where $\delta_{ij}$ is the Kronecker delta.

By the definitions of snake positions and snake modules, every snake module of type $A_n$ is a $U_q(\widehat{\mathfrak{g}})$-module with highest $l$-weight monomial of the form
\begin{align*} 
S^{(t)}_{k_{1}^{(i_{1},j_{1})}, k_{2}^{(i_{2},j_{2})}, \ldots, k_{m-1}^{(i_{m-1},j_{m-1})}, k_{m}^{(i_{m})}} :=
&\prod_{j=1}^{m} \left( \prod_{r=0}^{k_{j}-1} (i_{j})_{t+2r+\sum_{\ell=1}^{j-1}n_{\ell}} \right),
\end{align*}
where $t\in\mathbb{Z}$, $i_{j}\in I,\ k_{j} \geq 0,\ 1 \leq j \leq m$, $j_\ell \in \mathbb{Z}_{\geq 0}$, $1 \leq \ell \leq m-1$ and
\begin{align}\label{n_1 in type A_n}
n_{\ell} = 2k_{\ell} + |i_{\ell+1}-i_{\ell}| + 2j_{\ell}.
\end{align}

Every snake module of type $B_n$ is a $U_q(\widehat{\mathfrak{g}})$-module with highest $l$-weight monomial of the form
\begin{align*} 
S^{(t)}_{k_{1}^{(i_{1},j_{1})}, k_{2}^{(i_{2},j_{2})}, \ldots, k_{m-1}^{(i_{m-1},j_{m-1})}, k_{m}^{(i_{m})}} :=
&\prod_{j=1}^{m} \left( \prod_{r=0}^{k_{j}-1} (i_{j})_{t+2d_{i_{j}}r+\sum_{\ell=1}^{j-1}n_{\ell}} \right),
\end{align*}
where $t\in\mathbb{Z}$, $i_{j}\in I,\ k_{j} \geq 0,\ 1 \leq j \leq m$, $j_\ell \in \mathbb{Z}_{\geq 0}$, $1 \leq \ell \leq m-1$ and
\begin{align}\label{n_1 in type B_n}
n_{\ell} = 2d_{i_\ell}k_{\ell} + 2|i_{\ell+1}-i_{\ell}| +4-2d_{i_\ell}+ 4j_{\ell}+\varepsilon_{i_\ell,i_{\ell+1}}.
\end{align}
Let $S$ be a dominant monomial. We also use $\mathcal{S}$ to denote $L(S)$. In particular, we use
\[
\mathcal{S}^{(t)}_{k_{1}^{(i_{1},j_{1})}, k_{2}^{(i_{2},j_{2})}, \ldots, k_{m-1}^{(i_{m-1},j_{m-1})}, k_{m}^{(i_{m})}}
\]
to denote the finite-dimensional simple $U_{q}\widehat{\mathfrak{g}}$-module with highest $l$-weight monomial
\[
S^{(t)}_{k_{1}^{(i_{1},j_{1})}, k_{2}^{(i_{2},j_{2})}, \ldots, k_{m-1}^{(i_{m-1},j_{m-1})}, k_{m}^{(i_{m})}}.
\]

For simplicity, if $j_{\ell}=0$ for some $\ell$, $1 \leq \ell \leq m-1$, then we use
\[
S^{(t)}_{k_{1}^{(i_{1},j_{1})}, k_{2}^{(i_{2},j_{2})}, \ldots, k_{\ell}^{(i_{\ell})}, \ldots, k_{m-1}^{(i_{m-1},j_{m-1})}, k_{m}^{(i_{m})}}
\]
to denote $S^{(t)}_{k_{1}^{(i_{1},j_{1})}, k_{2}^{(i_{2},j_{2})}, \ldots, k_{\ell}^{(i_{\ell},0)}, \ldots, k_{m-1}^{(i_{m-1},j_{m-1})}, k_{m}^{(i_{m})}}$. In this notation, $\scalemath{0.92}{\mathcal{S}^{(t)}_{k_{1}^{(i_{1})}, k_{2}^{(i_{2})}, \ldots, k_{m-1}^{(i_{m-1})}, k_{m}^{(i_{m})}}}$ is a minimal snake module.

Let $\mathfrak{S}$ be the set of all prime snake modules and $\mathfrak{T}$ the set of all prime snakes.
We define a mapping
\begin{equation}\label{corresponding map}
\begin{split}
\varphi:  \mathfrak{S} & \to \mathfrak{T} \\
\mathcal{S} & \mapsto \text{ the snake of } \mathcal{S}.
\end{split}
\end{equation}
It is easy to see that the mapping $\varphi:  \mathfrak{S} \to \mathfrak{T}$ is a bijection.

\subsection{Neighboring points}\label{Section neighboring points}
The concept of neighboring points was introduced in Section $3$ of \cite{MY12b}. Let $(i,k)\in \mathcal{X}$ and $(i',k')\in \mathcal{X}$ such that
$(i',k')$ is in prime snake position with respect to $(i,k)$. The \textit{neighboring points} to the pair $(i,k)$, $(i',k')$ are two finite sequences $\mathbb{X}_{i,k}^{i',k'}$ and $\mathbb{Y}_{i,k}^{i',k'}$ of points in $\mathcal{X}$ defined as follows.

\begin{figure}[H]
\centerline{
\begin{tikzpicture}
\draw[fill] (1,1)--(2,2)circle (2pt)--(3,3)--(4,4) circle (2pt)--(5,3)circle (2pt)--(6,2);
\draw[fill] (1,3)--(2,2)--(3,1)circle (2pt)--(4,2)--(5,3)--(6,4);
\node [above] at (4,4) {$\iota(i,k)$};
\node [below] at (3,1) {$\iota(i',k')$};
\node [left]  at (1.9,1.9) {$\iota(\mathbb{X}_{i,k}^{i',k'})$};
\node [right] at (5.1,3) {$\iota(\mathbb{Y}_{i,k}^{i',k'})$};
\end{tikzpicture}}
\caption{The sequences $\mathbb{X}_{i,k}^{i',k'}$ and $\mathbb{Y}_{i,k}^{i',k'}$ are neighboring points to the pair $(i,k)$, $(i',k')$.}
\end{figure}
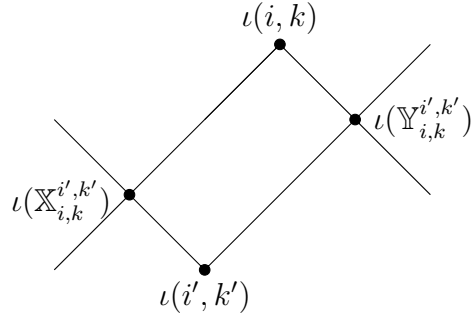

In the case of type $A_n$,
\begin{align*}
\mathbb{X}_{i,k}^{i',k'}=
\begin{cases}
((\frac{1}{2}(i+k+i'-k'),\frac{1}{2}(i+k-i'+k'))), & k+i > k'-i',\\
\emptyset, & k+i = k'-i', {\hskip 1in}
\end{cases}
\end{align*}
\begin{align*}
\mathbb{Y}_{i,k}^{i',k'}=
\begin{cases}
((\frac{1}{2}(i'+k'+i-k),\frac{1}{2}(i'+k'-i+k))), & k+n+1-i > k'-n-1+i',\\
\emptyset, & k+n+1-i = k'-n-1+i'.
\end{cases}
\end{align*}

In the case of type $B_n$,
\begin{gather}
\begin{align*}
(\mathbb{X}_{i,k}^{i',k'},\mathbb{Y}_{i,k}^{i',k'})=
\begin{cases}
(B_{i,k}^{i',k'},F_{i,k}^{i',k'}), & {\hskip -0.2em} i<n,\ 2n+k-2i \equiv 1 {\hskip -0.7em}\pmod 4, \text{ or } i=n,\ k \equiv 0 {\hskip -0.7em}\pmod 4,\\
(F_{i,k}^{i',k'},B_{i,k}^{i',k'}), & {\hskip -0.2em} i<n,\ 2n+k-2i \equiv 3 {\hskip -0.7em}\pmod 4, \text{ or } i=n,\ k \equiv 2 {\hskip -0.7em}\pmod 4,
\end{cases}
\end{align*}
\end{gather}
where
\begin{gather}
\begin{align*}
&B_{i,k}^{i',k'}=
\begin{cases}
\emptyset, & i<n,\ i'<n,\ k'-k=2i+2i',\\
((\frac{1}{4}(2i+k+2i'-k'),\frac{1}{2}(2i+k-2i'+k'))), & i<n,\ i'<n,\ k'-k<2i+2i',\\
\emptyset, & i<n,\ i'=n,\ k'-k=2i+2n-1,\\
((\frac{1}{4}(2i+k+2n-1-k'),\frac{1}{2}(2i+k-2n+1+k'))), & i<n,\ i'=n,\ k'-k<2i+2n-1,\\
((n,k'-2n+1+2i')), & i=n,\ i'<n,\\
\emptyset, & i=n,\ i'=n,
\end{cases}\\
&F_{i,k}^{i',k'}=
\begin{cases}
((\frac{1}{4}(2i'+k'+2i-k),\frac{1}{2}(2i'+k'-2i+k))), & i<n,\ i'<n,\ k'-k\ \leq 4n-4-2i-2i',\\
((n,k+2n-1-2i),(n,k'-2n+1+2i')), & i<n,\ i'<n,\ k'-k\ \geq 4n-2i-2i',\\
((n,k+2n-1-2i)), & i<n,\ i'=n,\\
((\frac{1}{4}(2n-1+k+2i'-k'),\frac{1}{2}(2n-1+k-2i'+k'))), & i=n,\ i'<n,\ k'-k<2n+2i'-1,\\
\emptyset, & i=n,\ i'<n,\ k'-k=2n+2i'-1,\\
((\frac{1}{4}(4n-2+k-k'), \frac{1}{2}(k+k'))), & i=n,\ i'=n.
\end{cases}
\end{align*}
\end{gather}

\subsection{$S$-systems of types $A_n$ and $B_n$}\label{definition of S-system}
For types $A_n$ and $B_n$, every prime snake module can be written as
\begin{align}\label{S_2}
\mathcal{S}^{(t)}_{k_{1}^{(i_{1},j_{1})}, k_{2}^{(i_{2},j_{2})}, \ldots,k_{m-1}^{(i_{m-1},j_{m-1})},k_{m}^{(i_{m})}},
\end{align}
where $m \geq 1$, $j_\ell \geq 0$, $1 \leq \ell \leq m-1$, if $j_{\ell}=0$, then $i_{\ell} \neq i_{\ell+1}$, $k_1, k_2, \ldots, k_m \in \mathbb{Z}_{\geq 1}$, $t \in \mathbb{Z}$.

Let $\mathcal{S}_2$ be the prime snake module (\ref{S_2}) and $\sgn(x)$ the sign function. We define $\mathcal{S}_1$ in Table \ref{definition S_1 in type A_{n}} (respectively, Table \ref{definition S_1 in type B_{n}}) for type $A_n$ (respectively, $B_n$).
\begin{table}[H]
\resizebox{1.0\width}{1.0\height}{
\begin{tabular}{|c|c|c|c|}
\hline %
\multicolumn{2}{|c|}{Conditions} & $\mathcal{S}_1$  & \\
\hline %
\multicolumn{2}{|c|}{ $m=1$ } & $\substack{\\ \mathcal{S}^{(t+2)}_{k_{1}^{(i_{1})}}}$ & (1)\\
\hline %
& $m=2$, $i_2 \equiv 1 {\hskip -0.2em}\pmod2$ & $\substack{\\
\mathcal{S}^{(t+2)}_{k_{1}^{(i_{1})}, {k_{2}}^{(i_{2}+\sgn(i_1-i_2))}}}$ & (2)\\
\cline{2-4}%
& $m=2$, $i_2 \equiv 0 {\hskip -0.2em}\pmod2$ & $\substack{\\
\mathcal{S}^{(t+2)}_{k_{1}^{(i_{1})}, {(k_{2}-1)}^{(i_{2}+\sgn(i_1-i_2))}}}$ & (3) \\
\cline{2-4}%
$j_1=0$ &$\substack{\\ j_2=0,\ m \geq 3,\ (i_2-i_1)(i_3-i_2)>0,\\ \underset{}{} j_\ell \geq 0,\ 3 \leq \ell \leq m-1 }$ &
$\substack{\\ \mathcal{S}^{(t+2)}_{k_{1}^{(i_{1})}, (k_{2}-1)^{(i_{2})}, k_{3}^{(i_{3},j_{3})}, k_{4}^{(i_{4},j_{4})}, \ldots, k_{m}^{(i_{m})}}}$ & (4)\\
\cline{2-4}%
& $\substack{\\ j_2 \geq 1,\ m\geq 3,\ (i_2-i_1)(i_3-i_2)\geq 0,\\ \underset{}{} j_\ell \geq 0,\ 3 \leq \ell \leq m-1}$
&$\substack{ \mathcal{S}^{(t+2)}_{k_{1}^{(i_{1})}, {k_{2}}^{(i_{2}+\sgn(i_1-i_2),j_2-1)}, k_{3}^{(i_{3},j_{3})}, k_{4}^{(i_{4},j_{4})}, \ldots, k_{m}^{(i_{m})}}}$ & (5) \\
\cline{2-4}%
& $\substack{ \\ j_2 \geq 0,\ m\geq 3,\ (i_2-i_1)(i_3-i_2)<0,\\ \underset{}{} j_\ell \geq 0,\ 3 \leq \ell \leq m-1}$
&$\substack{ \mathcal{S}^{(t+2)}_{k_{1}^{(i_{1})}, {k_{2}}^{(i_{2}+\sgn(i_1-i_2),j_2)}, k_{3}^{(i_{3},j_{3})}, k_{4}^{(i_{4},j_{4})}, \ldots, k_{m}^{(i_{m})}}}$ & (6) \\
\hline%
$ j_1 \geq 1$ &$\substack{\\ m\geq 2,\ j_\ell \geq 0,\ 2 \leq \ell \leq m-1 \\}$ & $\substack{\\ \mathcal{S}^{(t+2)}_{k_{1}^{(i_{1},j_{1}-1)}, k_{2}^{(i_{2},j_{2})},
k_{3}^{(i_{3},j_{3})}, k_{4}^{(i_{4},j_{4})}, \ldots,  k_{m}^{(i_{m})}}}$ & (7) \\
\hline %
\end{tabular}}
\caption{Definition of $\mathcal{S}_1$ in type $A_{n}$.}\label{definition S_1 in type A_{n}}
\end{table}
\begin{table}[H]
\resizebox{0.8\width}{0.8\height}{
\begin{tabular}{|c|c|c|c|}
\hline %
\multicolumn{2}{|c|}{Conditions} & $\mathcal{S}_1$  & \\
\hline %
\multicolumn{2}{|c|}{$m=1$} & $\substack{\\ \mathcal{S}^{(t+2d_{i_1})}_{k_{1}^{(i_{1})} }}$ & (1)\\
\hline %
\multicolumn{2}{|c|}{$j_1=0,\ i_1 \neq n,\ i_2=n,\ m=2$,\ \text{$k_2$ is odd}} & $\substack{\\ \mathcal{S}^{(t+4)}_{k_{1}^{(i_{1})}, (\frac{k_{2}-1}{2})^{(n-1)}}}$ & (2)\\
\hline %
\multicolumn{2}{|c|}{$j_1=0,\ i_1 \neq n,\ i_2=n,\ m=2$, \ \text{$k_2$ is even}} & $\substack{\\ \mathcal{S}^{(t+4)}_{{k_{1}}^{(i_{1})}, (\frac{k_{2}-2}{2})^{(n-1)}}}$ & (3)\\
\hline%
& $m=2$ & $\substack{\\ \mathcal{S}^{(t+4)}_{k_{1}^{(i_{1})}, {k_{2}}^{(i_{2}+\sgn(i_1-i_2))}}}$ & (4) \\
\cline{2-4}%
$j_1=0,\ i_1\neq n,$ & $\substack{\\ j_2=0,\ m\geq 3,\ (i_2-i_1)(i_3-i_2)>0,\\ \underset{}{} j_\ell \geq 0,\ 3\leq \ell \leq m-1 }$
&$\substack{\\ \mathcal{S}^{(t+4)}_{k_{1}^{(i_{1})},(k_{2}-1)^{(i_{2})}, k_{3}^{(i_{3},j_{3})}, k_{4}^{(i_{4},j_{4})}, \ldots, k_{m}^{(i_{m})}}}$ & (5)\\
\cline{2-4}%
$i_2 \neq n$ & $\substack{\\ j_2 \geq 1,\ m\geq 3,\ (i_2-i_1)(i_3-i_2)\geq 0,\\ \underset{}{} j_\ell \geq 0,\ 3 \leq \ell \leq m-1}$
&$\substack{\\ \mathcal{S}^{(t+4)}_{k_{1}^{(i_{1})}, {k_{2}}^{(i_{2}+\sgn(i_1-i_2),j_2-1)}, k_{3}^{(i_{3},j_{3})}, k_{4}^{(i_{4},j_{4})}, \ldots, k_{m}^{(i_{m})}}}$  & (6)\\
\cline{2-4}%
&$\substack{\\ j_2 \geq 0,\ m\geq 3,\ (i_2-i_1)(i_3-i_2)<0,\\ \underset{}{} j_\ell \geq 0,\ 3 \leq \ell \leq m-1}$
&$\substack{\\ \mathcal{S}^{(t+4)}_{k_{1}^{(i_{1})}, {k_{2}}^{(i_{2}+\sgn(i_1-i_2),j_2)}, k_{3}^{(i_{3},j_{3})}, k_{4}^{(i_{4},j_{4})}, \ldots, k_{m}^{(i_{m})}}}$  & (7)\\
\hline %
$j_1 \geq 1,\ i_1 \neq n$ & $\substack{ \\ m \geq 2,\ j_\ell \geq 0,\ 2\leq \ell \leq m-1 }$ &
$\substack{\\ \mathcal{S}^{(t+4)}_{k_{1}^{(i_{1},j_{1}-1)},k_{2}^{(i_{2},j_{2})}, k_{3}^{(i_{3},j_{3})}, k_{4}^{(i_{4},j_{4})}, \ldots,  k_{m}^{(i_{m})}}}$  & (8) \\
\hline%
& $m=2$ & $\substack{\\ \mathcal{S}^{(t+2)}_{(k_{1}+1)^{(n)}, (\frac{d_{i_{2}}}{d_{i_{2}+1}}k_{2})^{(i_{2}+1)}}}$  & (9) \\
\cline{2-4}%
& $\substack{ \\ j_2=0,\ m \geq 3,\ i_2 > i_3, \\ \underset{}{} j_\ell \geq 0,\ 3\leq \ell \leq m-1 }$
&$\substack{\\ \mathcal{S}^{(t+2)}_{(k_1+1)^{(n)}, (k_{2}-1)^{(i_{2})}, k_{3}^{(i_{3},j_{3})}, k_{4}^{(i_{4},j_{4})}, \ldots, k_{m}^{(i_{m})}}}$ & (10) \\
\cline{2-4}%
$j_1=0,\ i_1=n$  & $\substack{\\ j_2 \geq 1,\ m\geq 3,\ i_2 \geq i_3,\\ j_\ell \geq 0,\ 3 \leq \ell \leq m-1  }$ &
$\substack{\\ \mathcal{S}^{(t+2)}_{(k_{1}+1)^{(n)}, (\frac{d_{i_{2}}}{d_{i_{2}+1}}k_{2})^{(i_{2}+1,j_{2}-1)}, k_{3}^{(i_{3},j_{3})}, k_{4}^{(i_{4},j_{4})}, \ldots, k_{m}^{(i_{m})}}}$  & (11) \\
\cline{2-4}%
& $\substack{\\ j_2 \geq 0,\ m\geq 3,\ i_2<i_3,\\ j_\ell \geq 0,\ 3 \leq \ell \leq m-1  }$ &
$\substack{\\ \mathcal{S}^{(t+2)}_{(k_{1}+1)^{(n)}, (\frac{d_{i_{2}}}{d_{i_{2}+1}}k_{2})^{(i_{2}+1,j_{2})}, k_{3}^{(i_{3},j_{3})}, k_{4}^{(i_{4},j_{4})}, \ldots, k_{m}^{(i_{m})}}}$  & (12) \\
\hline%
&$\substack{\\ j_2=0,\ k_2 \text{ is odd},\ m =3}$ & $\substack{\\ \mathcal{S}^{(t+4)}_{k_{1}^{(i_{1})}, (\frac{k_{2}-1}{2})^{(n-1)}, 1^{(n)}, (\frac{d_{i_{3}}}{d_{i_{3}+1}}k_{3})^{(i_{3}+1)}}}$  & (13) \\
\cline{2-4}%
&$\substack{\\ j_2=0,\ k_2 \text{ is odd},\ m\geq 4, \\ \underset{}{} i_3 \leq i_4,\ j_\ell \geq 0,\ 3\leq \ell \leq m-1 }$
&$\substack{\\ \mathcal{S}^{(t+4)}_{k_{1}^{(i_{1})}, (\frac{k_{2}-1}{2})^{(n-1)}, 1^{(n)}, (\frac{d_{i_{3}}}{d_{i_{3}+1}}k_{3})^{(i_{3}+1,j_{3}-\delta_{i_3i_4})}, k_{4}^{(i_{4},j_{4})},\ldots, k_{m}^{(i_{m})}}}$  & (14) \\
\cline{2-4}%
$ j_1=0,\ i_2= n$ &$\substack{ \\ j_2=0,\ k_2 \text{ is odd},\ m\geq 4,\  \\ \underset{}{} i_3 > i_4,\ j_\ell \geq 0,\ 3 \leq \ell \leq m-1 }$
&$\substack{\\ \mathcal{S}^{(t+4)}_{k_{1}^{(i_{1})}, (\frac{k_{2}-1}{2})^{(n-1)}, 1^{(n)}, (k_{3}-1)^{(i_{3},j_{3})}, k_{4}^{(i_{4},j_{4})}, \ldots, k_{m}^{(i_{m})}}}$  & (15) \\
\cline{2-4}%
&$\substack{ \\ j_2 \geq 1,\ m\geq 3,\ k_2 \text{ is odd}, \\ \underset{}{} j_\ell \geq 0,\ 3\leq \ell \leq m-1 }$
& $\substack{\\ \mathcal{S}^{(t+4)}_{k_{1}^{(i_{1})}, (\frac{k_{2}-1}{2})^{(n-1)}, 1^{(n,j_{2}-1)}, k_{3}^{(i_{3},j_{3})}, k_{4}^{(i_{4},j_{4})}, \ldots, k_{m}^{(i_{m})}}}$  & (16) \\
\cline{2-4}%
&$\substack{ \\ j_2 \geq 0,\ m\geq 3,\ k_2 \text{ is even}, \\ \underset{}{} j_\ell \geq 0,\ 3\leq \ell \leq m-1 }$
& $\substack{\\ \mathcal{S}^{(t+4)}_{k_{1}^{(i_{1})}, (\frac{k_{2}}{2})^{(n-1,j_{2}-\delta_{ni_3})}, k_{3}^{(i_{3},j_{3})}, k_{4}^{(i_{4},j_{4})}, \ldots, k_{m}^{(i_{m})}}}$  & (17) \\
\hline %
$\ j_1 \geq 1,\ i_1=n$ & $\substack{ \\ m \geq 2,\ j_\ell \geq 0,\ \underset{}{} \underset{}{} 2 \leq \ell \leq m-1 }$
& $\substack{\\ \mathcal{S}^{(t+2)}_{(k_{1}+1)^{(n,j_{1}-1)}, k_{2}^{(i_{2},j_{2})}, k_{3}^{(i_{3},j_{3})}, k_{4}^{(i_{4},j_{4})}, \ldots, k_{m}^{(i_{m})}}}$ & (18) \\
\hline%
\end{tabular}}
\caption{Definition of $\mathcal{S}_1$ in type $B_{n}$.}\label{definition S_1 in type B_{n}}
\end{table}
Let  $\mathcal{X}_1 = \varphi(\mathcal{S}_1)$ and $\mathcal{X}_2 = \varphi(\mathcal{S}_2)$, where $\varphi$ is the mapping defined in (\ref{corresponding map}). We define
\begin{align}\label{S_3 S_4}
\mathcal{S}_3 = L((i_1)_{t}{\hskip -0.4em}\prod_{(i,k)\in \mathcal{X}_1}{\hskip -0.4em}i_k), \quad
\mathcal{S}_4 = L(\prod_{(i,k)\in \mathcal{X}_2 \backslash (i_1,t)}{\hskip -1em}i_k).
\end{align}
Let
\begin{align*}
&\mathcal{X}'_{i_1}=\{(i_1,t+2d_{i_1}j): 1 \leq j \leq k_1 \} \subset \mathcal{X}_1,\\
&\mathcal{X}_{i_1}=\{(i_1,t+2d_{i_1}j-2d_{i_1}): 1 \leq j \leq k_1 \} \subset \mathcal{X}_2,
\end{align*}
and
\begin{align*}
\mathbb{X}=\prod_{(i,k)\in \mathcal{X}_{i_1}}{\hskip -0.6em} \mathbb{X}_{i,k}^{i,k+2d_{i}}, \quad \mathbb{Y}=\prod_{(i,k)\in \mathcal{X}_{i_1}} {\hskip -0.6em}\mathbb{Y}_{i,k}^{i,k+2d_{i}}.
\end{align*}
If $m=1$, let
\begin{equation}\label{S_5 S_6}
\mathcal{S}_5 = L(\prod_{(i,k)\in \mathbb{X}}{\hskip -0.4em}i_k), \quad \mathcal{S}_6 = L(\prod_{(i,k)\in \mathbb{Y}}{\hskip -0.4em}i_k).
\end{equation}
If $m \geq 2$, we define $\mathcal{S}_5$, $\mathcal{S}_6$ as follows. In the case of type $A_n$, let
\begin{equation}\label{S_5 S_6 in type $A_n$}
\begin{split}
&\mathcal{S}_5 =
\begin{cases}
L((\prod_{(i,k)\in \mathbb{X}}i_k)(\prod_{(i,k)\in \mathcal{X}_1 \backslash \mathcal{X}'_{i_1}} i_k)),& {\hskip -0.4em} i_1 \leq i_2,\\
L((\prod_{(i,k)\in \mathbb{Y}}i_k)(\prod_{(i,k)\in \mathcal{X}_1 \backslash \mathcal{X}'_{i_1}} i_k)),& {\hskip -0.4em}i_1 > i_2,
\end{cases} \\
&\mathcal{S}_6 =
\begin{cases}
L((\prod_{(i,k)\in \mathbb{Y}}i_k)(\prod_{(i,k)\in \mathcal{X}_2 \backslash \mathcal{X}_{i_1}} i_k)), & {\hskip -0.4em}i_1 \leq i_2,\\
L((\prod_{(i,k)\in \mathbb{X}}i_k)(\prod_{(i,k)\in \mathcal{X}_2 \backslash \mathcal{X}_{i_1}} i_k)), & {\hskip -0.4em}i_1 > i_2.
\end{cases}
\end{split}
\end{equation}
In the case of type $B_n$, let
\begin{equation}\label{S_5 S_6 in type $B_n$}
\begin{split}
&\mathcal{S}_5 =
\begin{cases}
L((\prod_{(i,k)\in \mathbb{X}}i_k)(\prod_{(i,k)\in \mathcal{X}_1 \backslash \mathcal{X}'_{i_1}} i_k)),&{\hskip -0.5em} \pr_1(\iota(i_1,t)) < \pr_1(\iota(i_2,t+n_1)), \text{ or } i_1=i_2 \neq n, \text{ or } i_1=i_2=n,\ t+n_1\equiv 2{\hskip -0.5em}\pmod 4, \\
L((\prod_{(i,k)\in \mathbb{Y}}i_k)(\prod_{(i,k)\in \mathcal{X}_1 \backslash \mathcal{X}'_{i_1}} i_k)),&{\hskip -0.5em} \pr_1(\iota(i_1,t)) > \pr_1(\iota(i_2,t+n_1)), \text{ or } i_1=i_2=n,\ t+n_1\equiv 0{\hskip -0.5em}\pmod 4,
\end{cases} \\
&\mathcal{S}_6 =
\begin{cases}
L((\prod_{(i,k)\in \mathbb{Y}}i_k)(\prod_{(i,k)\in \mathcal{X}_2 \backslash \mathcal{X}_{i_1}} i_k)), & {\hskip -0.5em} \pr_1(\iota(i_1,t)) < \pr_1(\iota(i_2,t+n_1)), \text{ or } i_1=i_2 \neq n,\ \text{ or } i_1=i_2=n,\ t+n_1\equiv 2{\hskip -0.5em}\pmod 4, \\
L((\prod_{(i,k)\in \mathbb{X}}i_k)(\prod_{(i,k)\in \mathcal{X}_2 \backslash \mathcal{X}_{i_1}} i_k)), & {\hskip -0.5em} \pr_1(\iota(i_1,t)) > \pr_1(\iota(i_2,t+n_1)), \text{ or } i_1=i_2=n,\ t+n_1\equiv 0{\hskip -0.5em}\pmod 4,
\end{cases}
\end{split}
\end{equation}
where $n_1$ is defined in (\ref{n_1 in type B_n}) for type $B_n$, the mapping $\iota$ is defined in Section \ref{Snake positions and minimal snake positions}.

Now we are ready for our main results in this section.

\begin{theorem}\label{S-systems}
In type $A_n$ (respectively, $B_n$), let $\mathcal{S}_2$ be the prime snake module (\ref{S_2}). We have the following system of equations
\begin{align*}
[\mathcal{S}_1]  [\mathcal{S}_2] = [\mathcal{S}_3] [\mathcal{S}_4] + [\mathcal{S}_5] [\mathcal{S}_6], \label{general form of an equation in S-system}
\end{align*}
where $\mathcal{S}_1$ is defined in Table \ref{definition S_1 in type A_{n}} (respectively, Table \ref{definition S_1 in type B_{n}}), $\mathcal{S}_3$,
$\mathcal{S}_4$ are defined in (\ref{S_3 S_4}), $\mathcal{S}_5$, $\mathcal{S}_6$  are defined in (\ref{S_5 S_6}), (\ref{S_5 S_6 in type $A_n$}) (respectively, (\ref{S_5 S_6 in type $B_n$})).
\end{theorem}

We call the system of equations in Theorem \ref{S-systems} the \textit{$S$-system} for type $A_n$ (respectively, $B_n$). In particular when $m=1$, the system of equations in Theorem \ref{S-systems} is T-system for types $A_n$ and $B_n$ \cite{Her06, KNS94}. The equations in the $S$-systems are different from the equations in the extended $T$-systems \cite{MY12b}.

Theorem \ref{S-systems} will be proved in Section \ref{proof of S-systems}. The next example gives some equations in the $S$-system for types $A_3$ and $B_4$.

\begin{example}
The following are some equations in the $S$-system for type $A_3$.
\begin{gather}
\begin{align*}
[3_{-3}3_{-1}][3_{-5}3_{-3}] &= [3_{-5}3_{-3}3_{-1}][3_{-3}]+ [2_{-4}2_{-2}],\\
[3_{-3}2_{0}][3_{-5}1_{-1}] &= [3_{-5}3_{-3}2_{0}][1_{-1}] + [2_{0}][2_{-4}1_{-1}],\\
[3_{-3}3_{-1}][3_{-5}2_{-2}] &= [3_{-5}3_{-3}3_{-1}][2_{-2}] + [3_{-1}][2_{-4}2_{-2}],\\
[3_{-7}1_{-3}2_{0}][3_{-9}2_{-6}1_{-3}2_{0}]&=[3_{-9}3_{-7}1_{-3}2_{0}][2_{-6}1_{-3}2_{0}]+[1_{-3}2_{0}][2_{-8}2_{-6}1_{-3}2_{0}],\\
[2_{-4}2_{-2}2_{0}][2_{-6}1_{-3}2_{0}] &= [2_{-6}2_{-4}2_{-2}2_{0}][1_{-3}2_{0}] + [3_{-5}2_{-2}2_{0}][1_{-5}1_{-3}2_{0}],\\
[2_{-4}2_{-2}2_{0}][2_{-6}3_{-3}2_{0}] &= [2_{-6}2_{-4}2_{-2}2_{0}][3_{-3}2_{0}] + [1_{-5}2_{-2}2_{0}][3_{-5}3_{-3}2_{0}],\\
[2_{-2}2_{0}][2_{-4}2_{0}] &= [2_{-4}2_{-2}2_{0}][2_{0}] + [1_{-3}2_{0}][3_{-3}2_{0}],\\
[1_{-5}1_{-3}2_{0}][1_{-7}2_{-4}2_{0}] &= [1_{-7}1_{-5}1_{-3}2_{0}][2_{-4}2_{0}]+[1_{-3}2_{0}][2_{-6}2_{-4}2_{0}],\\
[3_{-5}3_{-3}2_{0}][3_{-7}2_{-4}2_{0}] &= [3_{-7}3_{-5}3_{-3}2_{0}][2_{-4}2_{0}]+[3_{-3}2_{0}][2_{-6}2_{-4}2_{0}].
\end{align*}
\end{gather}

The following are some equations in the $S$-system for type $B_4$.
\begin{gather}
\begin{align*}
[4_{-20}4_{-18}2_{-11}3_{-1}][4_{-22}3_{-17}2_{-11}3_{-1}]&=[4_{-22}4_{-20}4_{-18}2_{-11}3_{-1}][3_{-17}2_{-11}3_{-1}]+[4_{-18}2_{-11}3_{-1}]
[3_{-21}3_{-17}2_{-11}3_{-1}],\\
[4_{-20}4_{-18}2_{-11}2_{-7}3_{-1}][4_{-22}1_{-13}2_{-7}3_{-1}]&=[4_{-22}4_{-20}4_{-18}2_{-11}2_{-7}3_{-1}][1_{-13}2_{-7}3_{-1}]+[4_{-18}2_{-11}2_{-7}3_{-1}][3_{-21}
1_{-13}2_{-7}3_{-1}],\\
[4_{-8}4_{-6}3_{-1}][4_{-10}3_{-1}]&=[4_{-10}4_{-8}4_{-6}3_{-1}][3_{-1}]+[4_{-6}3_{-1}][3_{-9}3_{-1}],\\
[4_{-12}4_{-10}4_{-8}4_{-6}3_{-1}][4_{-14}3_{-9}3_{-1}]&=[4_{-14}4_{-12}4_{-10}4_{-8}4_{-6}3_{-1}][3_{-9}3_{-1}]+[4_{-10}4_{-8}4_{-6}3_{-1}][3_{-13}3_{-9}3_{-1}],\\
[4_{-16}4_{-14}3_{-9}3_{-1}][4_{-18}2_{-11}3_{-1}]&=[4_{-18}4_{-16}4_{-14}3_{-9}3_{-1}][2_{-11}3_{-1}]+[4_{-14}3_{-9}3_{-1}][3_{-17}2_{-11}3_{-1}],\\
[2_{-13}4_{-6}3_{-1}][2_{-17}4_{-10}2_{-3}]&=[2_{-17}2_{-13}4_{-6}3_{-1}][4_{-10}2_{-3}]+[1_{-15}4_{-6}3_{-1}][3_{-15}4_{-10}2_{-3}],\\
[3_{-15}4_{-10}4_{-8}4_{-6}3_{-1}][3_{-19}4_{-14}3_{-9}3_{-1}]&=[3_{-19}3_{-15}4_{-10}4_{-8}4_{-6}3_{-1}][4_{-14}3_{-9}3_{-1}]+[2_{-17}4_{-10}4_{-8}4_{-6}3_{-1}][4_{-18}4_{-16}4_{-14}3_{-9}3_{-1}],\\
[3_{-15}4_{-10}2_{-3}][3_{-19}4_{-14}3_{-9}2_{-3}]&=[3_{-19}3_{-15}4_{-10}2_{-3}][4_{-14}3_{-9}2_{-3}]+[2_{-17}4_{-10}2_{-3}][4_{-18}4_{-16}4_{-14}3_{-9}2_{-3}],\\
[3_{-11}4_{-6}3_{-1}][3_{-15}4_{-10}3_{-1}]&=[3_{-15}3_{-11}4_{-6}3_{-1}][4_{-10}3_{-1}]+[2_{-13}4_{-6}3_{-1}][4_{-14}4_{-12}4_{-10}3_{-1}],\\
[2_{-23}3_{-17}2_{-11}3_{-1}][2_{-27}4_{-20}4_{-18}2_{-11}3_{-1}]&=[2_{-27}2_{-23}3_{-17}2_{-11}3_{-1}][4_{-20}4_{-18}2_{-11}3_{-1}]+[1_{-25}3_{-17}2_{-11}3_{-1}][3_{-25}4_{-20}4_{-18}2_{-11}3_{-1}],\\
[2_{-23}3_{-17}4_{-8}2_{-1}][2_{-27}4_{-20}4_{-18}4_{-8}2_{-1}]&=[2_{-27}2_{-23}3_{-17}4_{-8}2_{-1}][4_{-20}4_{-18}4_{-8}2_{-1}]+[1_{-25}3_{-17}4_{-8}2_{-1}][3_{-25}4_{-20}4_{-18}4_{-8}2_{-1}].
\end{align*}
\end{gather}
\end{example}

Moreover, we have the following theorem.
\begin{theorem}\label{equation simple}
The modules in the summands on the right-hand side of each equation in Theorem \ref{S-systems} are simple.
\end{theorem}
Theorem \ref{equation simple} will be proved in Section \ref{proof irreducible}.

\subsection{The $s$-systems of types $A_n$ and $B_{n}$}
Let $\mathcal{S}$ be a $U_q\widehat{\mathfrak{g}}$-module. We use $\res(\mathcal{S})$ to denote the restriction of $\mathcal{S}$ to $U_q\mathfrak{g}$. Let $\chi(M)$ be the character of a $U_q \mathfrak{g}$-module $M$. We have a system of equations
\begin{equation*}
\begin{split}
\chi(\res(\mathcal{S}_1))\chi(\res(\mathcal{S}_2))=\chi(\res(\mathcal{S}_3))\chi (\res(\mathcal{S}_4))+ \chi(\res(\mathcal{S}_5))\chi(\res(\mathcal{S}_6)),
\end{split}
\end{equation*}
where $[\mathcal{S}_1]  [\mathcal{S}_2] = [\mathcal{S}_3] [\mathcal{S}_4] + [\mathcal{S}_5] [\mathcal{S}_6]$ are equations of the $S$-system for type $A_n$ (respectively, $B_n$). We call this system of equations the \textit{$s$-system} of type $A_n$ (respectively, $B_n$).

\section{Relation between $S$-systems and cluster algebras} \label{relation between M-systems and cluster algebras}
In this section, we show that every equation in the $S$-system of type $A_n$ (respectively, $B_n$) corresponds to a mutation in some cluster algebra $\mathscr{A}$ (respectively, $\mathscr{A}'$) and every prime snake module of type $A_n$ (respectively, $B_n$)
corresponds to some cluster variable in $\mathscr{A}$ (respectively, $\mathscr{A}'$). In particular, this proves that the Hernandez-Leclerc conjecture (Conjecture \ref{Hernandez-Leclerc conjecture}) is true for prime snake modules of types $A_{n}$ and $B_{n}$.

\subsection{Definition of cluster algebras $\mathscr{A}$ and $\mathscr{A}'$}\label{definition of cluster algebra}
We recall the definition of the cluster algebras introduced in \cite{HL13}. Let $\widetilde{V}=I \times \mathbb{Z}$ and let $\widetilde{\Gamma}$ be a quiver with the vertex set $\widetilde{V}$ whose arrows are given by $(i,r) \to (j,s)$
if $b_{ij} \neq 0$ and $s = r + b_{ij}$, where $B = (b_{ij})_{i,j \in I} = DC$ is defined in Section \ref{definition of quantum affine algebras}.

It is shown in \cite{HL13} that $\widetilde{\Gamma}$ has two isomorphic components. Let $\Gamma$ be one of the components and $V$ its vertex set. Define a mapping $\psi$ by $\psi(i, t) = (i, t+d_i)$ for $(i, t) \in V$. Let $W \subset I \times \mathbb{Z}$
be the image of $V$ under the map $\psi$ and let $G$ be the same quiver as $\Gamma$ but with vertices labeled by $W$. Denote $W^-=W \cap (I \times \mathbb{Z}_{\leq 0})$ and let $Q$ be the full sub-quiver of $G$ with vertex set $W^-$.

Let $\mathbf{z}^-=\{z_{i,t}: (i,t)\in W^{-}\}$ and let $\mathscr{A}$ be the cluster algebra defined by the initial seed $(\mathbf{z}^-, Q)$. For convenience, we denote by $Q'$ and $\mathscr{A}'$ the quiver $Q$ and the cluster algebra $\mathscr{A}$ in the case of type $B_n$, respectively.

In the case of type $A_{n}$, let
\begin{align*}
{\bf s} =
     \{s^{(-2k+2)}_{k^{(i)}} \mid i \text{ is } \text{even}, \ k \in \mathbb{Z}_{\geq 1} \}  \cup
     \{s^{(-2k+1)}_{k^{(i)}} \mid i \text { is } \text{odd}, \ k \in \mathbb{Z}_{\geq 1} \}.
\end{align*}
In the case of type $B_{n}$, let ${\bf s}'={\bf s}_{1} \cup {\bf s}_{2}$, where
\begin{align*}
{\bf s}_{1} =\{ s^{(-2k+2)}_{k^{(n)}}\mid k \in \mathbb{Z}_{\geq 1} \},
\end{align*}
\begin{align*}
{\bf s}_{2} =\{ s^{(-4k+3)}_{k^{(i)}},\ s^{(-4k+1)}_{k^{(i)}} \mid i\in \{1,\ldots,n-1\}, \ k \in \mathbb{Z}_{\geq 1} \}.
\end{align*}

Let $\mathscr{A}$ (respectively, $\mathscr{A}'$) be the cluster algebra defined by the initial seed $({\bf s}, Q)$ (respectively, $({\bf s}', Q')$). Here we identify ${\bf s}$ (respectively, ${\bf s}'$) with ${\bf z}^-$ as follows.
For $(i, t) \in W^-$, we identify $s^{(t)}_{k^{(i)}}$ with $z_{i, t}$. We say that $s^{(t)}_{k^{(i)}}$ is at the vertex $(i, t)$ and we say that the label of this vertex is $(i,t)$.

Let $\widetilde{Q}$ (respectively, $\widetilde{Q}'$) be a quiver which is mutation equivalent to $Q$ (respectively, $Q'$) in type $A_n$ (respectively, $B_n$). After we mutate at a vertex $v$ of $\widetilde{Q}$
(respectively, $\widetilde{Q}'$), the variable at $v$ is changed and the label of $v$ is changed.

In this paper, our mutation sequences satisfy the property: after we mutate a quiver using a mutation sequence, any two vertices in the current quiver we obtain have different labels. If we mutate at a vertex $v$ with the label $(i,t)$ in $\widetilde{Q}$ (respectively, $\widetilde{Q}'$), then the label of $v$ becomes $(i,t-2d_i)$. We also use $(i,t)$ to denote the vertex with the label $(i,t)$.

\subsection{Fundamental segments and distinguished factors} \label{Fundamental segments and distinguished factors}

\begin{definition}\label{define fundamental segment}
Let $\mathcal{S}$ be a prime snake module and $S$ its highest $l$-weight monomial. Then $S$ can be written as
\[
S=S^{(t)}_{k_{1}^{(i_{1},j_{1})}, k_{2}^{(i_{2},j_{2})}, \ldots, k_{m-1}^{(i_{m-1},j_{m-1})}, k_{m}^{(i_{m})}},
\]
where $m\geq 1$, $j_\ell \geq 0$, $1 \leq \ell \leq m-1$, if $j_{\ell}=0$, then $i_{\ell} \neq i_{\ell+1}$, $k_1, k_2, \ldots, k_m \in \mathbb{Z}_{\geq 1}$, $t \in \mathbb{Z}$.

Let
\begin{align*}
& \mathcal{FS}(S) = FS_1 \cup FS_2 \cup FS_3,
\end{align*}
where
\begin{gather}
\begin{align*}
FS_1 = & \{S^{(t_{m})}_{k_{m}^{(i_m)}}: t_{m}=t+\sum_{j=1}^{m-1}n_{j} \}, \\ 
FS_2 = & \{ S^{(t_\ell)}_{k_{\ell}^{(i_\ell)}, k_{\ell+1}^{(i_{\ell+1})}, \ldots, k_{\ell+r-1}^{(i_{\ell+r-1})}, 1^{(i_{\ell+r})}} : \ell+r \leq m, \ r \geq 1, \ \ell \geq 1, \ t_\ell=t+\sum_{j=1}^{\ell-1}n_{j}, i_{\ell} \leq i_{\ell-1}, \\
& \ \  i_{\ell+r} \geq i_{\ell+r+1}, \ i_{\ell}< \cdots < i_{\ell+r} \text{ or } i_{\ell} \geq i_{\ell-1}, \ i_{\ell+r} \leq i_{\ell+r+1}, \ i_{\ell} > \cdots > i_{\ell+r} \}, \\ 
FS_3 = & \{ S^{(t_\ell)}_{k_{\ell}^{(i_\ell,j_\ell)},1^{(i_{\ell+1})}} : 1 \leq \ell \leq m-1, \ j_{\ell} \geq 1, \ t_\ell=t+\sum_{j=1}^{\ell-1}n_{j}\}. 
\end{align*}
\end{gather}
We call $\mathcal{FS}(S)$ the set of fundamental segments of $S$.
\end{definition}

\begin{example}
In type $A_5$, we have
\begin{align*}
\mathcal{FS}(2_{-12}4_{-8}5_{-5}5_{-3}4_0)=FS_1 \cup FS_2 \cup FS_3,
\end{align*}
where $FS_1=\{ 4_0 \}$, $FS_2=\{5_{-5}5_{-3}4_0, 2_{-12}4_{-8}5_{-5} \}$, $FS_3=\emptyset$.

In type $A_{4}$, we have
\begin{align*}
\mathcal{FS}(2_{-16}3_{-13}3_{-11}2_{-8}2_{-4}3_{-1})=FS_1 \cup FS_2 \cup FS_3,
\end{align*}
where $FS_1=\{3_{-1} \}$, $FS_2=\{ 2_{-4}3_{-1}, 3_{-13}3_{-11}2_{-8}, 2_{-16}3_{-13} \}$, $FS_3=\{2_{-8}2_{-4} \}$.
\begin{align*}
\mathcal{FS}(2_{-30}2_{-26}1_{-23}1_{-21}2_{-18}3_{-15}2_{-12}2_{-10}2_{-6}4_{-2}4_0)=FS_1 \cup FS_2 \cup FS_3,
\end{align*}
where $FS_1=\{4_{-2}4_0 \}$, $FS_2=\{2_{-6}4_{-2}, 3_{-15}2_{-12}, 1_{-23}1_{-21}2_{-18}3_{-15}, 2_{-26}1_{-23} \}$,
$FS_3=\{ 2_{-12}2_{-10}2_{-6}, 2_{-30}2_{-26} \}$.

In type $B_{3}$, we have
\begin{align*}
\mathcal{FS}(1_{-31}2_{-25}2_{-17}3_{-12}3_{-6}2_{-1})=FS_1 \cup FS_2 \cup FS_3,
\end{align*}
where $FS_1=\{2_{-1} \}, FS_2=\{3_{-6}2_{-1},2_{-17}3_{-12}, 1_{-31}2_{-25} \}, FS_3=\{3_{-12}3_{-6}, 2_{-25}2_{-17}\}$.
\begin{align*}
\mathcal{FS}(2_{-43}2_{-35}2_{-31}1_{-25}3_{-18}3_{-8}3_{-2}3_0) = FS_1 \cup FS_2 \cup FS_3,
\end{align*}
where
\begin{gather}
\begin{align*}
FS_1=\{ 3_{-2}3_0 \}, FS_2=\{1_{-25}3_{-18}, 2_{-35}2_{-31}1_{-25}\}, FS_3=\{3_{-8}3_{-2}, 3_{-18}3_{-8}, 2_{-43}2_{-35}\}.
\end{align*}
\end{gather}
\end{example}

The following proposition is easy to prove.
\begin{proposition}\label{fundamental segments}
Let $\mathcal{S}$ be a prime snake module with the highest $l$-weight monomial $S$. Then $S$ is uniquely determined by $\mathcal{FS}(S)$.
\end{proposition}

\begin{definition}\label{the definition of distinguished factor}
Let $\mathcal{S}$ be a prime snake module with the highest $l$-weight monomial $S$. Let $M$ be a monomial in $\mathcal{FS}(S)$. The last factor of $M$ is called the distinguished factor of $M$.
\end{definition}

\begin{example}
In type $A_{4}$, let $S=2_{-16}3_{-13}3_{-11}2_{-8}2_{-4}3_{-1}$. Then the set of distinguished factors of $S$ is $\{ 3_{-1}, 2_{-4}, 2_{-8}, 3_{-13} \}$, see Figure \ref{distinguished factors 1}.

In type $B_{3}$, let $S=1_{-31}2_{-25}2_{-17}3_{-12}3_{-6}2_{-1}$. Then the set of distinguished factors of $S$ is
$\{ 2_{-1}, 3_{-6}, 3_{-12}, 2_{-17}, 2_{-25} \}$, see Figure \ref{distinguished factors 2}.

\begin{figure}[H]
\resizebox{.6\width}{.6\height}{
\begin{minipage}[b]{0.55\linewidth}
\centerline{
\begin{tikzpicture}
\draw[step=.5cm,gray,thin] (-0.5,0) grid (2,8);
\draw[fill] (0,8.3) circle (2pt) -- (0.5,8.3) circle (2pt) --(1,8.3) circle (2pt) --(1.5,8.3) circle (2pt);
\draw [fill] (0.9,7.4) rectangle (1.1,7.6);
\draw [fill] (0.4,5.9) rectangle (0.6,6.1);
\draw [fill] (0.4,3.9) rectangle (0.6,4.1);
\draw [fill] (1,2.5) circle (3pt);
\draw [fill] (0.9,1.4) rectangle (1.1,1.6);
\draw [fill] (0.5,0) circle (3pt);
\node [above] at (0,8.3) {$1$};
\node [above] at (0.5,8.3) {$2$};
\node [above] at (1,8.3) {$3$};
\node [above] at (1.5,8.3) {$4$};
\node [left] at (-0.5,8) {$0$};
\node [left] at (-0.5,7.5) {$-1$};
\node [left] at (-0.5,7) {$-2$};
\node [left] at (-0.5,6.5) {$-3$};
\node [left] at (-0.5,6) {$-4$};
\node [left] at (-0.5,5.5) {$-5$};
\node [left] at (-0.5,5) {$-6$};
\node [left] at (-0.5,4.5) {$-7$};
\node [left] at (-0.5,4) {$-8$};
\node [left] at (-0.5,3.5) {$-9$};
\node [left] at (-0.5,3) {$-10$};
\node [left] at (-0.5,2.5) {$-11$};
\node [left] at (-0.5,2) {$-12$};
\node [left] at (-0.5,1.5) {$-13$};
\node [left] at (-0.5,1) {$-14$};
\node [left] at (-0.5,0.5) {$-15$};
\node [left] at (-0.5,0) {$-16$};
\end{tikzpicture}}
\caption{The points $\blacksquare$ corresponding to distinguished factors of $2_{-16}3_{-13}3_{-11}2_{-8}2_{-4}3_{-1}$ in type $A_{4}$.}\label{distinguished factors 1}
\end{minipage}
\begin{minipage}[b]{0.55\linewidth}
\centerline{
\begin{tikzpicture}
\draw[step=.5cm,gray,thin] (0,0) grid (5,8);
\draw[fill] (1,8.3) circle (2pt)--(2,8.3) circle (2pt);
\draw[fill] (3,8.3) circle (2pt)--(4,8.3) circle (2pt);
\draw[fill] (2.5,8.3) circle (2pt);
\draw [fill] (2.9,7.65) rectangle (3.1,7.85);
\draw [fill] (2.4,6.4) rectangle (2.6,6.6);
\draw [fill] (2.4,4.9) rectangle (2.6,5.1);
\draw [fill] (1.9,3.65)rectangle (2.1,3.85);
\draw [fill] (1.9,1.65) rectangle (2.1,1.85);
\draw [fill] (1,0.25) circle (3pt);
\node [above] at (1,8.3) {$1$};
\node [above] at (2,8.3) {$2$};
\node [above] at (2.5,8.3) {$3$};
\node [above] at (3,8.3) {$2$};
\node [above] at (4,8.3) {$1$};
\node [left] at (0,8) {$0$};
\node [left] at (0,7.5) {$-2$};
\node [left] at (0,7) {$-4$};
\node [left] at (0,6.5) {$-6$};
\node [left] at (0,6) {$-8$};
\node [left] at (0,5.5) {$-10$};
\node [left] at (0,5) {$-12$};
\node [left] at (0,4.5) {$-14$};
\node [left] at (0,4) {$-16$};
\node [left] at (0,3.5) {$-18$};
\node [left] at (0,3) {$-20$};
\node [left] at (0,2.5) {$-22$};
\node [left] at (0,2) {$-24$};
\node [left] at (0,1.5) {$-26$};
\node [left] at (0,1) {$-28$};
\node [left] at (0,0.5) {$-30$};
\node [left] at (0,0) {$-32$};
\draw [double,->] (2.075,8.3)--(2.425,8.3);
\draw [double,->] (2.925,8.3)--(2.575,8.3);
\end{tikzpicture}}
\caption{The points $\blacksquare$ corresponding to distinguished factors of $1_{-31}2_{-25}2_{-17}3_{-12}3_{-6}2_{-1}$ in type $B_{3}$.}\label{distinguished factors 2}
\end{minipage}}
\end{figure}
\end{example}

\subsection{Distinguished sub-quivers}\label{distinguished quivers}
Let $\widetilde{Q}$ (respectively, $\widetilde{Q}'$) be a quiver which is mutation equivalent to $Q$ (respectively, $Q'$) and any two vertices in $\widetilde{Q}$ (respectively, $\widetilde{Q}'$) have different labels.
We define a subset $\mathcal{Y}$ of the set of vertices in $\widetilde{Q}$ (respectively, $\widetilde{Q}'$) as follows.  In the case of type $A_{n}$, let
\begin{align*}
\mathcal{Y}= \{(i,k)\in \widetilde{Q}: i-k \equiv 0 {\hskip -0.6em}\pmod 2\}.
\end{align*}
In the case of type $B_{n}$, let
\begin{align*}
\mathcal{Y}= \{(n,2k) \in \widetilde{Q}': k\in \mathbb{Z}_{\leq 0} \} \cup \{(i,k) \in \widetilde{Q}': i<n \text{ and } k \equiv 1{\hskip -0.6em} \pmod 2\}.
\end{align*}

For a quiver $\mathscr{L}$, we use $V(\mathscr{L})$ to denote the set of its vertices. We define a distinguished sub-quiver $\mathscr{L}^{\widetilde{Q}}_{i,t}$
(respectively, $\mathscr{L}^{\widetilde{Q}'}_{i,t}$) with respect to $(i,t)\in V(\widetilde{Q})$ (respectively, $(i,t)\in V(\widetilde{Q}')$) in type $A_n$ (respectively, $B_n$).

The mapping $\iota$ is defined in Section \ref{Snake positions and minimal snake positions}. In the case of type $A_{n}$, let
\begin{align*}
V(\mathscr{L}^{\widetilde{Q}}_{i,t}) = \{& \iota(j,y_j) \in \mathcal{Y}: j \in I,\ y_{i}=t-2,\ y_{j}=y_{j+1}+1,\ 1 \leq j \leq i-1, \\
& \text{and } y_{j+1}=y_{j}+1,\ i \leq j \leq n-1 \}.
\end{align*}

Figures \ref{the line (4,0)}, \ref{the line (5,-3)} illustrate distinguished sub-quivers $\mathscr{L}^{Q}_{4,0}$, $\mathscr{L}^{Q}_{5,-3}$ in the original quiver of type $A_8$ respectively.
\begin{figure}[H]
\resizebox{.7\width}{.7\height}{
\begin{minipage}[b]{0.7\linewidth}
\centerline{
\begin{tikzpicture}
\draw[step=.5cm,gray,thin] (-0.5,0) grid (4,5);
\draw[fill] (0,5.3) circle (2pt)--(0.5,5.3) circle (2pt);
\draw[fill] (0.5,5.3) circle (2pt)--(1,5.3) circle (2pt);
\draw[fill] (1,5.3) circle (2pt)--(1.5,5.3) circle (2pt);
\draw[fill] (1.5,5.3) circle (2pt)--(2,5.3) circle (2pt);
\draw[fill] (2,5.3) circle (2pt)--(2.5,5.3) circle (2pt);
\draw[fill] (2.5,5.3) circle (2pt)--(3,5.3)circle (2pt);
\draw[fill] (3,5.3) circle (2pt)--(3.5,5.3)circle (2pt);
\draw [fill] (0.5,5) circle (1.5pt);
\draw [fill] (1,4.5) circle (1.5pt);
\draw [fill] (1.5,4) circle (1.5pt);
\draw [fill] (2,4.5) circle (1.5pt);
\draw [fill] (2.5,5) circle (1.5pt);
\node[above] at (0,5.3) {1};
\node[above] at  (0.5,5.3) {2};
\node[above] at  (1,5.3)   {3};
\node[above] at  (1.5,5.3) {4};
\node[above] at  (2,5.3)   {5};
\node[above] at  (2.5,5.3) {6};
\node[above] at  (3,5.3)   {7};
\node[above] at  (3.5,5.3) {8};
\node [left] at (-0.5,5) {0};
\node [left] at (-0.5,4.5) {-1};
\node [left] at (-0.5,4) {-2};
\node [left] at (-0.5,3.5) {-3};
\node [left] at (-0.5,3) {-4};
\node [left] at (-0.5,2.5) {-5};
\node [left] at (-0.5,2) {-6};
\node [left] at (-0.5,1.5) {-7};
\node [left] at (-0.5,1) {-8};
\node [left] at (-0.5,0.5) {-9};
\node [left] at (-0.5,0) {-10};
\begin{scope}[thin, every node/.style={sloped,allow upside down}]
\draw (0.5,5)--node {\midarrow}(1,4.5);
\draw (1,4.5)--node {\midarrow}(1.5,4);
\draw (2,4.5)--node {\midarrow}(1.5,4);
\draw (2.5,5)--node {\midarrow}(2,4.5);
\end{scope}
\end{tikzpicture}}
\caption{The distinguished sub-quiver $\mathscr{L}^{Q}_{4,0}$ in the original quiver of type $A_8$.}\label{the line (4,0)}
\end{minipage}
\begin{minipage}[b]{0.7\linewidth}
\centerline{
\begin{tikzpicture}
\draw[step=.5cm,gray,thin] (-0.5,0) grid (4,5);
\draw[fill] (0,5.3) circle (2pt)--(0.5,5.3) circle (2pt);
\draw[fill] (0.5,5.3) circle (2pt)--(1,5.3) circle (2pt);
\draw[fill] (1,5.3) circle (2pt)--(1.5,5.3) circle (2pt);
\draw[fill] (1.5,5.3) circle (2pt)--(2,5.3) circle (2pt);
\draw[fill] (2,5.3) circle (2pt)--(2.5,5.3) circle (2pt);
\draw[fill] (2.5,5.3) circle (2pt)--(3,5.3)circle (2pt);
\draw[fill] (3,5.3) circle (2pt)--(3.5,5.3)circle (2pt);
\draw [fill] (0,4.5) circle (1.5pt);
\draw [fill] (0.5,4) circle (1.5pt);
\draw [fill] (1,3.5) circle (1.5pt);
\draw [fill] (1.5,3) circle (1.5pt);
\draw [fill] (2,2.5) circle (1.5pt);
\draw [fill] (2.5,3) circle (1.5pt);
\draw [fill] (3,3.5) circle (1.5pt);
\draw [fill] (3.5,4) circle (1.5pt);
\node[above] at (0,5.3) {1};
\node[above] at  (0.5,5.3) {2};
\node[above] at  (1,5.3)   {3};
\node[above] at  (1.5,5.3) {4};
\node[above] at  (2,5.3)   {5};
\node[above] at  (2.5,5.3) {6};
\node[above] at  (3,5.3)   {7};
\node[above] at  (3.5,5.3) {8};
\node [left] at (-0.5,5) {0};
\node [left] at (-0.5,4.5) {-1};
\node [left] at (-0.5,4) {-2};
\node [left] at (-0.5,3.5) {-3};
\node [left] at (-0.5,3) {-4};
\node [left] at (-0.5,2.5) {-5};
\node [left] at (-0.5,2) {-6};
\node [left] at (-0.5,1.5) {-7};
\node [left] at (-0.5,1) {-8};
\node [left] at (-0.5,0.5) {-9};
\node [left] at (-0.5,0) {-10};
\begin{scope}[thin, every node/.style={sloped,allow upside down}]
\draw (0,4.5)--node {\midarrow}(0.5,4);
\draw (0.5,4)--node {\midarrow}(1,3.5);
\draw (1,3.5)--node {\midarrow}(1.5,3);
\draw (1.5,3)--node {\midarrow}(2,2.5);
\draw (2.5,3)--node {\midarrow}(2,2.5);
\draw (3,3.5)--node {\midarrow}(2.5,3);
\draw (3.5,4)--node {\midarrow}(3,3.5);
\end{scope}
\end{tikzpicture}}
\caption{The distinguished sub-quiver $\mathscr{L}^{Q}_{5,-3}$ in the original quiver of type $A_8$.}\label{the line (5,-3)}
\end{minipage}}
\end{figure}

In the case of type $B_{n}$,
\begin{itemize}
\item for all $i=n$, let
\begin{align*}
V(\mathscr{L}^{\widetilde{Q}'}_{n, \ell})= &\left\{ \iota(j,y_{j}) \in \mathcal{Y}: j\in I,\ y_{n}\in\{ \ell-4, \ell-6\},\ y_{n-1}=y_{n}+3, \right. \\ & \left.\ y_{j}=y_{j+1}+2,\ 1 \leq j \leq n-2 \right\};
\end{align*}

\item for all $i<n$, let
\begin{align*}
V(\mathscr{L}^{\widetilde{Q}'}_{i, \ell})=&\left\{\iota(j,y_{j}) \in \mathcal{Y}: j\in I,\ y_{i}=\ell-4,\ y_j=y_{j+1}+2, 1 \leq j \leq i-1, \right.\\
& \left. y_{j+1}= y_{j}+2,\ i \leq j \leq n-2,\  y_{n}=y_{n-1}+1 \right\} \cup \\
& \left\{ \iota(j,y_{j}) \in \mathcal{Y}: j\in I,\ y_{n}=\ell+2n-2i-7,\ y_{n-1}=y_{n}+3, \right. \\ & \left.\ y_{j}=y_{j+1}+2,\ 1 \leq j \leq n-2 \right\}.
\end{align*}
\end{itemize}

Figures \ref{the line (4, -4)}, \ref{the line (4,-6)}, \ref{the line (2,-7)}, \ref{the line (2,-5)} illustrate the distinguished sub-quivers $\mathscr{L}^{Q'}_{4,-4}$, $\mathscr{L}^{Q'}_{4,-6}$, $\scalemath{0.92}{\mathscr{L}^{Q'}_{2,-7}}$, $\mathscr{L}^{Q'}_{2,-5}$ in the original quiver of type $B_4$ respectively.
\begin{figure}[H]
\resizebox{.6\width}{.6\height}{
\begin{minipage}[b]{0.8\linewidth}
\centerline{
\begin{tikzpicture}
\draw[step=.5cm,gray,thin] (-1,0) grid (6,8);
\draw[fill] (0,8.3) circle (2pt)--(1,8.3) circle (2pt);
\draw[fill] (1,8.3) circle (2pt)--(2,8.3) circle (2pt);
\draw[fill] (3,8.3) circle (2pt)--(4,8.3) circle (2pt);
\draw[fill] (4,8.3) circle (2pt)--(5,8.3) circle (2pt);
\draw[fill] (2.5,8.3) circle (2pt);
\node [above] at (0,8.3) {1};
\node [above] at (1,8.3) {2};
\node [above] at (2,8.3) {3};
\node [above] at (2.5,8.3) {4};
\node [above] at (3,8.3) {3};
\node [above] at (4,8.3) {2};
\node [above] at (5,8.3) {1};
\draw [double,->] (2.075,8.3)--(2.425,8.3);
\draw [double,->] (2.925,8.3)--(2.575,8.3);
\node [left] at (-1,8) {0};
\node [left] at (-1,7) {-2};
\node [left] at (-1,6) {-4};
\node [left] at (-1,5) {-6};
\node [left] at (-1,4) {-8};
\node [left] at (-1,3) {-10};
\node [left] at (-1,2) {-12};
\node [left] at (-1,1) {-14};
\node [left] at (-1,0) {-16};
\draw[fill] (0,7.5) circle (1.5pt);
\draw[fill] (1,6.5) circle (1.5pt);
\draw[fill] (2,5.5) circle (1.5pt);
\draw[fill] (2.5,4) circle (1.5pt);
\draw[fill] (2.5,3) circle (1.5pt);
\draw[fill] (3,4.5) circle (1.5pt);
\draw[fill] (4,5.5) circle (1.5pt);
\draw[fill] (5,6.5) circle (1.5pt);
\begin{scope}[thin, every node/.style={sloped,allow upside down}]
\draw (0,7.5)--node {\midarrow}(1.5,6);
\draw (1.5,6)--node {\midarrow}(2,5.5);
\draw (2,5.5)--node {\midarrow}(2.5,4);
\draw (2.5,3)--node {\midarrow}(2.5,4);
\draw (3,4.5)--node {\midarrow}(2.5,3);
\draw (3.5,5)--node {\midarrow}(3,4.5);
\draw (5,6.5)--node {\midarrow}(3.5,5);
\end{scope}
\end{tikzpicture}}
\caption{The distinguished sub-quiver $\mathscr{L}^{Q'}_{4, -4}$ in the original quiver of type $B_4$.} \label{the line (4, -4)}
\end{minipage}
\begin{minipage}[b]{0.8\linewidth}
\centerline{
\begin{tikzpicture}
\draw[step=.5cm,gray,thin] (-1,0) grid (6,8);
\draw[fill] (0,8.3) circle (2pt)--(1,8.3) circle (2pt);
\draw[fill] (1,8.3) circle (2pt)--(2,8.3) circle (2pt);
\draw[fill] (3,8.3) circle (2pt)--(4,8.3) circle (2pt);
\draw[fill] (4,8.3) circle (2pt)--(5,8.3) circle (2pt);
\draw[fill] (2.5,8.3) circle (2pt);
\node [above] at (0,8.3) {1};
\node [above] at (1,8.3) {2};
\node [above] at (2,8.3) {3};
\node [above] at (2.5,8.3) {4};
\node [above] at (3,8.3) {3};
\node [above] at (4,8.3) {2};
\node [above] at (5,8.3) {1};
\draw [double,->] (2.075,8.3)--(2.425,8.3);
\draw [double,->] (2.925,8.3)--(2.575,8.3);
\node [left] at (-1,8) {0};
\node [left] at (-1,7) {-2};
\node [left] at (-1,6) {-4};
\node [left] at (-1,5) {-6};
\node [left] at (-1,4) {-8};
\node [left] at (-1,3) {-10};
\node [left] at (-1,2) {-12};
\node [left] at (-1,1) {-14};
\node [left] at (-1,0) {-16};
\draw[fill] (0,5.5) circle (1.5pt);
\draw[fill] (1,4.5) circle (1.5pt);
\draw[fill] (2,3.5) circle (1.5pt);
\draw[fill] (2.5,3) circle (1.5pt);
\draw[fill] (2.5,2) circle (1.5pt);
\draw[fill] (3,4.5) circle (1.5pt);
\draw[fill] (4,5.5) circle (1.5pt);
\draw[fill] (5,6.5) circle (1.5pt);
\begin{scope}[thin, every node/.style={sloped,allow upside down}]
\draw (0,5.5)--node {\midarrow}(1.5,4);
\draw (1.5,4)--node {\midarrow}(2,3.5);
\draw (2,3.5)--node {\midarrow}(2.5,2);
\draw (2.5,2)--node {\midarrow}(2.5,3);
\draw (3,4.5)--node {\midarrow}(2.5,3);
\draw (3.5,5)--node {\midarrow}(3,4.5);
\draw (5,6.5)--node {\midarrow}(3.5,5);
\end{scope}
\end{tikzpicture}}
\caption{The distinguished sub-quiver $\mathscr{L}^{Q'}_{4, -6}$ in the original quiver of type $B_4$.}\label{the line (4,-6)}
\end{minipage}}
\end{figure}

\begin{figure}[H]
\resizebox{.6\width}{.6\height}{
\begin{minipage}[b]{0.8\linewidth}
\centerline{
\begin{tikzpicture}
\draw[step=.5cm,gray,thin] (-1,0) grid (6,8);
\draw[fill] (0,8.3) circle (2pt)--(1,8.3) circle (2pt);
\draw[fill] (1,8.3) circle (2pt)--(2,8.3) circle (2pt);
\draw[fill] (3,8.3) circle (2pt)--(4,8.3) circle (2pt);
\draw[fill] (4,8.3) circle (2pt)--(5,8.3) circle (2pt);
\draw[fill] (2.5,8.3) circle (2pt);
\node [above] at (0,8.3) {1};
\node [above] at (1,8.3) {2};
\node [above] at (2,8.3) {3};
\node [above] at (2.5,8.3) {4};
\node [above] at (3,8.3) {3};
\node [above] at (4,8.3) {2};
\node [above] at (5,8.3) {1};
\draw [double,->] (2.075,8.3)--(2.425,8.3);
\draw [double,->] (2.925,8.3)--(2.575,8.3);
\node [left] at (-1,8) {0};
\node [left] at (-1,7) {-2};
\node [left] at (-1,6) {-4};
\node [left] at (-1,5) {-6};
\node [left] at (-1,4) {-8};
\node [left] at (-1,3) {-10};
\node [left] at (-1,2) {-12};
\node [left] at (-1,1) {-14};
\node [left] at (-1,0) {-16};
\draw[fill] (0,3.5) circle (1.5pt);
\draw[fill] (1,2.5) circle (1.5pt);
\draw[fill] (2,3.5) circle (1.5pt);
\draw[fill] (2.5,3) circle (1.5pt);
\draw[fill] (2.5,4) circle (1.5pt);
\draw[fill] (3,4.5) circle (1.5pt);
\draw[fill] (4,5.5) circle (1.5pt);
\draw[fill] (5,6.5) circle (1.5pt);
\begin{scope}[thin, every node/.style={sloped,allow upside down}]
\draw (0,3.5)--(0.5,3);
\draw (0.5,3)--node {\midarrow}(1,2.5);
\draw (2,3.5)--(1.5,3);
\draw (1.5,3)--node {\midarrow}(1,2.5);
\draw (2.5,4)--node {\midarrow}(2,3.5);
\draw (2.5,3)--node {\midarrow}(2.5,4);
\draw (3,4.5)--node {\midarrow}(2.5,3);
\draw (3.5,5)--node {\midarrow}(3,4.5);
\draw (5,6.5)--node {\midarrow}(3.5,5);
\end{scope}
\end{tikzpicture}}
\caption{The distinguished sub-quiver $\mathscr{L}^{Q'}_{2,-7}$ in the original quiver of type $B_4$.}\label{the line (2,-7)}
\end{minipage}
\begin{minipage}[b]{0.8\linewidth}
\centerline{
\begin{tikzpicture}
\draw[step=.5cm,gray,thin] (-1,0) grid (6,8);
\draw[fill] (0,8.3) circle (2pt)--(1,8.3) circle (2pt);
\draw[fill] (1,8.3) circle (2pt)--(2,8.3) circle (2pt);
\draw[fill] (3,8.3) circle (2pt)--(4,8.3) circle (2pt);
\draw[fill] (4,8.3) circle (2pt)--(5,8.3) circle (2pt);
\draw[fill] (2.5,8.3) circle (2pt);
\node [above] at (0,8.3) {1};
\node [above] at (1,8.3) {2};
\node [above] at (2,8.3) {3};
\node [above] at (2.5,8.3) {4};
\node [above] at (3,8.3) {3};
\node [above] at (4,8.3) {2};
\node [above] at (5,8.3) {1};
\draw [double,->] (2.075,8.3)--(2.425,8.3);
\draw [double,->] (2.925,8.3)--(2.575,8.3);
\node [left] at (-1,8) {0};
\node [left] at (-1,7) {-2};
\node [left] at (-1,6) {-4};
\node [left] at (-1,5) {-6};
\node [left] at (-1,4) {-8};
\node [left] at (-1,3) {-10};
\node [left] at (-1,2) {-12};
\node [left] at (-1,1) {-14};
\node [left] at (-1,0) {-16};
\draw[fill] (0,7.5) circle (1.5pt);
\draw[fill] (1,6.5) circle (1.5pt);
\draw[fill] (2,5.5) circle (1.5pt);
\draw[fill] (2.5,4) circle (1.5pt);
\draw[fill] (2.5,5) circle (1.5pt);
\draw[fill] (3,4.5) circle (1.5pt);
\draw[fill] (4,3.5) circle (1.5pt);
\draw[fill] (5,4.5) circle (1.5pt);
\begin{scope}[thin, every node/.style={sloped,allow upside down}]
\draw (0,7.5)--node {\midarrow}(1.5,6);
\draw (1.5,6)--node {\midarrow}(2,5.5);
\draw (2,5.5)--node {\midarrow}(2.5,4);
\draw (2.5,4)--node {\midarrow}(2.5,5);
\draw (2.5,5)--node {\midarrow}(3,4.5);
\draw (3,4.5)--(3.5,4);
\draw (3.5,4)--node {\midarrow}(4,3.5);
\draw (5,4.5)--(4.5,4);
\draw (4.5,4)--node {\midarrow}(4,3.5);
\end{scope}
\end{tikzpicture}}
\caption{The distinguished sub-quiver $\mathscr{L}^{Q'}_{2,-5}$ in the original quiver of type $B_4$.}\label{the line (2,-5)}
\end{minipage}}
\end{figure}

We define two sub-quivers $\mathfrak{l}(\mathscr{L}^{\widetilde{Q}}_{i,t})$ and $\mathfrak{r}(\mathscr{L}^{\widetilde{Q}}_{i,t})$ as follows.
\begin{align*}
V(\mathfrak{l}(\mathscr{L}^{\widetilde{Q}}_{i,t}))=V(\mathscr{L}^{\widetilde{Q}}_{i,t}) \cap \{\iota(j,y_j) \in \mathcal{Y}:\pr\nolimits_1(\iota(j,y_j))<\pr\nolimits_ 1(\iota(i,t)) \},
\end{align*}
\begin{align*}
V(\mathfrak{r}(\mathscr{L}^{\widetilde{Q}}_{i,t}))=V(\mathscr{L}^{\widetilde{Q}}_{i,t}) \cap \{\iota(j,y_j) \in \mathcal{Y}:\pr\nolimits_1(\iota(j,y_j))>\pr\nolimits_1(\iota(i,t)) \}.
\end{align*}
Similarly, we may define the sub-quivers $\mathfrak{l}(\mathscr{L}^{\widetilde{Q}'}_{i,t})$ and $\mathfrak{r}(\mathscr{L}^{\widetilde{Q}'}_{i,t})$.

\subsection{Mutation sequences with respect to a quiver}
Let $\widetilde{Q}$ (respectively, $\widetilde{Q}'$) be a quiver which is mutation equivalent to $Q$ (respectively, $Q'$) and any two vertices in $\widetilde{Q}$ (respectively, $\widetilde{Q}'$) have different labels.
By saying that we mutate $C_{i,t}$ of $\widetilde{Q}$ (respectively, $\widetilde{Q}'$), we mean that we mutate at the vertex of $\widetilde{Q}$ (respectively, $\widetilde{Q}'$) which has the label $(i, t)$ in the $i$th column and so on until the vertex at infinity in the $i$th column.

For convenience, in the case of type $B_n$, let
\begin{gather}
\begin{align*}
&V'(\mathfrak{l}(\mathscr{L}^{\widetilde{Q}'}_{i,\ell}))=
\begin{cases}
V(\mathfrak{l}(\mathscr{L}^{\widetilde{Q}'}_{i,\ell})), & (i\neq n,\ 2n-2i+\ell \equiv 1 {\hskip -0.6em}\pmod4) \text{ or } i=n, \\
V(\mathfrak{l}(\mathscr{L}^{\widetilde{Q}'}_{i,\ell}))-\{(2n-1, \ell+2n-2i-7)\}, & i\neq n,\ 2n-2i+\ell \equiv 3 {\hskip -0.6em}\pmod4,
\end{cases}\\
&V'(\mathfrak{r}(\mathscr{L}^{\widetilde{Q}'}_{i,\ell}))=
\begin{cases}
V(\mathfrak{r}(\mathscr{L}^{\widetilde{Q}'}_{i,\ell}))-\{(2n-1, \ell+2n-2i-7)\}, & i\neq n,\ 2n-2i+\ell \equiv 1 {\hskip -0.6em} \pmod4,\\
V(\mathfrak{r}(\mathscr{L}^{\widetilde{Q}'}_{i,\ell})), & (i\neq n,\ 2n-2i+\ell \equiv 3 {\hskip -0.6em} \pmod 4) \text{ or } i=n.
\end{cases}
\end{align*}
\end{gather}

For type $A_n$ (respectively, $B_n$), suppose that
\[
V(\mathfrak{l}(\mathscr{L}^{\widetilde{Q}}_{i,\ell})) \ (respectively, \ V'(\mathfrak{l}(\mathscr{L}^{\widetilde{Q}'}_{i,\ell}))) = \{(j_1,t_1), (j_2,t_2), \ldots, (j_m,t_m)\},
\]
where $j_1 < j_2 < \cdots < j_m < \pr_1(\iota(i,\ell))$. We say that we mutate $\mathfrak{l}(\mathscr{L}^{\widetilde{Q}}_{i,\ell})$ (respectively, $\mathfrak{l}(\mathscr{L}^{\widetilde{Q}'}_{i,\ell})$), we mean that we mutate $C_{j_1,t_1}$, $C_{j_2,t_2}$, $\ldots$, $C_{j_m,t_m}$. Suppose that
\[
V(\mathfrak{r}(\mathscr{L}^{\widetilde{Q}}_{i,\ell})) \ (respectively, \ V'(\mathfrak{r}(\mathscr{L}^{\widetilde{Q}'}_{i,\ell}))) = \{(j_1,t_1), (j_2,t_2), \ldots, (j_m,t_m)\},
\]
where $j_1 > j_2 > \cdots > j_m > \pr_1(\iota(i,\ell))$. We say that we mutate $\mathfrak{r}(\mathscr{L}^{\widetilde{Q}}_{i,\ell})$ (respectively, $\mathfrak{r}(\mathscr{L}^{\widetilde{Q}'}_{i,\ell})$), we mean that we mutate $C_{j_1,t_1}$, $C_{j_2,t_2}$, $\ldots$, $C_{j_m,t_m}$.

We say that we mutate $\mathscr{L}^{\widetilde{Q}}_{i,\ell}$ (respectively, $\mathscr{L}^{\widetilde{Q}'}_{i,\ell}$), we mean that we mutate $\mathfrak{l}(\mathscr{L}^{\widetilde{Q}}_{i,\ell})$ (respectively, $\mathfrak{l}(\mathscr{L}^{\widetilde{Q}'}_{i,\ell})$), $\mathfrak{r}(\mathscr{L}^{\widetilde{Q}}_{i,\ell})$ (respectively, $\mathfrak{r}(\mathscr{L}^{\widetilde{Q}'}_{i,\ell})$), and then mutate $C_{i,\ell-2}$ (respectively, $C_{i,\ell-4}$).

\subsection{Definitions of the mappings $\tau_l$, $\tau_r$, $\tau$.}
Let $\widetilde{Q}_1$ (respectively, $\widetilde{Q}'_1$) be a quiver which is mutation equivalent to $\widetilde{Q}$ (respectively, $\widetilde{Q}'$).

Let $\mathfrak{L}=\{ \mathscr{L}^{\widetilde{Q}}_{i,t} : (i,t) \in \mathcal{Y}\}$ (respectively, $\{ \mathscr{L}^{\widetilde{Q}'}_{i,t} : (i,t) \in \mathcal{Y}\}$). We define three mappings $\tau_l$, $\tau_r$, $\tau$ on $\mathfrak{L}$ as follows. In the case of type $A_{n}$, let
\begin{align}\label{tau_l in type A_n}
\tau_l(\mathscr{L}^{\widetilde{Q}}_{i, t}) &  = \mathscr{L}^{\widetilde{Q}_1}_{i-1, t-1},
\end{align}
\begin{align}\label{tau_r in type A_n}
\tau_r(\mathscr{L}^{\widetilde{Q}}_{i, t}) &  = \mathscr{L}^{\widetilde{Q}_1}_{i+1, t-1},
\end{align}
\begin{align}\label{tau in type A_n}
\tau(\mathscr{L}^{\widetilde{Q}}_{i, t})   &  = \mathscr{L}^{\widetilde{Q}_1}_{i, t-2},
\end{align}
where the quivers $\widetilde{Q}_1$'s in (\ref{tau_l in type A_n}), (\ref{tau_r in type A_n}), (\ref{tau in type A_n}) are obtained from $\widetilde{Q}$ by mutating $\mathfrak{l}(\mathscr{L}^{\widetilde{Q}}_{i, t})$, $\mathfrak{r}(\mathscr{L}^{\widetilde{Q}}_{i, t})$ and $\mathscr{L}^{\widetilde{Q}}_{i, t}$,  respectively.

In the case of type $B_{n}$, let
\begin{align}\label{tau_l in type B_n}
\tau_l(\mathscr{L}^{\widetilde{Q}'}_{i, t})=
\begin{cases}
\mathscr{L}^{\widetilde{Q}'_1}_{i-1, t-2}, & i \leq n-1,\ 2n+t-2i \equiv 1 {\hskip -0.6em}\pmod 4,\\
\mathscr{L}^{\widetilde{Q}'_1}_{n-1, t-3}, & i = n,\ t \equiv 2 {\hskip -0.6em}\pmod 4,\\
\mathscr{L}^{\widetilde{Q}'_1}_{n-1, t-1}, & i = n,\ t \equiv 0 {\hskip -0.6em}\pmod 4,\\
\mathscr{L}^{\widetilde{Q}'_1}_{n, t-3},   & i = n-1,\ 2n+t-2i \equiv 3 {\hskip -0.6em}\pmod 4,  \\
\mathscr{L}^{\widetilde{Q}'_1}_{i+1, t-2}, & i < n-1,\ 2n+t-2i \equiv 3 {\hskip -0.6em}\pmod 4, \\
\end{cases}
\end{align}
\begin{align}\label{tau_r in type B_n}
\tau_r(\mathscr{L}^{\widetilde{Q}'}_{i, t})=
\begin{cases}
\mathscr{L}^{\widetilde{Q}'_1}_{i+1, t-2}, & i < n-1,\ 2n+t-2i \equiv 1 {\hskip -0.6em}\pmod 4,\\
\mathscr{L}^{\widetilde{Q}'_1}_{n, t-3},   & i=n-1,\ 2n+t-2i \equiv 1 {\hskip -0.6em}\pmod 4, \\
\mathscr{L}^{\widetilde{Q}'_1}_{n-1, t-1}, & i = n,\ t \equiv 2 {\hskip -0.6em}\pmod 4,\\
\mathscr{L}^{\widetilde{Q}'_1}_{n-1, t-3}, & i = n,\ t \equiv 0 {\hskip -0.6em}\pmod 4,\\
\mathscr{L}^{\widetilde{Q}'_1}_{i-1, t-2}, & i \leq n-1 \text{ and } 2n+t-2i \equiv 3 {\hskip -0.6em}\pmod 4, \\
\end{cases}
\end{align}
\begin{align}\label{tau in type B_n}
\tau(\mathscr{L}^{\widetilde{Q}'}_{i, t})=\mathscr{L}^{\widetilde{Q}'_1}_{i, t-4},
\end{align}
where the quivers $\widetilde{Q}'_1$'s in (\ref{tau_l in type B_n}), (\ref{tau_r in type B_n}), (\ref{tau in type B_n}) are obtained from $\widetilde{Q}'$ by mutating $\mathfrak{l}(\mathscr{L}^{\widetilde{Q}'}_{i, t})$, $\mathfrak{r}(\mathscr{L}^{\widetilde{Q}'}_{i, t})$ and $\mathscr{L}^{\widetilde{Q}'}_{i, t}$,  respectively.

Let $\tau_{l}^0(\mathscr{L}^{\widetilde{Q}}_{i,t})=\mathscr{L}^{\widetilde{Q}}_{i,t}$ and  $\tau_{l}^{m}(\mathscr{L}^{\widetilde{Q}}_{i,t})=\tau_l(\tau_{l}^{m-1}(\mathscr{L}^{\widetilde{Q}}_{i,t}))$. We use the following convention: if $m < 0 $, then $\tau_{l}^{m}(\mathscr{L}^{\widetilde{Q}}_{i,t})=\emptyset$. The quivers $\tau_{r}^{m}(\mathscr{L}^{\widetilde{Q}}_{i,t})$, $\tau^{m}(\mathscr{L}^{\widetilde{Q}}_{i,t})$, $\tau_{l}^{m}(\mathscr{L}^{\widetilde{Q}'}_{i,t})$, $\tau_{r}^{m}(\mathscr{L}^{\widetilde{Q}'}_{i,t})$, $\tau^{m}(\mathscr{L}^{\widetilde{Q}'}_{i,t})$ are defined similarly.

By definition, $\tau_{l}^{m}(\mathscr{L}^{\widetilde{Q}}_{i,t})$ (respectively, $\tau_{l}^{m}(\mathscr{L}^{\widetilde{Q}'}_{i,t})$) is a sub-quiver of some quiver $\widetilde{Q}_m$ (respectively, $\widetilde{Q}'_m$) which is mutation equivalent to $\widetilde{Q}$ (respectively, $\widetilde{Q}'$). For simplicity, we write $\mathscr{L}_{i',t'}=\varphi^{m}(\mathscr{L}_{i,t})$ for $\mathscr{L}^{\widetilde{Q}_m}_{i',t'}=\varphi^{m}(\mathscr{L}^{\widetilde{Q}}_{i,t})$ and $\mathscr{L}^{\widetilde{Q}'_m}_{i',t'}=\varphi^{m}(\mathscr{L}^{\widetilde{Q}'}_{i,t})$, where $\varphi \in \{ \tau_l, \ \tau_r, \ \tau \}$.

\begin{example}
Figures \ref{A_8} and \ref{B_4} illustrate the mappings $\tau_l$, $\tau_r$, $\tau$ in the original quivers of types $A_{8}$, $B_{4}$, respectively.
\begin{figure}[H]
\resizebox{.7\width}{.7\height}{
\begin{minipage}[b]{0.35\linewidth}
\centerline{
\begin{tikzpicture}
\draw[step=.5cm,gray,thin] (-0.5,0) grid (4,5);
\draw[fill] (0,5.3) circle (2pt)--(0.5,5.3) circle (2pt);
\draw[fill] (0.5,5.3) circle (2pt)--(1,5.3) circle (2pt);
\draw[fill] (1,5.3) circle (2pt)--(1.5,5.3) circle (2pt);
\draw[fill] (1.5,5.3) circle (2pt)--(2,5.3) circle (2pt);
\draw[fill] (2,5.3) circle (2pt)--(2.5,5.3) circle (2pt);
\draw[fill] (2.5,5.3) circle (2pt)--(3,5.3)circle (2pt);
\draw[fill] (3,5.3) circle (2pt)--(3.5,5.3)circle (2pt);
\draw [fill] (0,4.5) circle (1.5pt);
\draw [fill,red] (0,3.5) circle (1.5pt);
\draw [fill] (0.5,4) circle (1.5pt);
\draw [fill,red] (0.5,3) circle (1.5pt);
\draw [fill] (1,3.5) circle (1.5pt);
\draw [fill,red] (1,2.5) circle (1.5pt);
\draw [fill] (1.5,3) circle (1.5pt);
\draw [fill,red] (1.5,2) circle (1.5pt);
\draw [fill] (2,2.5) circle (1.5pt);
\draw [fill] (2.5,3) circle (1.5pt);
\draw [fill] (3,3.5) circle (1.5pt);
\draw [fill] (3.5,4) circle (1.5pt);
\node[above] at (1,3.8) {$\bf{\mathscr{L}^{Q}_{5,-3}}$};
\node[red,left] at (1.6,2) {$\bf{\mathscr{L}^{Q_1}_{4,-4}}$};
\node[above] at (0,5.3) {$1$};
\node[above] at  (0.5,5.3) {$2$};
\node[above] at  (1,5.3)   {$3$};
\node[above] at  (1.5,5.3) {$4$};
\node[above] at  (2,5.3)   {$5$};
\node[above] at  (2.5,5.3) {$6$};
\node[above] at  (3,5.3)   {$7$};
\node[above] at  (3.5,5.3) {$8$};
\node [left] at (-0.5,5) {$0$};
\node [left] at (-0.5,4.5) {$-1$};
\node [left] at (-0.5,4) {$-2$};
\node [left] at (-0.5,3.5) {$-3$};
\node [left] at (-0.5,3) {$-4$};
\node [left] at (-0.5,2.5) {$-5$};
\node [left] at (-0.5,2) {$-6$};
\node [left] at (-0.5,1.5) {$-7$};
\node [left] at (-0.5,1) {$-8$};
\node [left] at (-0.5,0.5) {$-9$};
\node [left] at (-0.5,0) {$-10$};
\node [left] at (2,-0.6) {$(\text{a})$};
\begin{scope}[thin, every node/.style={sloped,allow upside down}]
\draw (0,4.5)--node {\midarrow}(0.5,4);
\draw (0.5,4)--node {\midarrow}(1,3.5);
\draw (1,3.5)--node {\midarrow}(1.5,3);
\draw (1.5,3)--node {\midarrow}(2,2.5);
\draw[red] (0,3.5)--node {\midarrow}(0.5,3);
\draw[red] (0.5,3)--node {\midarrow}(1,2.5);
\draw[red] (1,2.5)--node {\midarrow}(1.5,2);
\draw[red] (2,2.5)--node {\midarrow}(1.5,2);
\draw (2.5,3)--node {\midarrow}(2,2.5);
\draw (3,3.5)--node {\midarrow}(2.5,3);
\draw(3.5,4)--node {\midarrow}(3,3.5);
\end{scope}
\end{tikzpicture}}
\end{minipage}
\begin{minipage}[b]{0.35\linewidth}
\centerline{
\begin{tikzpicture}
\draw[step=.5cm,gray,thin] (-0.5,0) grid (4,5);
\draw[fill] (0,5.3) circle (2pt)--(0.5,5.3) circle (2pt);
\draw[fill] (0.5,5.3) circle (2pt)--(1,5.3) circle (2pt);
\draw[fill] (1,5.3) circle (2pt)--(1.5,5.3) circle (2pt);
\draw[fill] (1.5,5.3) circle (2pt)--(2,5.3) circle (2pt);
\draw[fill] (2,5.3) circle (2pt)--(2.5,5.3) circle (2pt);
\draw[fill] (2.5,5.3) circle (2pt)--(3,5.3)circle (2pt);
\draw[fill] (3,5.3) circle (2pt)--(3.5,5.3)circle (2pt);
\draw [fill] (0,4.5) circle (1.5pt);
\draw [fill] (0.5,4) circle (1.5pt);
\draw [fill] (1,3.5) circle (1.5pt);
\draw [fill] (1.5,3) circle (1.5pt);
\draw [fill] (2,2.5) circle (1.5pt);
\draw [fill] (2.5,3) circle (1.5pt);
\draw [fill,red] (2.5,2) circle (1.5pt);
\draw [fill] (3,3.5) circle (1.5pt);
\draw [fill,red] (3,2.5) circle (1.5pt);
\draw [fill] (3.5,4) circle (1.5pt);
\draw [fill,red] (3.5,3) circle (1.5pt);
\node[above] at (3,3.85) {$\bf{\mathscr{L}^{Q}_{5,-3}}$};
\node[red,right] at (2.4,1.9) {$\bf{\mathscr{L}^{Q_1}_{6,-4}}$};
\node[above] at (0,5.3) {$1$};
\node[above] at  (0.5,5.3) {$2$};
\node[above] at  (1,5.3)   {$3$};
\node[above] at  (1.5,5.3) {$4$};
\node[above] at  (2,5.3)   {$5$};
\node[above] at  (2.5,5.3) {$6$};
\node[above] at  (3,5.3)   {$7$};
\node[above] at  (3.5,5.3) {$8$};
\node [left] at (-0.5,5) {$0$};
\node [left] at (-0.5,4.5) {$-1$};
\node [left] at (-0.5,4) {$-2$};
\node [left] at (-0.5,3.5) {$-3$};
\node [left] at (-0.5,3) {$-4$};
\node [left] at (-0.5,2.5) {$-5$};
\node [left] at (-0.5,2) {$-6$};
\node [left] at (-0.5,1.5) {$-7$};
\node [left] at (-0.5,1) {$-8$};
\node [left] at (-0.5,0.5) {$-9$};
\node [left] at (-0.5,0) {$-10$};
\node [left] at (2,-0.6) {$(\text{b})$};
\begin{scope}[thin, every node/.style={sloped,allow upside down}]
\draw (0,4.5)--node {\midarrow}(0.5,4);
\draw (0.5,4)--node {\midarrow}(1,3.5);
\draw (1,3.5)--node {\midarrow}(1.5,3);
\draw (1.5,3)--node {\midarrow}(2,2.5);
\draw[red] (3.5,3)--node {\midarrow}(3,2.5);
\draw[red] (3,2.5)--node {\midarrow}(2.5,2);
\draw[red] (2,2.5)--node {\midarrow}(2.5,2);
\draw (2.5,3)--node {\midarrow}(2,2.5);
\draw (3,3.5)--node {\midarrow}(2.5,3);
\draw (3.5,4)--node {\midarrow}(3,3.5);
\end{scope}
\end{tikzpicture}}
\end{minipage}
\begin{minipage}[b]{0.35\linewidth}
\centerline{
\begin{tikzpicture}
\draw[step=.5cm,gray,thin] (-0.5,0) grid (4,5);
\draw[fill] (0,5.3) circle (2pt)--(0.5,5.3) circle (2pt);
\draw[fill] (0.5,5.3) circle (2pt)--(1,5.3) circle (2pt);
\draw[fill] (1,5.3) circle (2pt)--(1.5,5.3) circle (2pt);
\draw[fill] (1.5,5.3) circle (2pt)--(2,5.3) circle (2pt);
\draw[fill] (2,5.3) circle (2pt)--(2.5,5.3) circle (2pt);
\draw[fill] (2.5,5.3) circle (2pt)--(3,5.3)circle (2pt);
\draw[fill] (3,5.3) circle (2pt)--(3.5,5.3)circle (2pt);
\draw [fill] (0,4.5) circle (1.5pt);
\draw [fill,red] (0,3.5) circle (1.5pt);
\draw [fill] (0.5,4) circle (1.5pt);
\draw [fill,red] (0.5,3) circle (1.5pt);
\draw [fill] (1,3.5) circle (1.5pt);
\draw [fill,red] (1,2.5) circle (1.5pt);
\draw [fill] (1.5,3) circle (1.5pt);
\draw [fill,red] (1.5,2) circle (1.5pt);
\draw [fill] (2,2.5) circle (1.5pt);
\draw [fill,red] (2,1.5) circle (1.5pt);
\draw [fill] (2.5,3) circle (1.5pt);
\draw [fill,red] (2.5,2) circle (1.5pt);
\draw [fill] (3,3.5) circle (1.5pt);
\draw [fill,red] (3,2.5) circle (1.5pt);
\draw [fill] (3.5,4) circle (1.5pt);
\draw [fill,red] (3.5,3) circle (1.5pt);
\node[above] at (2,3) {$\bf{\mathscr{L}^{Q}_{5,-3}}$};
\node[red,below] at (2,1.5) {$\bf{\mathscr{L}^{Q_1}_{5,-5}}$};
\node[above] at (0,5.3) {$1$};
\node[above] at  (0.5,5.3) {$2$};
\node[above] at  (1,5.3)   {$3$};
\node[above] at  (1.5,5.3) {$4$};
\node[above] at  (2,5.3)   {$5$};
\node[above] at  (2.5,5.3) {$6$};
\node[above] at  (3,5.3)   {$7$};
\node[above] at  (3.5,5.3) {$8$};
\node [left] at (-0.5,5) {$0$};
\node [left] at (-0.5,4.5) {$-1$};
\node [left] at (-0.5,4) {$-2$};
\node [left] at (-0.5,3.5) {$-3$};
\node [left] at (-0.5,3) {$-4$};
\node [left] at (-0.5,2.5) {$-5$};
\node [left] at (-0.5,2) {$-6$};
\node [left] at (-0.5,1.5) {$-7$};
\node [left] at (-0.5,1) {$-8$};
\node [left] at (-0.5,0.5) {$-9$};
\node [left] at (-0.5,0) {$-10$};
\node [left] at (2,-0.6) {$(\text{c})$};
\begin{scope}[thin, every node/.style={sloped,allow upside down}]
\draw (0,4.5)--node {\midarrow}(0.5,4);
\draw (0.5,4)--node {\midarrow}(1,3.5);
\draw (1,3.5)--node {\midarrow}(1.5,3);
\draw (1.5,3)--node {\midarrow}(2,2.5);
\draw[red] (0,3.5)--node {\midarrow}(0.5,3);
\draw[red] (0.5,3)--node {\midarrow}(1,2.5);
\draw[red] (1,2.5)--node {\midarrow}(1.5,2);
\draw[red] (3.5,3)--node {\midarrow}(3,2.5);
\draw[red] (1.5,2)--node {\midarrow}(2,1.5);
\draw[red] (2.5,2)--node {\midarrow}(2,1.5);
\draw[red] (3,2.5)--node {\midarrow}(2.5,2);
\draw (2.5,3)--node {\midarrow}(2,2.5);
\draw (3,3.5)--node {\midarrow}(2.5,3);
\draw (3.5,4)--node {\midarrow}(3,3.5);
\end{scope}
\end{tikzpicture}}
\end{minipage}}
\caption{In the original quiver of type $A_8$: $(a)$ $\mathscr{L}_{4,-4}=\tau_l(\mathscr{L}_{5,-3})$; $(b)$ $\mathscr{L}_{6,-4}=\tau_r(\mathscr{L}_{5,-3})$; $(c)$ $\mathscr{L}_{5,-5}=\tau(\mathscr{L}_{5,-3})$.}\label{A_8}
\end{figure}
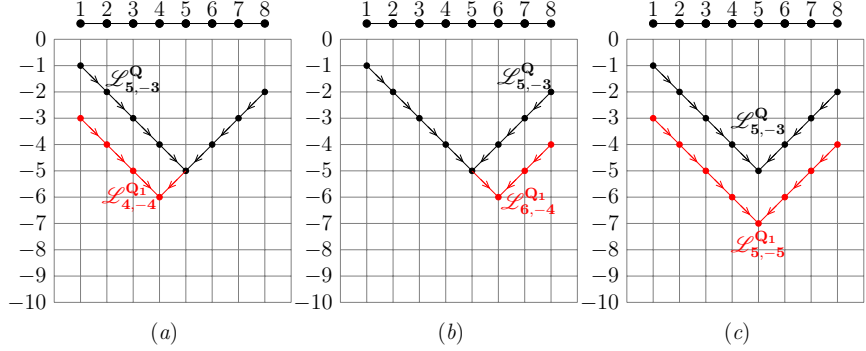

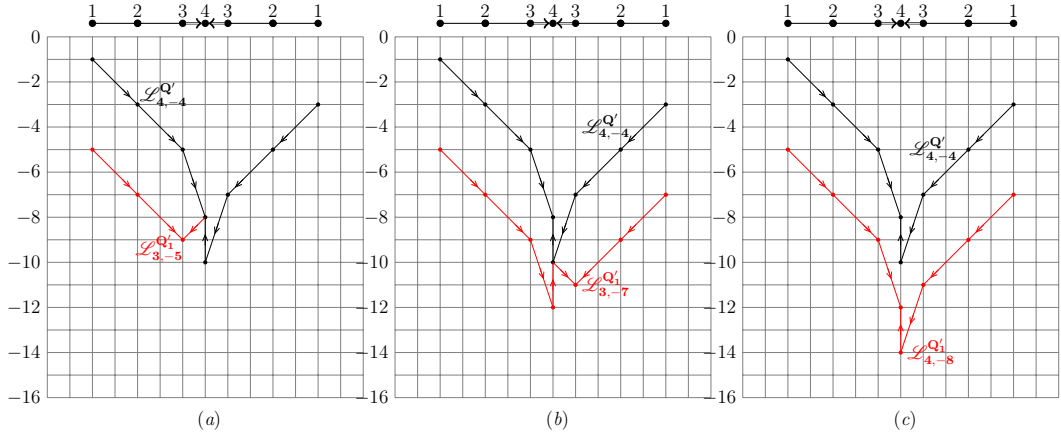
\begin{figure}[H]
\resizebox{.6\width}{.6\height}{
\begin{minipage}[b]{0.5\linewidth}
\centerline{
\begin{tikzpicture}
\draw[step=.5cm,gray,thin] (-1,0) grid (6,8);
\draw[fill] (0,8.3) circle (2pt)--(1,8.3) circle (2pt);
\draw[fill] (1,8.3) circle (2pt)--(2,8.3) circle (2pt);
\draw[fill] (3,8.3) circle (2pt)--(4,8.3) circle (2pt);
\draw[fill] (4,8.3) circle (2pt)--(5,8.3) circle (2pt);
\draw[fill] (2.5,8.3) circle (2pt);
\node [above] at (0,8.3) {$1$};
\node [above] at (1,8.3) {$2$};
\node [above] at (2,8.3) {$3$};
\node [above] at (2.5,8.3) {$4$};
\node [above] at (3,8.3) {$3$};
\node [above] at (4,8.3) {$2$};
\node [above] at (5,8.3) {$1$};
\draw [double,->] (2.075,8.3)--(2.425,8.3);
\draw [double,->] (2.925,8.3)--(2.575,8.3);
\node [left] at (-1,8) {$0$};
\node [left] at (-1,7) {$-2$};
\node [left] at (-1,6) {$-4$};
\node [left] at (-1,5) {$-6$};
\node [left] at (-1,4) {$-8$};
\node [left] at (-1,3) {$-10$};
\node [left] at (-1,2) {$-12$};
\node [left] at (-1,1) {$-14$};
\node [left] at (-1,0) {$-16$};
\node [left] at (3,-0.5) {$(\text{a})$};
\draw[fill] (0,7.5) circle (1pt);
\draw[fill,red] (0,5.5) circle (1pt);
\draw[fill] (1,6.5) circle (1pt);
\draw[fill,red] (1,4.5) circle (1pt);
\draw[fill] (2,5.5) circle (1pt);
\draw[fill,red] (2,3.5) circle (1pt);
\draw[fill] (2.5,4) circle (1pt);
\draw[fill] (2.5,3) circle (1pt);
\draw[fill] (3,4.5) circle (1pt);
\draw[fill] (4,5.5) circle (1pt);
\draw[fill] (5,6.5) circle (1pt);
\node[right] at (0.9,6.7) {$\bf{\mathscr{L}^{Q'}_{4,-4}}$};
\node[red,left] at (2.15,3.3) {$\bf{\mathscr{L}^{Q'_1}_{3,-5}}$};
\begin{scope}[thin, every node/.style={sloped,allow upside down}]
\draw (0,7.5)--node{\midarrow}(1.5,6);
\draw[red] (0,5.5)--node{\midarrow}(1.5,4);
\draw (1.5,6)--node{\midarrow}(2,5.5);
\draw[red] (1.5,4)--node{\midarrow}(2,3.5);
\draw (2,5.5)--node{\midarrow}(2.5,4);
\draw[red] (2.5,4)--node{\midarrow}(2,3.5);
\draw (2.5,3)--node{\midarrow}(2.5,4);
\draw (3,4.5)--node{\midarrow}(2.5,3);
\draw (3.5,5)--node{\midarrow}(3,4.5);
\draw (5,6.5)--node{\midarrow}(3.5,5);
\end{scope}
\end{tikzpicture}}
\end{minipage}
\begin{minipage}[b]{0.5\linewidth}
\centerline{
\begin{tikzpicture}
\draw[step=.5cm,gray,thin] (-1,0) grid (6,8);
\draw[fill] (0,8.3) circle (2pt)--(1,8.3) circle (2pt);
\draw[fill] (1,8.3) circle (2pt)--(2,8.3) circle (2pt);
\draw[fill] (3,8.3) circle (2pt)--(4,8.3) circle (2pt);
\draw[fill] (4,8.3) circle (2pt)--(5,8.3) circle (2pt);
\draw[fill] (2.5,8.3) circle (2pt);
\node [above] at (0,8.3) {$1$};
\node [above] at (1,8.3) {$2$};
\node [above] at (2,8.3) {$3$};
\node [above] at (2.5,8.3) {$4$};
\node [above] at (3,8.3) {$3$};
\node [above] at (4,8.3) {$2$};
\node [above] at (5,8.3) {$1$};
\draw [double,->] (2.075,8.3)--(2.425,8.3);
\draw [double,->] (2.925,8.3)--(2.575,8.3);
\node [left] at (-1,8) {$0$};
\node [left] at (-1,7) {$-2$};
\node [left] at (-1,6) {$-4$};
\node [left] at (-1,5) {$-6$};
\node [left] at (-1,4) {$-8$};
\node [left] at (-1,3) {$-10$};
\node [left] at (-1,2) {$-12$};
\node [left] at (-1,1) {$-14$};
\node [left] at (-1,0) {$-16$};
\node [left] at (3,-0.5) {$(\text{b})$};
\draw[fill] (0,7.5) circle (1pt);
\draw[fill,red] (0,5.5) circle (1pt);
\draw[fill] (1,6.5) circle (1pt);
\draw[fill,red] (1,4.5) circle (1pt);
\draw[fill] (2,5.5) circle (1pt);
\draw[fill,red] (2,3.5) circle (1pt);
\draw[fill] (2.5,4) circle (1pt);
\draw[fill] (2.5,3) circle (1pt);
\draw[fill,red] (2.5,2) circle (1pt);
\draw[fill] (3,4.5) circle (1pt);
\draw[fill,red] (3,2.5) circle (1pt);
\draw[fill] (4,5.5) circle (1pt);
\draw[fill,red] (4,3.5) circle (1pt);
\draw[fill] (5,6.5) circle (1pt);
\draw[fill,red] (5,4.5) circle (1pt);
\node[right] at (3,6) {$\bf{\mathscr{L}^{Q'}_{4,-4}}$};
\node[red,right] at (3,2.5) {$\bf{\mathscr{L}^{Q'_1}_{3,-7}}$};
\begin{scope}[thin, every node/.style={sloped,allow upside down}]
\draw (0,7.5)--node {\midarrow}(1.5,6);
\draw[red] (0,5.5)--node {\midarrow}(1.5,4);
\draw (1.5,6)--node {\midarrow}(2,5.5);
\draw[red] (1.5,4)--node {\midarrow}(2,3.5);
\draw (2,5.5)--node {\midarrow}(2.5,4);
\draw (2.5,3)--node {\midarrow}(2.5,4);
\draw[red] (2.5,2)--node {\midarrow}(2.5,3);
\draw[red] (2,3.5)--node {\midarrow}(2.5,2);
\draw[red] (2.5,3)--node {\midarrow}(3,2.5);
\draw[red] (5,4.5)--node {\midarrow}(3.5,3);
\draw[red] (3.5,3)--node {\midarrow}(3,2.5);
\draw (3,4.5)--node {\midarrow}(2.5,3);
\draw (3.5,5)--node {\midarrow}(3,4.5);
\draw (5,6.5)--node {\midarrow}(3.5,5);
\end{scope}
\end{tikzpicture}}
\end{minipage}
\begin{minipage}[b]{0.5\linewidth}
\centerline{
\begin{tikzpicture}
\draw[step=.5cm,gray,thin] (-1,0) grid (6,8);
\draw[fill] (0,8.3) circle (2pt)--(1,8.3) circle (2pt);
\draw[fill] (1,8.3) circle (2pt)--(2,8.3) circle (2pt);
\draw[fill] (3,8.3) circle (2pt)--(4,8.3) circle (2pt);
\draw[fill] (4,8.3) circle (2pt)--(5,8.3) circle (2pt);
\draw[fill] (2.5,8.3) circle (2pt);
\node [above] at (0,8.3) {$1$};
\node [above] at (1,8.3) {$2$};
\node [above] at (2,8.3) {$3$};
\node [above] at (2.5,8.3) {$4$};
\node [above] at (3,8.3) {$3$};
\node [above] at (4,8.3) {$2$};
\node [above] at (5,8.3) {$1$};
\draw [double,->] (2.075,8.3)--(2.425,8.3);
\draw [double,->] (2.925,8.3)--(2.575,8.3);
\node [left] at (-1,8) {$0$};
\node [left] at (-1,7) {$-2$};
\node [left] at (-1,6) {$-4$};
\node [left] at (-1,5) {$-6$};
\node [left] at (-1,4) {$-8$};
\node [left] at (-1,3) {$-10$};
\node [left] at (-1,2) {$-12$};
\node [left] at (-1,1) {$-14$};
\node [left] at (-1,0) {$-16$};
\node [left] at (3,-0.5) {$(\text{c})$};
\draw[fill] (0,7.5) circle (1pt);
\draw[fill,red] (0,5.5) circle (1pt);
\draw[fill] (1,6.5) circle (1pt);
\draw[fill,red] (1,4.5) circle (1pt);
\draw[fill] (2,5.5) circle (1pt);
\draw[fill,red] (2,3.5) circle (1pt);
\draw[fill] (2.5,4) circle (1pt);
\draw[fill] (2.5,3) circle (1pt);
\draw[fill,red] (2.5,2) circle (1pt);
\draw[fill,red] (2.5,1) circle (1pt);
\draw[fill] (3,4.5) circle (1pt);
\draw[fill,red] (3,2.5) circle (1pt);
\draw[fill] (4,5.5) circle (1pt);
\draw[fill,red] (4,3.5) circle (1pt);
\draw[fill] (5,6.5) circle (1pt);
\draw[fill,red] (5,4.5) circle (1pt);
\node[right] at (2.55,5.5) {$\bf{\mathscr{L}^{Q'}_{4,-4}}$};
\node[red,right] at (2.5,1) {$\bf{\mathscr{L}^{Q'_1}_{4,-8}}$};
\begin{scope}[thin, every node/.style={sloped,allow upside down}]
\draw (0,7.5)--node {\midarrow}(1.5,6);
\draw[red] (0,5.5)--node {\midarrow}(1.5,4);
\draw (1.5,6)--node {\midarrow}(2,5.5);
\draw[red] (1.5,4)--node {\midarrow}(2,3.5);
\draw (2,5.5)--node {\midarrow}(2.5,4);
\draw (2.5,3)--node {\midarrow}(2.5,4);
\draw[red] (2.5,1)--node {\midarrow}(2.5,2);
\draw[red] (2,3.5)--node {\midarrow}(2.5,2);
\draw[red] (3,2.5)--node {\midarrow}(2.5,1);
\draw[red] (5,4.5)--node {\midarrow}(3.5,3);
\draw[red] (3.5,3)--node {\midarrow}(3,2.5);
\draw (3,4.5)--node {\midarrow}(2.5,3);
\draw (3.5,5)--node {\midarrow}(3,4.5);
\draw (5,6.5)--node {\midarrow}(3.5,5);
\end{scope}
\end{tikzpicture}}
\end{minipage}}
\caption{In the original quiver of type $B_4$: $(a)$ $\mathscr{L}_{3, -5}=\tau_l(\mathscr{L}_{4, -4})$; $(b)$ $\mathscr{L}_{3, -7}=\tau_r(\mathscr{L}_{4, -4})$; $(c)$ $\mathscr{L}_{4, -8}=\tau(\mathscr{L}_{4, -4})$.}\label{B_4}
\end{figure}
\end{example}

\subsection{Mutation sequences of Kirillov--Reshetikhin modules}\label{mutation sequences of Kirillov--Reshetikhin modules}
In \cite{HL13}, Hernandez and Leclerc defined a sequence of mutations for every Kirillov--Reshetikhin module whose highest weight monomial $m$ satisfies the property: $(i, t) \in W^-$ for every factor $Y_{i,q^t}$ in $m$. We recall the mutation sequences for Kirillov--Reshetikhin modules introduced in \cite{HL13}.

Let $\widetilde{Q}$ (respectively, $\widetilde{Q}'$) be a quiver which is mutation equivalent to $Q$ (respectively, $Q'$) defined in Section \ref{definition of cluster algebra} and any two vertices in $\widetilde{Q}$ (respectively, $\widetilde{Q}'$) have different labels. By saying that we mutate $C_{i,t}$ of $\widetilde{Q}$ (respectively, $\widetilde{Q}'$), we mean that we mutate at the vertex of $\widetilde{Q}$ (respectively, $\widetilde{Q}'$) which has the label $(i, t)$ in the $i$th column and so on until the vertex at infinity in the $i$th column.

Let $\seq_j$, $m_1 \leq j \leq m_2$, be mutation sequences. We use
\[
\overleftarrow{\prod_{j=m_2}^{m_1}}\seq\nolimits_{j} \text{ and } \left((\seq\nolimits_j)_{m_1\leq j \leq m_2} \text{ or } \prod_{j=m_1}^{m_2} \seq\nolimits_j\right)
\]
to denote the mutation sequences
\[
\seq\nolimits_{m_2}, \seq\nolimits_{m_2-1}, \ldots, \seq\nolimits_{m_1} \text{ and } \seq\nolimits_{m_1}, \seq\nolimits_{m_1+1}, \ldots, \seq\nolimits_{m_2},
\]
respectively.

Consider the Kirillov--Reshetikhin module $\mathcal{S}^{(t)}_{k^{(i)}}$, $t \leq 0$, $k \in \mathbb{Z}_{\geq 1}$, $i \in I$. In the case of type $A_n$ (\cite{HL13}, Section 3), we mutate
\begin{align*}
\prod_{\ell=0}^{j-1}\left(\prod_{r=1}^{\lfloor \frac{n}{2} \rfloor}C_{2r, -2\ell}\right) \left(\prod_{r=1}^{\lceil \frac{n}{2} \rceil}C_{2r-1, -2\ell-1}\right)
\end{align*}
starting from the quiver $Q$, where $j$ is defined by the formula
\begin{align*}
t=
\begin{cases}
-2k-2j+2, & i \in 2\mathbb{Z} \cap I, \\
-2k-2j+1, & i \in (2\mathbb{Z}+1) \cap I.
\end{cases}
\end{align*}
We use $Q_0$ to denote the current quiver. Then we obtain the Kirillov--Reshetikhin module $\mathcal{S}^{(t)}_{k^{(i)}}$, $t \leq 0$, $k \in \mathbb{Z}_{\geq 1}$, $i \in I$, at the vertex $(i, t)$ of $Q_0$.

In the case of type $B_n$ (\cite{HL13}, Section 3), let
\begin{gather}
\begin{align*}
KR(n,\ell)=C_{n, -4\ell} \left(\prod_{r=0}^{\lfloor \frac{n}{2}-1 \rfloor}C_{n-1-2r, -4\ell-1}\right) \left(\prod_{r=0}^{\lceil \frac{n}{2}-2 \rceil} C_{n-2-2r, -4\ell-3}\right) C_{n, -4\ell-2} \left(\prod_{r=0}^{\lceil \frac{n}{2}-2 \rceil} C_{n-2-2r, -4\ell-1}\right)\left(\prod_{r=0}^{\lfloor \frac{n}{2}-1 \rfloor}C_{n-1-2r, -4\ell-3}\right).
\end{align*}
\end{gather}
When $i \neq n$, $t=-4k-2j+3$, we mutate $\prod\limits_{\ell=0}^{\lfloor \frac{j}{2}-1 \rfloor }KR(n,\ell)$ starting from the quiver $Q'$. When $i =n$, $t=-2k-2j+4$, we mutate $\prod\limits_{\ell=0}^{\lceil \frac{j}{2}-2 \rceil }KR(n,\ell)$ if $j$ is odd, and mutate
\[
\left(\prod_{\ell=0}^{\lceil \frac{j}{2}-2 \rceil }KR(n,\ell)\right)\left(C_{n, -2j+4} \left(\prod_{r=0}^{\lfloor \frac{n}{2}-1 \rfloor}C_{n-1-2r, -2j+3}\right) \left(\prod_{r=0}^{\lceil \frac{n}{2}-2 \rceil} C_{n-2-2r, -2j+1}\right)\right)
\]
if $j$ is even, starting from the quiver $Q'$. We use $Q_0'$ to denote the current quiver. Then we obtain the Kirillov--Reshetikhin module $\mathcal{S}^{(t)}_{k^{(i)}}$, $t \leq 0$, $k \in \mathbb{Z}_{\geq 1}$, $i \in I$, at the vertex $(i, t)$ of $Q_0'$.

\subsection{Mutation sequences for prime snake modules of types $A_{n}$ and $B_{n}$}\label{mutation sequences for prime snake modules}
Let $\mathcal{S}$ be a prime snake module with the highest $l$-weight monomial $S$. Then $S$ can be written as
\begin{align*}
S = S^{(t)}_{k_{1}^{(i_{1},j_{1})}, k_{2}^{(i_{2},j_{2})}, \ldots,k_{m-1}^{(i_{m-1},j_{m-1})},k_{m}^{(i_{m})}},
\end{align*}
where $m \geq 1$, $j_\ell \geq 0$, $1 \leq \ell \leq m-1$, if $j_{\ell}=0$, then $i_{\ell} \neq i_{\ell+1}$, $\scalemath{0.92}{k_1, \ldots, k_m \in \mathbb{Z}_{\geq 1},  t \in \mathbb{Z}}$.

By Definition \ref{define fundamental segment}, we have
\[
\mathcal{FS}(S)=FS_1(S)\cup FS_2(S) \cup FS_3(S),
\]
where $FS_1=\{M_1=S_{{k_m}^{(i_m)}}^{(t_m)}\}$, $FS_2(S) \cup FS_3(S)=\{M_2, \ldots, M_{q}\}$. We reorder the elements in $FS_2(S) \cup FS_3(S)$ such that the distinguished factor $(l_p)_{s_p}$ of $M_p$ and the distinguished factor $(l_{p+1})_{s_{p+1}}$ of $M_{p+1}$ satisfy $s_p > s_{p+1}$, $2 \leq p \leq q-1$.

Let $Q_1,\ Q_2,\ \ldots,\ Q_h$ be quivers in a mutation sequence. Let $(i_\ell,s_\ell) \in V(Q_\ell)$, $1 \leq \ell \leq h$.
For simplicity, we write $\mathscr{L}_{i_\ell,s_\ell}$ for $\mathscr{L}^{Q_\ell}_{i_\ell,s_\ell}$ and write $C_{i_\ell,s_\ell}$ for $C^{Q_\ell}_{i_\ell,s_\ell}$.

Let $M_p \in \mathcal{FS}(S)$, $1 \leq p \leq q$, and let $(l_p)_{s_p}$ be the distinguished factor of $M_p$. Using the mutation sequence defined in Section \ref{mutation sequences of Kirillov--Reshetikhin modules} starting from the initial quiver $Q$ in type $A_n$ (respectively, $Q'$ in type $B_n$) defined in Section \ref{definition of cluster algebra}, we can obtain a quiver $Q_0$ (respectively, $Q_0'$) and obtain the module $L(M_1) = L(S_{{k_m}^{(i_m)}}^{(t_m)})$ at the vertex $(i_m, t_m)$ of $Q_0$ (respectively, $Q_0'$).

In the following, we define mutation sequences $\seq_2$, $\seq_3$, $\ldots$, $\seq_q$ starting from the quiver $Q_0$ (respectively, $Q'_0$) of type $A_n$ (respectively, $B_n$) such that after we mutate $\seq_1$, $\seq_2$, $\ldots$, $\seq_q$, we obtain the prime snake module $\mathcal{S}=L(S)$ at the vertex $(i_1,t)$.

The following is the case of type $A_n$.
\begin{enumerate}[(1)]
\item Suppose that $M_p \in FS_2(S)$. Then there are some $\ell$, $r$ such that
\begin{align*}
M_p=S^{(t_\ell)}_{k_{\ell}^{(i_\ell)}, k_{\ell+1}^{(i_{\ell+1})},\ldots, k_{\ell+r-1}^{(i_{\ell+r-1})}, 1^{(i_{\ell+r})}}.
\end{align*}
If the sequence $(i_u)_{\ell \leq u \leq \ell+r}$ is in decreasing order (respectively, in increasing order), then we mutate $(\mathfrak{r}(\tau^{h}_r(\mathscr{L}_{l_p,s_p})))_{0\leq h \leq n-l_p-1}$
(respectively, $\scalemath{0.88}{(\mathfrak{l}(\tau^{h}_l(\mathscr{L}_{l_p,s_p})))_{0\leq h \leq l_p-2}}$). If $r \geq 2$, we continue to mutate
    \begin{align*}
    \overleftarrow{\prod_{u=\ell+r-1}^{\ell+1}} \left\{\prod_{j=1+\sum_{i=u+1}^{\ell+r-1}k_{i}}^{\sum_{i=u}^{\ell+r-1}k_i} \left(C_{i_{u}+1, s_{i_{u}+1}-2j}, C_{i_{u}+2, s_{i_{u}+2}-2j}, \ldots, C_{n, s_n-2j}\right)\right\}
    \end{align*}
    \begin{align*}
    \text{ (respectively, } \overleftarrow{\prod_{u=\ell+r-1}^{\ell+1}} \left\{ \prod_{j=1+\sum_{i=u+1}^{\ell+r-1}k_{i}}^{\sum_{i=u}^{\ell+r-1}k_i} \left(C_{i_{u}-1, s_{i_{u}-1}-2j}, C_{i_{u}-2, s_{i_{u}-2}-2j}, \ldots, C_{1, s_1-2j}\right)\right\}),
    \end{align*}
where $i$, $s_i$ satisfy $\mathscr{L}_{i, s_i} = \tau_{r}^{i-l_p}(\mathscr{L}_{l_p, s_p})$, $i_{\ell+r} \leq i \leq n$
(\text{respectively,} $\mathscr{L}_{i, s_i}=\tau_{l}^{l_p-i}(\mathscr{L}_{l_p, s_p})$, $1 \leq i \leq i_{\ell+r})$.

\item Suppose that $M_p\in FS_3(S)$. Then there is some $\ell$ such that $M_p=S^{(t_\ell)}_{k_{\ell}^{(i_\ell,j_\ell)},1^{(i_{\ell+1})}}$.
\begin{itemize}
\item If $i_{\ell} \geq l_p$, then we mutate
     \[
     (\mathfrak{r}(\tau^{h}_r(\mathscr{L}_{l_p,s_p})))_{0\leq h \leq i_\ell-l_p-1},\quad (\tau^{h}(\tau^{i_\ell-l_p}_r(\mathscr{L}_{l_p,s_p})))_{0 \leq h \leq j_\ell-1}.
     \]
\item If $i_{\ell} \leq l_p$, then we mutate
     \[
     (\mathfrak{l}(\tau^{h}_l(\mathscr{L}_{l_p,s_p})))_{0\leq h \leq l_p-i_\ell-1},\quad  (\tau^{h}(\tau^{l_p-i_\ell}_l(\mathscr{L}_{l_p,s_p})))_{0\leq h \leq j_\ell-1}.
     \]
\end{itemize}
\end{enumerate}

The following is the case of type $B_{n}$.
\begin{enumerate}[(1)]
\item Suppose that $M_p \in FS_2(S)$. Then there are some $\ell$, $r$ such that
\begin{align*}
M_p=S^{(t_\ell)}_{k_{\ell}^{(i_\ell)}, k_{\ell+1}^{(i_{\ell+1})}, \ldots, k_{\ell+r-1}^{(i_{\ell+r-1})}, 1^{(i_{\ell+r})}}.
\end{align*}
Suppose that the sequence $(i_u)_{\ell \leq u \leq \ell+r}$ is in increasing order.
If $l_p\neq n$, $2n+s_p-2l_p \equiv 1 \pmod 4$ or $l_p=n$, $s_p \equiv 0 \pmod 4$ (respectively, $l_p \neq n$, $2n+s_p-2l_p \equiv 3 \pmod 4$ or $l_p=n$, $s_p \equiv 2 \pmod 4$), then we mutate
    \[
    (\mathfrak{l}(\tau^{h}_l(\mathscr{L}_{l_p,s_p})))_{0\leq h \leq l_p-2} \quad (\text{respectively, } (\mathfrak{r}(\tau^{h}_r(\mathscr{L}_{l_p,s_p})))_{0\leq h \leq l_p-2}).
    \]
If $r \geq 2$, then we continue to mutate
\begin{align*}
\overleftarrow{\prod_{u=\ell+r-1}^{\ell+1}} \left\{ \prod_{j=1+\sum_{i=u+1}^{\ell+r-1}k_{i} }^{\sum_{i=u}^{\ell+r-1}k_i} \left(C_{i_{u}-1, s_{i_{u}-1}-4j}, C_{i_{u}-2, s_{i_{u}-2}-4j}, \ldots, C_{1, s_1-4j}\right)\right\},
\end{align*}
where $i$, $s_i$ satisfy $\mathscr{L}_{i, s_i}=\tau_{l}^{l_p-i}(\mathscr{L}_{l_p, s_p})$, $1 \leq i \leq i_{\ell+r}$.

\item Suppose that $M_p \in FS_2(S)$. Then there are some $\ell$, $r$ such that
\begin{align*}
M_p=S^{(t_\ell)}_{k_{\ell}^{(i_\ell)}, k_{\ell+1}^{(i_{\ell+1})}, \ldots, k_{\ell+r-1}^{(i_{\ell+r-1})}, 1^{(i_{\ell+r})}}.
\end{align*}
Suppose that the sequence $(i_u)_{\ell \leq u \leq \ell+r}$ is in decreasing order.
\begin{itemize}
\item If $2n+s_p-2l_p \equiv 1 \pmod 4$, then we mutate $(\mathfrak{r}(\tau^{h}_r(\mathscr{L}_{l_p,s_p})))_{0\leq h \leq n-l_p-1}$. If $r \geq 2$, then we continue to mutate
\begin{gather}
\begin{align*}
\overleftarrow{\prod_{u=\ell+r-1}^{\ell+1}} \left\{ \prod_{j=1+\sum_{i=u+1}^{\ell+r-1}k_{i} }^{\sum_{i=u}^{\ell+r-1}k_i} \left(\mathfrak{r}(\mathscr{L}_{n, s_n-4j+4}), C_{i_{u}+1, s_{i_{u}+1}-4j}, C_{i_{u}+2, s_{i_{u}+2}-4j}, \ldots, C_{n, s_n-4j}\right)\right\},
\end{align*}
\end{gather}
where $i$, $s_i$ satisfy $\mathscr{L}_{i, s_i}=\tau_{r}^{i-l_p}(\mathscr{L}_{l_p, s_p})$, $i_{\ell+r} \leq i \leq n$.

\item If $2n+s_p-2l_p \equiv 3 \pmod 4$, then we mutate $(\mathfrak{l}(\tau^{h}_l(\mathscr{L}_{l_p,s_p})))_{0\leq h \leq n-l_p-1}$. If $r \geq 2$, then we continue to mutate
\begin{gather}
\begin{align*}
\overleftarrow{\prod_{u=\ell+r-1}^{\ell+1}} \left\{ \prod_{j=1+\sum_{i=u+1}^{\ell+r-1}k_{i} }^{\sum_{i=u}^{\ell+r-1}k_i} \left(\mathfrak{l}(\mathscr{L}_{n, s_n-4j+4}), C_{i_{u}+1, s_{i_{u}+1}-4j}, C_{i_{u}+2, s_{i_{u}+2}-4j}, \ldots, C_{n, s_n-4j}\right)\right\},
\end{align*}
\end{gather}
where $i$, $s_i$ satisfy $\mathscr{L}_{i, s_i}=\tau_{l}^{i-l_p}(\mathscr{L}_{l_p, s_p})$, $i_{\ell+r} \leq i \leq n$.
\end{itemize}

\item Suppose that $M_p \in FS_3(S)$. Then there is some $\ell$ such that $M_p=S^{(t_\ell)}_{k_{\ell}^{(i_\ell,j_\ell)},1^{(i_{\ell+1})}}$.
\begin{itemize}
\item If $l_p \neq n$, $2n+s_p-2l_p \equiv 1 \pmod 4$, $i_{\ell} \geq l_p$ or $l_p \neq n$, $2n+s_p-2l_p \equiv 3 \pmod 4$, $i_{\ell} \leq l_p$ or $l_p=n$, $s_p \equiv 2 \pmod 4$, then we mutate
    \[
     (\mathfrak{r}(\tau^{h}_r(\mathscr{L}_{l_p,s_p})))_{0\leq h \leq |i_\ell-l_p|-1}, \quad (\tau^{h}(\tau^{|l_p-i_\ell|}_r(\mathscr{L}_{l_p,s_p})))_{0 \leq h \leq j_\ell-1}.
     \]
\item If $l_p \neq n$, $2n+s_p-2l_p \equiv 1 \pmod 4$, $i_{\ell} \leq l_p$ or $l_p \neq n$, $2n+s_p-2l_p \equiv 3 \pmod 4$, $i_{\ell} \geq l_p$ or $l_p=n$, $s_p \equiv 0 \pmod 4$, then we mutate
    \[
    (\mathfrak{l}(\tau^{h}_l(\mathscr{L}_{l_p,s_p})))_{0\leq h \leq |i_\ell-l_p|-1}, \quad (\tau^{h}(\tau^{|l_p-i_\ell|}_l(\mathscr{L}_{l_p,s_p})))_{0 \leq h \leq j_\ell-1}.
    \]
\end{itemize}
\end{enumerate}

\begin{remark}
Minimal affinizations are modules $\mathcal{S}^{(t)}_{k_{1}^{(i_1)}, k_{2}^{(i_{2})}, \ldots, k_{r}^{(i_{r})}}$ which satisfy $i_1 < \cdots < i_r$ or $i_1 > \cdots > i_r$. The mutation sequences above in the case (1) in type $A_n$ (respectively, the cases (1), (2) in type $B_n$) are mutation sequences for minimal affinizations which satisfy $i_1 < \cdots < i_r$ or $i_1 > \cdots > i_r$. These mutation sequences are defined in the same cluster algebra. In \cite{ZDLL16}, the mutation sequences for minimal affinizations which satisfy $i_1 < \cdots < i_r$ are defined in a cluster algebra $\mathscr{A}$ and the mutation sequences for minimal affinizations which satisfy $i_1 > \cdots > i_r$ are defined in another cluster algebra $\widetilde{\mathscr{A}}$ which is dual to $\mathscr{A}$.
\end{remark}

\subsection{The equations in the $S$-system of type $A_n$ (respectively, $B_{n}$) correspond to mutations in the cluster algebra $\mathscr{A}$ (respectively, $\mathscr{A}'$)}
In this section, we give the relation between prime snake modules and cluster variables.

Let $\bf{S}$ be the set of prime snake modules and
\begin{gather}
\begin{align*}
\mathfrak{s}=\{s^{(t)}_{k_{1}^{(i_{1},j_1)}, k_{2}^{(i_{2},j_2)}, \ldots, k_{m-1}^{(i_{m-1},j_{m-1})}, k_{m}^{(i_{m})}}: j_\ell \geq 0,\ 1 \leq \ell \leq m-1,\ k_1, \ldots, k_m \in \mathbb{Z}_{\geq 1},\ t \in \mathbb{Z}\}.
\end{align*}
\end{gather}

There is a bijection from $\bf{S}$ to $\mathfrak{s}$ given by
\begin{equation*}
\begin{split}
\psi: \bf{S}& \longrightarrow \mathfrak{s} \\
\mathcal{S}^{(t)}_{k_{1}^{(i_{1},j_1)}, k_{2}^{(i_{2},j_2)}, k_{3}^{(i_{3},j_{3})}, \ldots, k_{m}^{(i_{m})}}& \longmapsto s^{(t)}_{k_{1}^{(i_{1},j_1)}, k_{2}^{(i_{2},j_2)}, k_{3}^{(i_{3},j_{3})}, \ldots, k_{m}^{(i_{m})}}.
\end{split}
\end{equation*}
We apply $\psi$ to the equations $[\mathcal{S}_1][\mathcal{S}_2]=[\mathcal{S}_3][\mathcal{S}_4]+[\mathcal{S}_5][\mathcal{S}_6]$ in the $S$-system for type $A_n$ (respectively, $B_n$). Then we have a new system of equations:
\begin{align}\label{mutation equations 1}
s_1s_2=s_3s_4+s_5s_6,
\end{align}
where $s_i= \psi(\mathcal{S}_i)$, $1 \leq i \leq 6$.

On the other hand, in type $A_n$ (respectively, $B_n$), for every mutation sequence in Section \ref{mutation sequences for prime snake modules}, at each step, the variable at a vertex we mutate is $\psi(\mathcal{S}_1)$, where $\mathcal{S}_1$ is one of the modules in Table \ref{definition S_1 in type A_{n}} (respectively, Table \ref{definition S_1 in type B_{n}}). For every $\mathcal{S}_1$ in Table \ref{definition S_1 in type A_{n}} (respectively, Table \ref{definition S_1 in type B_{n}}), there is a corresponding $\mathcal{S}_2$ in (\ref{S_2}) and corresponding $\mathcal{S}_3$, $\mathcal{S}_4$, $\mathcal{S}_5$, $\mathcal{S}_6$ defined in Section \ref{definition of S-system}.

After a mutation at a vertex with variable $x = \psi(\mathcal{S}_1)=s_1$, let the new variable at this vertex be $x' = \psi(\mathcal{S}_2)=s_2$. Then the exchange relation is $x x' = s_1 s_2 = s_3s_4 + s_5 s_6$. Therefore we have the following theorem.

\begin{theorem} \label{set of exchange relations corresponds to S systems}
The set of exchange relations of the mutations of type $A_n$ (respectively, $B_n$) in Section \ref{mutation sequences for prime snake modules} correpsonds to the $S$-system in type $A_n$ (respectively, $B_n$).
\end{theorem}
Indeed, originally the equations in $S$-systems are found by analyzing the mutation sequences in Section \ref{mutation sequences for prime snake modules}.

The following theorem follows from Theorems \ref{prime snake modules}, \ref{real snake modules}, \ref{path description of q-characters}, and \ref{set of exchange relations corresponds to S systems}.
\begin{theorem}
The Hernandez-Leclerc conjecture (Conjecture \ref{Hernandez-Leclerc conjecture}) is true for prime snake modules of types $A_{n}$ and $B_{n}$: prime snake modules of type $A_{n}$ (respectively, $B_{n}$) correspond to cluster variables in the corresponding cluster algebra.
\end{theorem}

\section{Examples of mutation sequences for some prime snake modules}\label{examples of mutation sequences for some snake modules}
In this section, we give some examples of mutation sequences for some prime snake modules.
\begin{example}\label{example 1}
Suppose that $U_q \widehat{\mathfrak{g}}$ is of type $A_{5}$. (1) Let $S=2_{-12}4_{-8}5_{-5}5_{-3}4_{0}$. By Definition \ref{define fundamental segment},
\begin{align*}
\mathcal{FS}(S)=\{4_0, \ 5_{-5}5_{-3}4_{0}, \ 2_{-12}4_{-8}5_{-5} \}.
\end{align*}
The set of distinguished factors of $S$ is $\{4_0,\ 5_{-5}\}$. The mutation sequence for $S$ is
\begin{align}\label{mutation sequence in example 6.1}
\mathfrak{r}(\mathscr{L}_{4,0}), \ \mathfrak{l}(\mathscr{L}_{5,-5}), \ \mathfrak{l}(\tau_l(\mathscr{L}_{5,-5})), \ \mathfrak{l}(\tau^{2}_l(\mathscr{L}_{5,-5})), \ \mathfrak{l}(\tau^{3}_l(\mathscr{L}_{5,-5})), \ C_{3,-9}, \ C_{2,-10}, \ C_{1,-11}.
\end{align}
We obtain the prime snake module $\mathcal{S}=L(S)$ at the vertex which has the label $(2,-12)$, see Figure \ref{Q_{1}}.

The initial quivers in this section are the initial quivers in \cite{HL13}. The mutation sequences in this section are similar to the mutation sequences given in \cite{HL13}. In \cite{HL13}, the mutation sequences produce  Kirillov--Reshetikhin modules. In this paper, the mutation sequences defined in Section \ref{mutation sequences for prime snake modules} produce prime snake modules. A very useful device for mutating quivers is due to B. Keller's Java applet \cite{K12}.

\begin{figure}[H]
\centering
\resizebox{.6\width}{.6\height}{
\begin{minipage}[b]{0.8\textwidth}
\centerline{
\begin{xy}
(10,70) *+{s^{(-1)}_{1^{(1)}}} ="0",
(10,60) *+{s^{(-3)}_{2^{(1)}}} ="1",
(10,50) *+{s^{(-5)}_{3^{(1)}}} ="2",
(10,40) *+{s^{(-7)}_{4^{(1)}}} ="3",
(10,30) *+{s^{(-9)}_{5^{(1)}}} ="4",
(10,20) *+{s^{(-11)}_{6^{(1)}}} ="5",
(10,10) *+{s^{(-13)}_{7^{(1)}}} ="6",
(10,0) *+{\vdots} ="7",
(25,75) *+{s^{(0)}_{1^{(2)}}} ="8",
(25,65) *+{s^{(-2)}_{2^{(2)}}} ="9",
(25,55) *+{s^{(-4)}_{3^{(2)}}} ="10",
(25,45) *+{s^{(-6)}_{4^{(2)}}} ="11",
(25,35) *+{s^{(-8)}_{5^{(2)}}} ="12",
(25,25) *+{s^{(-10)}_{6^{(2)}}} ="13",
(25,15) *+{s^{(-12)}_{7^{(2)}}} ="14",
(25,5)  *+{\vdots} ="15",
(40,70) *+{s^{(-1)}_{1^{(3)}}} ="16",
(40,60) *+{s^{(-3)}_{2^{(3)}}} ="17",
(40,50) *+{s^{(-5)}_{3^{(3)}}} ="18",
(40,40) *+{s^{(-7)}_{4^{(3)}}} ="19",
(40,30) *+{s^{(-9)}_{5^{(3)}}} ="20",
(40,20) *+{s^{(-11)}_{6^{(3)}}} ="21",
(40,10) *+{s^{(-13)}_{7^{(3)}}} ="22",
(40,0)  *+{\vdots} ="23",
(55,75) *+{s^{(0)}_{1^{(4)}}} ="24",
(55,65) *+{s^{(-2)}_{2^{(4)}}} ="25",
(55,55) *+{s^{(-4)}_{3^{(4)}}} ="26",
(55,45) *+{s^{(-6)}_{4^{(4)}}} ="27",
(55,35) *+{s^{(-8)}_{5^{(4)}}} ="28",
(55,25) *+{s^{(-10)}_{6^{(4)}}} ="29",
(55,15) *+{s^{(-12)}_{7^{(4)}}} ="30",
(55,5)  *+{\vdots} ="31",
(70,70) *+{s^{(-1)}_{1^{(5)}}} ="32",
(70,60) *+{s^{(-3)}_{2^{(5)}}} ="33",
(70,50) *+{s^{(-5)}_{3^{(5)}}} ="34",
(70,40) *+{s^{(-7)}_{4^{(5)}}} ="35",
(70,30) *+{s^{(-9)}_{5^{(5)}}} ="36",
(70,20) *+{s^{(-11)}_{6^{(5)}}} ="37",
(70,10) *+{s^{(-13)}_{7^{(5)}}} ="38",
(70,0)  *+{\vdots} ="39",
(40,-10) *+{(a)} ="40",
"1", {\ar"0"},
"8", {\ar"0"},
"0", {\ar"9"},
"2", {\ar"1"},
"9", {\ar"1"},
"1", {\ar"10"},
"3", {\ar"2"},
"10", {\ar"2"},
"2", {\ar"11"},
"4", {\ar"3"},
"11", {\ar"3"},
"3", {\ar"12"},
"5", {\ar"4"},
"12", {\ar"4"},
"4", {\ar"13"},
"6", {\ar"5"},
"13", {\ar"5"},
"5", {\ar"14"},
"7", {\ar"6"},
"14", {\ar"6"},
"6", {\ar"15"},
"15", {\ar"7"},
"9", {\ar"8"},
"8", {\ar"16"},
"10", {\ar"9"},
"16", {\ar"9"},
"9", {\ar"17"},
"11", {\ar"10"},
"17", {\ar"10"},
"10", {\ar"18"},
"12", {\ar"11"},
"18", {\ar"11"},
"11", {\ar"19"},
"13", {\ar"12"},
"19", {\ar"12"},
"12", {\ar"20"},
"14", {\ar"13"},
"20", {\ar"13"},
"13", {\ar"21"},
"15", {\ar"14"},
"21", {\ar"14"},
"14", {\ar"22"},
"22", {\ar"15"},
"15", {\ar"23"},
"17", {\ar"16"},
"24", {\ar"16"},
"16", {\ar"25"},
"18", {\ar"17"},
"25", {\ar"17"},
"17", {\ar"26"},
"19", {\ar"18"},
"26", {\ar"18"},
"18", {\ar"27"},
"20", {\ar"19"},
"27", {\ar"19"},
"19", {\ar"28"},
"21", {\ar"20"},
"28", {\ar"20"},
"20", {\ar"29"},
"22", {\ar"21"},
"29", {\ar"21"},
"21", {\ar"30"},
"23", {\ar"22"},
"30", {\ar"22"},
"22", {\ar"31"},
"31", {\ar"23"},
"25", {\ar"24"},
"24", {\ar"32"},
"26", {\ar"25"},
"32", {\ar"25"},
"25", {\ar"33"},
"27", {\ar"26"},
"33", {\ar"26"},
"26", {\ar"34"},
"28", {\ar"27"},
"34", {\ar"27"},
"27", {\ar"35"},
"29", {\ar"28"},
"35", {\ar"28"},
"28", {\ar"36"},
"30", {\ar"29"},
"36", {\ar"29"},
"29", {\ar"37"},
"31", {\ar"30"},
"37", {\ar"30"},
"30", {\ar"38"},
"38", {\ar"31"},
"31", {\ar"39"},
"33", {\ar"32"},
"34", {\ar"33"},
"35", {\ar"34"},
"36", {\ar"35"},
"37", {\ar"36"},
"38", {\ar"37"},
"39", {\ar"38"},
\end{xy}}
\end{minipage}
\begin{minipage}[b]{.6\textwidth}
\centerline{
\begin{xy}
(10,70) *+{s^{(-1)}_{1^{(1)}}} ="0",
(10,60) *+{s^{(-3)}_{2^{(1)}}} ="1",
(10,50) *+{s^{(-5)}_{3^{(1)}}} ="2",
(10,40) *+{s^{(-7)}_{4^{(1)}}} ="3",
(10,30) *+{s^{(-9)}_{5^{(1)}}} ="4",
(10,20) *+{s^{(-11)}_{6^{(1)}}} ="5",
(10,10) *+{s^{(-13)}_{7^{(1)}}} ="6",
(10,0) *+{\vdots} ="7",
(25,75) *+{s^{(0)}_{1^{(2)}}} ="8",
(25,65) *+{s^{(-2)}_{2^{(2)}}} ="9",
(25,55) *+{s^{(-4)}_{3^{(2)}}} ="10",
(25,45) *+{s^{(-6)}_{4^{(2)}}} ="11",
(25,35) *+{s^{(-8)}_{5^{(2)}}} ="12",
(25,25) *+{s^{(-10)}_{6^{(2)}}} ="13",
(25,15) *+{s^{(-12)}_{7^{(2)}}} ="14",
(25,5)  *+{\vdots} ="15",
(40,70) *+{s^{(-1)}_{1^{(3)}}} ="16",
(40,60) *+{s^{(-3)}_{2^{(3)}}} ="17",
(40,50) *+{s^{(-5)}_{3^{(3)}}} ="18",
(40,40) *+{s^{(-7)}_{4^{(3)}}} ="19",
(40,30) *+{s^{(-9)}_{5^{(3)}}} ="20",
(40,20) *+{s^{(-11)}_{6^{(3)}}} ="21",
(40,10) *+{s^{(-13)}_{7^{(3)}}} ="22",
(40,0)  *+{\vdots} ="23",
(55,75) *+{s^{(0)}_{1^{(4)}}} ="24",
(55,65) *+{s^{(-2)}_{2^{(4)}}} ="25",
(55,55) *+{s^{(-4)}_{3^{(4)}}} ="26",
(55,45) *+{s^{(-6)}_{4^{(4)}}} ="27",
(55,35) *+{s^{(-8)}_{5^{(4)}}} ="28",
(55,25) *+{s^{(-10)}_{6^{(4)}}} ="29",
(55,15) *+{s^{(-12)}_{7^{(4)}}} ="30",
(55,5)  *+{\vdots} ="31",
(75,70) *+{s^{(-3)}_{1^{(5)},1^{(4)}}} ="32",
(75,60) *+[red]{s^{(-5)}_{2^{(5)},1^{(4)}}} ="33",
(75,50) *+{s^{(-7)}_{3^{(5)},1^{(4)}}} ="34",
(75,40) *+{s^{(-9)}_{4^{(5)},1^{(4)}}} ="35",
(75,30) *+{s^{(-11)}_{5^{(5)},1^{(4)}}} ="36",
(75,20) *+{s^{(-13)}_{6^{(5)},1^{(4)}}} ="37",
(75,10) *+{s^{(-15)}_{7^{(5)},1^{(4)}}} ="38",
(75,0)  *+{\vdots} ="39",
(40,-10) *+{(b)} ="40",
"1", {\ar"0"},
"8", {\ar"0"},
"0", {\ar"9"},
"2", {\ar"1"},
"9", {\ar"1"},
"1", {\ar"10"},
"3", {\ar"2"},
"10", {\ar"2"},
"2", {\ar"11"},
"4", {\ar"3"},
"11", {\ar"3"},
"3", {\ar"12"},
"5", {\ar"4"},
"12", {\ar"4"},
"4", {\ar"13"},
"6", {\ar"5"},
"13", {\ar"5"},
"5", {\ar"14"},
"7", {\ar"6"},
"14", {\ar"6"},
"6", {\ar"15"},
"15", {\ar"7"},
"9", {\ar"8"},
"8", {\ar"16"},
"10", {\ar"9"},
"16", {\ar"9"},
"9", {\ar"17"},
"11", {\ar"10"},
"17", {\ar"10"},
"10", {\ar"18"},
"12", {\ar"11"},
"18", {\ar"11"},
"11", {\ar"19"},
"13", {\ar"12"},
"19", {\ar"12"},
"12", {\ar"20"},
"14", {\ar"13"},
"20", {\ar"13"},
"13", {\ar"21"},
"15", {\ar"14"},
"21", {\ar"14"},
"14", {\ar"22"},
"22", {\ar"15"},
"15", {\ar"23"},
"17", {\ar"16"},
"24", {\ar"16"},
"16", {\ar"25"},
"18", {\ar"17"},
"25", {\ar"17"},
"17", {\ar"26"},
"19", {\ar"18"},
"26", {\ar"18"},
"18", {\ar"27"},
"20", {\ar"19"},
"27", {\ar"19"},
"19", {\ar"28"},
"21", {\ar"20"},
"28", {\ar"20"},
"20", {\ar"29"},
"22", {\ar"21"},
"29", {\ar"21"},
"21", {\ar"30"},
"23", {\ar"22"},
"30", {\ar"22"},
"22", {\ar"31"},
"31", {\ar"23"},
"32", {\ar"24"},
"26", {\ar"25"},
"25", {\ar"32"},
"27", {\ar"26"},
"32", {\ar"26"},
"26", {\ar"33"},
"28", {\ar"27"},
"33", {\ar"27"},
"27", {\ar"34"},
"29", {\ar"28"},
"34", {\ar"28"},
"28", {\ar"35"},
"30", {\ar"29"},
"35", {\ar"29"},
"29", {\ar"36"},
"31", {\ar"30"},
"36", {\ar"30"},
"30", {\ar"37"},
"37", {\ar"31"},
"31", {\ar"38"},
"33", {\ar"32"},
"34", {\ar"33"},
"35", {\ar"34"},
"36", {\ar"35"},
"37", {\ar"36"},
"38", {\ar"37"},
"39", {\ar"38"},
\end{xy}}
\end{minipage}}
\end{figure}
\begin{figure}[H]
\centering
\resizebox{.6\width}{.6\height}{
\begin{minipage}[b]{0.8\textwidth}
\centerline{
\begin{xy}
(8,70) *+{s^{(-1)}_{1^{(1)}}} ="0",
(8,60) *+{s^{(-5)}_{1^{(1)},2^{(2)}}} ="1",
(8,50) *+{s^{(-7)}_{2^{(1)},2^{(2)}}} ="2",
(8,40) *+{s^{(-9)}_{3^{(1)},2^{(2)}}} ="3",
(8,30) *+{s^{(-11)}_{4^{(1)},2^{(2)}}} ="4",
(8,20) *+{s^{(-13)}_{5^{(1)},2^{(2)}}} ="5",
(8,10) *+{s^{(-15)}_{6^{(1)},2^{(2)}}} ="6",
(8,0) *+{\vdots} ="7",
(25,75) *+{s^{(0)}_{1^{(2)}}} ="8",
(25,65) *+{s^{(-2)}_{2^{(2)}}} ="9",
(25,55) *+{s^{(-6)}_{1^{(2)},2^{(3)}}} ="10",
(25,45) *+{s^{(-8)}_{2^{(2)},2^{(3)}}} ="11",
(25,35) *+{s^{(-10)}_{3^{(2)},2^{(3)}}} ="12",
(25,25) *+{s^{(-12)}_{4^{(2)},2^{(3)}}}="13",
(25,15) *+{s^{(-14)}_{5^{(2)},2^{(3)}}}="14",
(25,5)  *+{\vdots} ="15",
(42,70) *+{s^{(-1)}_{1^{(3)}}} ="16",
(42,60) *+{s^{(-3)}_{2^{(3)}}} ="17",
(42,50) *+{s^{(-7)}_{1^{(3)},3^{(4)}}}="18",
(42,40) *+{s^{(-9)}_{2^{(3)},3^{(4)}}}="19",
(42,30) *+{s^{(-11)}_{3^{(3)},3^{(4)}}}="20",
(42,20) *+{s^{(-13)}_{4^{(3)},3^{(4)}}}="21",
(42,10) *+{s^{(-15)}_{5^{(3)},3^{(4)}}}="22",
(42,0)  *+{\vdots} ="23",
(62,75) *+{s^{(0)}_{1^{(4)}}} ="24",
(62,65) *+{s^{(-2)}_{2^{(4)}}} ="25",
(62,55) *+{s^{(-4)}_{3^{(4)}}} ="26",
(62,45) *+[red]{s^{(-8)}_{1^{(4)},2^{(5)},1^{(4)}}} ="27",
(62,35) *+{s^{(-10)}_{2^{(4)},2^{(5)},1^{(4)}}} ="28",
(62,25) *+{s^{(-12)}_{3^{(4)},2^{(5)},1^{(4)}}} ="29",
(62,15) *+{s^{(-14)}_{4^{(4)},2^{(5)},1^{(4)}}} ="30",
(62,5)  *+{\vdots} ="31",
(82,70) *+{s^{(-3)}_{1^{(5)},1^{(4)}}} ="32",
(82,60) *+{s^{(-5)}_{2^{(5)},1^{(4)}}} ="33",
(82,50) *+{s^{(-7)}_{3^{(5)},1^{(4)}}} ="34",
(82,40) *+{s^{(-9)}_{4^{(5)},1^{(4)}}} ="35",
(82,30) *+{s^{(-11)}_{5^{(5)},1^{(4)}}} ="36",
(82,20) *+{s^{(-13)}_{6^{(5)},1^{(4)}}} ="37",
(82,10) *+{s^{(-15)}_{7^{(5)},1^{(4)}}} ="38",
(82,0)  *+{\vdots} ="39",
(40,-10) *+{(c)} ="40",
"0", {\ar@/_2pc/"7"},
"8", {\ar"0"},
"2", {\ar"1"},
"1", {\ar"9"},
"1", {\ar"10"},
"17", {\ar"1"},
"3", {\ar"2"},
"10", {\ar"2"},
"2", {\ar"11"},
"4", {\ar"3"},
"11", {\ar"3"},
"3", {\ar"12"},
"5", {\ar"4"},
"12", {\ar"4"},
"4", {\ar"13"},
"6", {\ar"5"},
"13", {\ar"5"},
"5", {\ar"14"},
"7", {\ar"6"},
"14", {\ar"6"},
"6", {\ar"15"},
"9", {\ar"8"},
"8", {\ar"16"},
"16", {\ar"9"},
"9", {\ar"17"},
"11", {\ar"10"},
"10", {\ar"17"},
"10", {\ar"18"},
"26", {\ar"10"},
"12", {\ar"11"},
"18", {\ar"11"},
"11", {\ar"19"},
"13", {\ar"12"},
"19", {\ar"12"},
"12", {\ar"20"},
"14", {\ar"13"},
"20", {\ar"13"},
"13", {\ar"21"},
"15", {\ar"14"},
"21", {\ar"14"},
"14", {\ar"22"},
"22", {\ar"15"},
"15", {\ar"23"},
"17", {\ar"16"},
"24", {\ar"16"},
"16", {\ar"25"},
"25", {\ar"17"},
"17", {\ar"26"},
"19", {\ar"18"},
"18", {\ar"26"},
"18", {\ar"27"},
"33", {\ar@/^1pc/"18"},
"20", {\ar"19"},
"27", {\ar"19"},
"19", {\ar"28"},
"21", {\ar"20"},
"28", {\ar"20"},
"20", {\ar"29"},
"22", {\ar"21"},
"29", {\ar"21"},
"21", {\ar"30"},
"23", {\ar"22"},
"30", {\ar"22"},
"22", {\ar"31"},
"32", {\ar"24"},
"26", {\ar"25"},
"25", {\ar"32"},
"32", {\ar"26"},
"26", {\ar"33"},
"28", {\ar"27"},
"27", {\ar"33"},
"34", {\ar"27"},
"27", {\ar"35"},
"29", {\ar"28"},
"35", {\ar"28"},
"28", {\ar"36"},
"30", {\ar"29"},
"36", {\ar"29"},
"29", {\ar"37"},
"31", {\ar"30"},
"37", {\ar"30"},
"30", {\ar"38"},
"38", {\ar"31"},
"33", {\ar"32"},
"35", {\ar"34"},
"36", {\ar"35"},
"37", {\ar"36"},
"38", {\ar"37"},
"39", {\ar"38"},
\end{xy}}
\end{minipage}
\begin{minipage}[b]{0.7\textwidth}
\centerline{
\begin{xy}
(5,70) *+{s^{(-1)}_{1^{(1)}}} ="0",
(5,60) *+{s^{(-7)}_{1^{(1)},2^{(3)}}} ="1",
(5,50) *+{s^{(-9)}_{2^{(1)},2^{(3)}}} ="2",
(5,40) *+{s^{(-11)}_{3^{(1)},2^{(3)}}} ="3",
(5,30) *+{s^{(-13)}_{4^{(1)},2^{(3)}}} ="4",
(5,20) *+{s^{(-15)}_{5^{(1)},2^{(3)}}}="5",
(5,10) *+{s^{(-17)}_{6^{(1)},2^{(3)}}}="6",
(5,0) *+{\vdots} ="7",
(25,75) *+{s^{(0)}_{1^{(2)}}} ="8",
(25,65) *+{s^{(-2)}_{2^{(2)}}} ="9",
(25,55) *+{s^{(-8)}_{1^{(2)},2^{(4)}}} ="10",
(25,45) *+{s^{(-10)}_{2^{(2)},2^{(4)}}} ="11",
(25,35) *+{s^{(-12)}_{3^{(2)},2^{(4)}}} ="12",
(25,25) *+{s^{(-14)}_{4^{(2)},2^{(4)}}}="13",
(25,15) *+{s^{(-16)}_{5^{(2)},2^{(4)}}}="14",
(25,5)  *+{\vdots} ="15",
(42,70) *+{s^{(-1)}_{1^{(3)}}} ="16",
(42,60) *+{s^{(-3)}_{2^{(3)}}} ="17",
(42,50) *+[red]{s^{(-9)}_{1^{(3)},2^{(5)},1^{(4)}}}="18",
(42,40) *+{s^{(-11)}_{2^{(3)},2^{(5)},1^{(4)}}}="19",
(42,30) *+{s^{(-13)}_{3^{(3)},2^{(5)},1^{(4)}}}="20",
(42,20) *+{s^{(-15)}_{4^{(3)},2^{(5)},1^{(4)}}}="21",
(42,10) *+{s^{(-17)}_{5^{(3)},2^{(5)},1^{(4)}}}="22",
(42,0)  *+{\vdots} ="23",
(62,75) *+{s^{(0)}_{1^{(4)}}} ="24",
(62,65) *+{s^{(-2)}_{2^{(4)}}} ="25",
(62,55) *+{s^{(-4)}_{3^{(4)}}} ="26",
(62,45) *+[red]{s^{(-8)}_{1^{(4)},2^{(5)},1^{(4)}}} ="27",
(62,35) *+{s^{(-10)}_{2^{(4)},2^{(5)},1^{(4)}}} ="28",
(62,25) *+{s^{(-12)}_{3^{(4)},2^{(5)},1^{(4)}}} ="29",
(62,15) *+{s^{(-14)}_{4^{(4)},2^{(5)},1^{(4)}}} ="30",
(62,5)  *+{\vdots} ="31",
(82,70) *+{s^{(-3)}_{1^{(5)},1^{(4)}}} ="32",
(82,60) *+{s^{(-5)}_{2^{(5)},1^{(4)}}} ="33",
(82,50) *+{s^{(-7)}_{3^{(5)},1^{(4)}}} ="34",
(82,40) *+{s^{(-9)}_{4^{(5)},1^{(4)}}} ="35",
(82,30) *+{s^{(-11)}_{5^{(5)},1^{(4)}}} ="36",
(82,20) *+{s^{(-13)}_{6^{(5)},1^{(4)}}} ="37",
(82,10) *+{s^{(-15)}_{7^{(5)},1^{(4)}}} ="38",
(82,0)  *+{\vdots} ="39",
(40,-10) *+{(d)} ="40",
"8", {\ar"0"},
"0", {\ar"9"},
"23", {\ar"0"},
"0", {\ar"31"},
"2", {\ar"1"},
"1", {\ar"10"},
"1", {\ar"17"},
"26", {\ar"1"},
"3", {\ar"2"},
"10", {\ar"2"},
"2", {\ar"11"},
"4", {\ar"3"},
"11", {\ar"3"},
"3", {\ar"12"},
"5", {\ar"4"},
"12", {\ar"4"},
"4", {\ar"13"},
"6", {\ar"5"},
"13", {\ar"5"},
"5", {\ar"14"},
"7", {\ar"6"},
"9", {\ar"7"},
"14", {\ar"7"},
"7", {\ar"15"},
"9", {\ar"8"},
"8", {\ar"16"},
"16", {\ar"9"},
"11", {\ar"10"},
"10", {\ar"18"},
"10", {\ar"26"},
"33", {\ar@/^1pc/"10"},
"12", {\ar"11"},
"18", {\ar"11"},
"11", {\ar"19"},
"13", {\ar"12"},
"19", {\ar"12"},
"12", {\ar"20"},
"14", {\ar"13"},
"20", {\ar"13"},
"13", {\ar"21"},
"15", {\ar"14"},
"21", {\ar"14"},
"14", {\ar"22"},
"22", {\ar"15"},
"15", {\ar"23"},
"17", {\ar"16"},
"24", {\ar"16"},
"16", {\ar"25"},
"25", {\ar"17"},
"17", {\ar"26"},
"19", {\ar"18"},
"27", {\ar"18"},
"18", {\ar"28"},
"18", {\ar@/_1pc/"33"},
"20", {\ar"19"},
"28", {\ar"19"},
"19", {\ar"29"},
"21", {\ar"20"},
"29", {\ar"20"},
"20", {\ar"30"},
"22", {\ar"21"},
"30", {\ar"21"},
"21", {\ar"31"},
"23", {\ar"22"},
"31", {\ar"23"},
"32", {\ar"24"},
"26", {\ar"25"},
"25", {\ar"32"},
"32", {\ar"26"},
"26", {\ar"33"},
"28", {\ar"27"},
"34", {\ar"27"},
"27", {\ar"35"},
"29", {\ar"28"},
"35", {\ar"28"},
"28", {\ar"36"},
"30", {\ar"29"},
"36", {\ar"29"},
"29", {\ar"37"},
"31", {\ar"30"},
"37", {\ar"30"},
"30", {\ar"38"},
"38", {\ar"31"},
"33", {\ar"32"},
"35", {\ar"34"},
"36", {\ar"35"},
"37", {\ar"36"},
"38", {\ar"37"},
"39", {\ar"38"},
\end{xy}}
\end{minipage}}
\end{figure}
\begin{figure}[H]
\centering
\resizebox{.55\width}{.55\height}{
\begin{minipage}[b]{0.8\textwidth}
\centerline{
\begin{xy}
(5,70) *+{s^{(-1)}_{1^{(1)}}} ="0",
(5,60) *+{s^{(-9)}_{1^{(1)},3^{(4)}}} ="1",
(5,50) *+{s^{(-11)}_{2^{(1)},3^{(4)}}} ="2",
(5,40) *+{s^{(-13)}_{3^{(1)},3^{(4)}}} ="3",
(5,30) *+{s^{(-15)}_{4^{(1)},3^{(4)}}} ="4",
(5,20) *+{s^{(-17)}_{5^{(1)},3^{(4)}}}="5",
(5,10) *+{s^{(-19)}_{6^{(1)},3^{(4)}}}="6",
(5,0) *+{\vdots} ="7",
(25,75) *+{s^{(0)}_{1^{(2)}}} ="8",
(25,65) *+{s^{(-2)}_{2^{(2)}}} ="9",
(25,55) *+[red]{s^{(-10)}_{1^{(2)},2^{(5)},1^{(4)}}} ="10",
(25,45) *+{s^{(-12)}_{2^{(2)},2^{(5)},1^{(4)}}} ="11",
(25,35) *+{s^{(-14)}_{3^{(2)},2^{(5)},1^{(4)}}} ="12",
(25,25) *+{s^{(-16)}_{4^{(2)},2^{(5)},1^{(4)}}}="13",
(25,15) *+{s^{(-18)}_{5^{(2)},2^{(5)},1^{(4)}}}="14",
(25,5)  *+{\vdots} ="15",
(42,70) *+{s^{(-1)}_{1^{(3)}}} ="16",
(42,60) *+{s^{(-3)}_{2^{(3)}}} ="17",
(42,50) *+[red]{s^{(-9)}_{1^{(3)},2^{(5)},1^{(4)}}}="18",
(42,40) *+{s^{(-11)}_{2^{(3)},2^{(5)},1^{(4)}}}="19",
(42,30) *+{s^{(-13)}_{3^{(3)},2^{(5)},1^{(4)}}}="20",
(42,20) *+{s^{(-15)}_{4^{(3)},2^{(5)},1^{(4)}}}="21",
(42,10) *+{s^{(-17)}_{5^{(3)},2^{(5)},1^{(4)}}}="22",
(42,0)  *+{\vdots} ="23",
(62,75) *+{s^{(0)}_{1^{(4)}}} ="24",
(62,65) *+{s^{(-2)}_{2^{(4)}}} ="25",
(62,55) *+{s^{(-4)}_{3^{(4)}}} ="26",
(62,45) *+[red]{s^{(-8)}_{1^{(4)},2^{(5)},1^{(4)}}} ="27",
(62,35) *+{s^{(-10)}_{2^{(4)},2^{(5)},1^{(4)}}} ="28",
(62,25) *+{s^{(-12)}_{3^{(4)},2^{(5)},1^{(4)}}} ="29",
(62,15) *+{s^{(-14)}_{4^{(4)},2^{(5)},1^{(4)}}} ="30",
(62,5)  *+{\vdots} ="31",
(82,70) *+{s^{(-3)}_{1^{(5)},1^{(4)}}} ="32",
(82,60) *+{s^{(-5)}_{2^{(5)},1^{(4)}}} ="33",
(82,50) *+{s^{(-7)}_{3^{(5)},1^{(4)}}} ="34",
(82,40) *+{s^{(-9)}_{4^{(5)},1^{(4)}}} ="35",
(82,30) *+{s^{(-11)}_{5^{(5)},1^{(4)}}} ="36",
(82,20) *+{s^{(-13)}_{6^{(5)},1^{(4)}}} ="37",
(82,10) *+{s^{(-15)}_{7^{(5)},1^{(4)}}} ="38",
(82,0)  *+{\vdots} ="39",
(40,-10) *+{(e)} ="40",
"8", {\ar"0"},
"0", {\ar"9"},
"23", {\ar"0"},
"0", {\ar"31"},
"2", {\ar"1"},
"1", {\ar"10"},
"1", {\ar"26"},
"33", {\ar@/_2pc/"1"},
"3", {\ar"2"},
"10", {\ar"2"},
"2", {\ar"11"},
"4", {\ar"3"},
"11", {\ar"3"},
"3", {\ar"12"},
"5", {\ar"4"},
"12", {\ar"4"},
"4", {\ar"13"},
"6", {\ar"5"},
"7", {\ar"6"},
"13", {\ar"6"},
"6", {\ar"14"},
"14", {\ar"7"},
"7", {\ar"15"},
"17", {\ar"7"},
"9", {\ar"8"},
"8", {\ar"16"},
"15", {\ar@/^3pc/"9"},
"16", {\ar"9"},
"9", {\ar"17"},
"9", {\ar"23"},
"11", {\ar"10"},
"18", {\ar"10"},
"10", {\ar"19"},
"10", {\ar"33"},
"12", {\ar"11"},
"19", {\ar"11"},
"11", {\ar"20"},
"13", {\ar"12"},
"20", {\ar"12"},
"12", {\ar"21"},
"14", {\ar"13"},
"21", {\ar"13"},
"13", {\ar"22"},
"15", {\ar"14"},
"22", {\ar"14"},
"14", {\ar"23"},
"23", {\ar"15"},
"17", {\ar"16"},
"24", {\ar"16"},
"16", {\ar"25"},
"25", {\ar"17"},
"19", {\ar"18"},
"27", {\ar"18"},
"18", {\ar"28"},
"20", {\ar"19"},
"28", {\ar"19"},
"19", {\ar"29"},
"21", {\ar"20"},
"29", {\ar"20"},
"20", {\ar"30"},
"22", {\ar"21"},
"30", {\ar"21"},
"21", {\ar"31"},
"23", {\ar"22"},
"31", {\ar"23"},
"32", {\ar"24"},
"26", {\ar"25"},
"25", {\ar"32"},
"32", {\ar"26"},
"26", {\ar"33"},
"28", {\ar"27"},
"34", {\ar"27"},
"27", {\ar"35"},
"29", {\ar"28"},
"35", {\ar"28"},
"28", {\ar"36"},
"30", {\ar"29"},
"36", {\ar"29"},
"29", {\ar"37"},
"31", {\ar"30"},
"37", {\ar"30"},
"30", {\ar"38"},
"38", {\ar"31"},
"33", {\ar"32"},
"35", {\ar"34"},
"36", {\ar"35"},
"37", {\ar"36"},
"38", {\ar"37"},
"39", {\ar"38"},
\end{xy}}
\end{minipage}
\begin{minipage}[b]{0.7\textwidth}
\centerline{
\begin{xy}
(0,70) *+{s^{(-1)}_{1^{(1)}}} ="0",
(0,60) *+[red]{s^{(-11)}_{1^{(1)},2^{(5)},1^{(4)}}} ="1",
(0,50) *+{s^{(-13)}_{2^{(2)},2^{(5)},1^{(4)}}} ="2",
(0,40) *+{s^{(-15)}_{3^{(2)},2^{(5)},1^{(4)}}} ="3",
(0,30) *+{s^{(-17)}_{4^{(2)},2^{(5)},1^{(4)}}} ="4",
(0,20) *+{s^{(-19)}_{5^{(2)},2^{(5)},1^{(4)}}}="5",
(0,10) *+{s^{(-21)}_{6^{(2)},2^{(5)},1^{(4)}}}="6",
(0,0) *+{\vdots} ="7",
(25,75) *+{s^{(0)}_{1^{(2)}}} ="8",
(25,65) *+{s^{(-2)}_{2^{(2)}}} ="9",
(25,55) *+[red]{s^{(-10)}_{1^{(2)},2^{(5)},1^{(4)}}} ="10",
(25,45) *+{s^{(-12)}_{2^{(2)},2^{(5)},1^{(4)}}} ="11",
(25,35) *+{s^{(-14)}_{3^{(2)},2^{(5)},1^{(4)}}} ="12",
(25,25) *+{s^{(-16)}_{4^{(2)},2^{(5)},1^{(4)}}}="13",
(25,15) *+{s^{(-18)}_{5^{(2)},2^{(5)},1^{(4)}}}="14",
(25,5)  *+{\vdots} ="15",
(42,70) *+{s^{(-1)}_{1^{(3)}}} ="16",
(42,60) *+{s^{(-3)}_{2^{(3)}}} ="17",
(42,50) *+[red]{s^{(-9)}_{1^{(3)},2^{(5)},1^{(4)}}}="18",
(42,40) *+{s^{(-11)}_{2^{(3)},2^{(5)},1^{(4)}}}="19",
(42,30) *+{s^{(-13)}_{3^{(3)},2^{(5)},1^{(4)}}}="20",
(42,20) *+{s^{(-15)}_{4^{(3)},2^{(5)},1^{(4)}}}="21",
(42,10) *+{s^{(-17)}_{5^{(3)},2^{(5)},1^{(4)}}}="22",
(42,0)  *+{\vdots} ="23",
(62,75) *+{s^{(0)}_{1^{(4)}}} ="24",
(62,65) *+{s^{(-2)}_{2^{(4)}}} ="25",
(62,55) *+{s^{(-4)}_{3^{(4)}}} ="26",
(62,45) *+[red]{s^{(-8)}_{1^{(4)},2^{(5)},1^{(4)}}} ="27",
(62,35) *+{s^{(-10)}_{2^{(4)},2^{(5)},1^{(4)}}} ="28",
(62,25) *+{s^{(-12)}_{3^{(4)},2^{(5)},1^{(4)}}} ="29",
(62,15) *+{s^{(-14)}_{4^{(4)},2^{(5)},1^{(4)}}} ="30",
(62,5)  *+{\vdots} ="31",
(82,70) *+{s^{(-3)}_{1^{(5)},1^{(4)}}} ="32",
(82,60) *+{s^{(-5)}_{2^{(5)},1^{(4)}}} ="33",
(82,50) *+{s^{(-7)}_{3^{(5)},1^{(4)}}} ="34",
(82,40) *+{s^{(-9)}_{4^{(5)},1^{(4)}}} ="35",
(82,30) *+{s^{(-11)}_{5^{(5)},1^{(4)}}} ="36",
(82,20) *+{s^{(-13)}_{6^{(5)},1^{(4)}}} ="37",
(82,10) *+{s^{(-15)}_{7^{(5)},1^{(4)}}} ="38",
(82,0)  *+{\vdots} ="39",
(40,-10) *+{(f)} ="40",
"8", {\ar"0"},
"0", {\ar"9"},
"23", {\ar"0"},
"0", {\ar"31"},
"2", {\ar"1"},
"10", {\ar"1"},
"1", {\ar"11"},
"1", {\ar@/^1pc/"33"},
"3", {\ar"2"},
"11", {\ar"2"},
"2", {\ar"12"},
"4", {\ar"3"},
"12", {\ar"3"},
"3", {\ar"13"},
"5", {\ar"4"},
"6", {\ar"5"},
"13", {\ar"5"},
"5", {\ar"14"},
"7", {\ar"6"},
"14", {\ar"6"},
"6", {\ar"15"},
"15", {\ar"7"},
"7", {\ar"17"},
"26", {\ar"7"},
"9", {\ar"8"},
"8", {\ar"16"},
"15", {\ar@/^3pc/"9"},
"16", {\ar"9"},
"9", {\ar"17"},
"9", {\ar"23"},
"11", {\ar"10"},
"18", {\ar"10"},
"10", {\ar"19"},
"12", {\ar"11"},
"19", {\ar"11"},
"11", {\ar"20"},
"13", {\ar"12"},
"20", {\ar"12"},
"12", {\ar"21"},
"14", {\ar"13"},
"21", {\ar"13"},
"13", {\ar"22"},
"15", {\ar"14"},
"22", {\ar"14"},
"14", {\ar"23"},
"17", {\ar"15"},
"23", {\ar"15"},
"17", {\ar"16"},
"24", {\ar"16"},
"16", {\ar"25"},
"25", {\ar"17"},
"17", {\ar"26"},
"19", {\ar"18"},
"27", {\ar"18"},
"18", {\ar"28"},
"20", {\ar"19"},
"28", {\ar"19"},
"19", {\ar"29"},
"21", {\ar"20"},
"29", {\ar"20"},
"20", {\ar"30"},
"22", {\ar"21"},
"30", {\ar"21"},
"21", {\ar"31"},
"23", {\ar"22"},
"31", {\ar"23"},
"32", {\ar"24"},
"26", {\ar"25"},
"25", {\ar"32"},
"32", {\ar"26"},
"28", {\ar"27"},
"34", {\ar"27"},
"27", {\ar"35"},
"29", {\ar"28"},
"35", {\ar"28"},
"28", {\ar"36"},
"30", {\ar"29"},
"36", {\ar"29"},
"29", {\ar"37"},
"31", {\ar"30"},
"37", {\ar"30"},
"30", {\ar"38"},
"38", {\ar"31"},
"33", {\ar"32"},
"35", {\ar"34"},
"36", {\ar"35"},
"37", {\ar"36"},
"38", {\ar"37"},
"39", {\ar"38"},
\end{xy}}
\end{minipage}}
\end{figure}
\begin{figure}[H]
\centering
\resizebox{.55\width}{.55\height}{
\begin{minipage}[b]{0.8\textwidth}
\centerline{
\begin{xy}
(0,70) *+{s^{(-1)}_{1^{(1)}}} ="0",
(0,60) *+{s^{(-11)}_{1^{(1)},2^{(5)},1^{(4)}}} ="1",
(0,50) *+{s^{(-13)}_{2^{(2)},2^{(5)},1^{(4)}}} ="2",
(0,40) *+{s^{(-15)}_{3^{(2)},2^{(5)},1^{(4)}}} ="3",
(0,30) *+{s^{(-17)}_{4^{(2)},2^{(5)},1^{(4)}}} ="4",
(0,20) *+{s^{(-19)}_{5^{(2)},2^{(5)},1^{(4)}}}="5",
(0,10) *+{s^{(-21)}_{6^{(2)},2^{(5)},1^{(4)}}}="6",
(0,0) *+{\vdots} ="7",
(25,75) *+{s^{(0)}_{1^{(2)}}} ="8",
(25,65) *+{s^{(-2)}_{2^{(2)}}} ="9",
(25,55) *+{s^{(-10)}_{1^{(2)},2^{(5)},1^{(4)}}} ="10",
(25,45) *+{s^{(-12)}_{2^{(2)},2^{(5)},1^{(4)}}} ="11",
(25,35) *+{s^{(-14)}_{3^{(2)},2^{(5)},1^{(4)}}} ="12",
(25,25) *+{s^{(-16)}_{4^{(2)},2^{(5)},1^{(4)}}}="13",
(25,15) *+{s^{(-18)}_{5^{(2)},2^{(5)},1^{(4)}}}="14",
(25,5)  *+{\vdots} ="15",
(50,70) *+{s^{(-1)}_{1^{(3)}}} ="16",
(50,60) *+{s^{(-3)}_{2^{(3)}}} ="17",
(50,50) *+[red]{s^{(-11)}_{1^{(3)},1^{(4)},2^{(5)},1^{(4)}}}="18",
(50,40) *+{s^{(-13)}_{2^{(3)},1^{(4)},2^{(5)},1^{(4)}}}="19",
(50,30) *+{s^{(-15)}_{3^{(3)},1^{(4)},2^{(5)},1^{(4)}}}="20",
(50,20) *+{s^{(-17)}_{4^{(3)},1^{(4)},2^{(5)},1^{(4)}}}="21",
(50,10) *+{s^{(-19)}_{5^{(3)},1^{(4)},2^{(5)},1^{(4)}}}="22",
(50,0)  *+{\vdots} ="23",
(75,75) *+{s^{(0)}_{1^{(4)}}} ="24",
(75,65) *+{s^{(-2)}_{2^{(4)}}} ="25",
(75,55) *+{s^{(-4)}_{3^{(4)}}} ="26",
(75,45) *+[red]{s^{(-8)}_{1^{(4)},2^{(5)},1^{(4)}}} ="27",
(75,35) *+{s^{(-10)}_{2^{(4)},2^{(5)},1^{(4)}}} ="28",
(75,25) *+{s^{(-12)}_{3^{(4)},2^{(5)},1^{(4)}}} ="29",
(75,15) *+{s^{(-14)}_{4^{(4)},2^{(5)},1^{(4)}}} ="30",
(75,5)  *+{\vdots} ="31",
(95,70) *+{s^{(-3)}_{1^{(5)},1^{(4)}}} ="32",
(95,60) *+{s^{(-5)}_{2^{(5)},1^{(4)}}} ="33",
(95,50) *+{s^{(-7)}_{3^{(5)},1^{(4)}}} ="34",
(95,40) *+{s^{(-9)}_{4^{(5)},1^{(4)}}} ="35",
(95,30) *+{s^{(-11)}_{5^{(5)},1^{(4)}}} ="36",
(95,20) *+{s^{(-13)}_{6^{(5)},1^{(4)}}} ="37",
(95,10) *+{s^{(-15)}_{7^{(5)},1^{(4)}}} ="38",
(95,0)  *+{\vdots} ="39",
(40,-10) *+{(g)} ="40",
"8", {\ar"0"},
"22", {\ar"0"},
"0", {\ar"23"},
"0", {\ar"31"},
"2", {\ar"1"},
"10", {\ar"1"},
"1", {\ar"11"},
"1", {\ar@/^1pc/"33"},
"3", {\ar"2"},
"11", {\ar"2"},
"2", {\ar"12"},
"4", {\ar"3"},
"12", {\ar"3"},
"3", {\ar"13"},
"5", {\ar"4"},
"6", {\ar"5"},
"13", {\ar"5"},
"5", {\ar"14"},
"7", {\ar"6"},
"14", {\ar"6"},
"6", {\ar"15"},
"15", {\ar"7"},
"7", {\ar"17"},
"26", {\ar"7"},
"9", {\ar"8"},
"8", {\ar"16"},
"16", {\ar"9"},
"9", {\ar"17"},
"23", {\ar"9"},
"11", {\ar"10"},
"10", {\ar"18"},
"27", {\ar@/_1pc/"10"},
"12", {\ar"11"},
"18", {\ar"11"},
"11", {\ar"19"},
"13", {\ar"12"},
"19", {\ar"12"},
"12", {\ar"20"},
"14", {\ar"13"},
"20", {\ar"13"},
"13", {\ar"21"},
"15", {\ar"14"},
"21", {\ar"14"},
"14", {\ar"22"},
"17", {\ar"15"},
"22", {\ar"15"},
"15", {\ar"23"},
"17", {\ar"16"},
"24", {\ar"16"},
"16", {\ar"25"},
"25", {\ar"17"},
"17", {\ar"26"},
"19", {\ar"18"},
"18", {\ar"27"},
"28", {\ar"18"},
"18", {\ar"29"},
"20", {\ar"19"},
"29", {\ar"19"},
"19", {\ar"30"},
"21", {\ar"20"},
"30", {\ar"20"},
"20", {\ar"31"},
"22", {\ar"21"},
"23", {\ar"22"},
"31", {\ar"22"},
"32", {\ar"24"},
"26", {\ar"25"},
"25", {\ar"32"},
"32", {\ar"26"},
"34", {\ar"27"},
"27", {\ar"35"},
"29", {\ar"28"},
"35", {\ar"28"},
"28", {\ar"36"},
"30", {\ar"29"},
"36", {\ar"29"},
"29", {\ar"37"},
"31", {\ar"30"},
"37", {\ar"30"},
"30", {\ar"38"},
"38", {\ar"31"},
"33", {\ar"32"},
"35", {\ar"34"},
"36", {\ar"35"},
"37", {\ar"36"},
"38", {\ar"37"},
"39", {\ar"38"},
\end{xy}}
\end{minipage}
\begin{minipage}[b]{0.8\linewidth}
\centerline{
\begin{xy}
(0,70) *+{s^{(-1)}_{1^{(1)}}} ="0",
(0,60) *+{s^{(-11)}_{1^{(1)},2^{(5)},1^{(4)}}} ="1",
(0,50) *+{s^{(-13)}_{2^{(2)},2^{(5)},1^{(4)}}} ="2",
(0,40) *+{s^{(-15)}_{3^{(2)},2^{(5)},1^{(4)}}} ="3",
(0,30) *+{s^{(-17)}_{4^{(2)},2^{(5)},1^{(4)}}} ="4",
(0,20) *+{s^{(-19)}_{5^{(2)},2^{(5)},1^{(4)}}}="5",
(0,10) *+{s^{(-21)}_{6^{(2)},2^{(5)},1^{(4)}}}="6",
(0,0) *+{\vdots} ="7",
(25,75) *+{s^{(0)}_{1^{(2)}}} ="8",
(25,65) *+{s^{(-2)}_{2^{(2)}}} ="9",
(25,55) *+[red]{s^{(-12)}_{1^{(2)},1^{(4)},2^{(5)},1^{(4)}}} ="10",
(25,45) *+{s^{(-14)}_{2^{(2)},1^{(4)},2^{(5)},1^{(4)}}} ="11",
(25,35) *+{s^{(-16)}_{3^{(2)},1^{(4)},2^{(5)},1^{(4)}}} ="12",
(25,25) *+{s^{(-18)}_{4^{(2)},1^{(4)},2^{(5)},1^{(4)}}}="13",
(25,15) *+{s^{(-20)}_{5^{(2)},1^{(4)},2^{(5)},1^{(4)}}}="14",
(25,5)  *+{\vdots} ="15",
(50,70) *+{s^{(-1)}_{1^{(3)}}} ="16",
(50,60) *+{s^{(-3)}_{2^{(3)}}} ="17",
(50,50) *+[red]{s^{(-11)}_{1^{(3)},1^{(4)},2^{(5)},1^{(4)}}}="18",
(50,40) *+{s^{(-13)}_{2^{(3)},1^{(4)},2^{(5)},1^{(4)}}}="19",
(50,30) *+{s^{(-15)}_{3^{(3)},1^{(4)},2^{(5)},1^{(4)}}}="20",
(50,20) *+{s^{(-17)}_{4^{(3)},1^{(4)},2^{(5)},1^{(4)}}}="21",
(50,10) *+{s^{(-19)}_{5^{(3)},1^{(4)},2^{(5)},1^{(4)}}}="22",
(50,0)  *+{\vdots} ="23",
(75,75) *+{s^{(0)}_{1^{(4)}}} ="24",
(75,65) *+{s^{(-2)}_{2^{(4)}}} ="25",
(75,55) *+{s^{(-4)}_{3^{(4)}}} ="26",
(75,45) *+[red]{s^{(-8)}_{1^{(4)},2^{(5)},1^{(4)}}} ="27",
(75,35) *+{s^{(-10)}_{2^{(4)},2^{(5)},1^{(4)}}} ="28",
(75,25) *+{s^{(-12)}_{3^{(4)},2^{(5)},1^{(4)}}} ="29",
(75,15) *+{s^{(-14)}_{4^{(4)},2^{(5)},1^{(4)}}} ="30",
(75,5)  *+{\vdots} ="31",
(95,70) *+{s^{(-3)}_{1^{(5)},1^{(4)}}} ="32",
(95,60) *+{s^{(-5)}_{2^{(5)},1^{(4)}}} ="33",
(95,50) *+{s^{(-7)}_{3^{(5)},1^{(4)}}} ="34",
(95,40) *+{s^{(-9)}_{4^{(5)},1^{(4)}}} ="35",
(95,30) *+{s^{(-11)}_{5^{(5)},1^{(4)}}} ="36",
(95,20) *+{s^{(-13)}_{6^{(5)},1^{(4)}}} ="37",
(95,10) *+{s^{(-15)}_{7^{(5)},1^{(4)}}} ="38",
(95,0)  *+{\vdots} ="39",
(40,-10) *+{(h)} ="40",
"8", {\ar"0"},
"22", {\ar"0"},
"0", {\ar"23"},
"0", {\ar"31"},
"2", {\ar"1"},
"1", {\ar"10"},
"27", {\ar@/_1pc/"1"},
"1", {\ar@/^1pc/"33"},
"3", {\ar"2"},
"10", {\ar"2"},
"2", {\ar"11"},
"4", {\ar"3"},
"11", {\ar"3"},
"3", {\ar"12"},
"5", {\ar"4"},
"6", {\ar"5"},
"12", {\ar"5"},
"5", {\ar"13"},
"7", {\ar"6"},
"13", {\ar"6"},
"6", {\ar"14"},
"14", {\ar"7"},
"7", {\ar"15"},
"26", {\ar"7"},
"9", {\ar"8"},
"8", {\ar"16"},
"16", {\ar"9"},
"9", {\ar"17"},
"23", {\ar"9"},
"11", {\ar"10"},
"18", {\ar"10"},
"10", {\ar"19"},
"10", {\ar@/^2pc/"27"},
"12", {\ar"11"},
"19", {\ar"11"},
"11", {\ar"20"},
"13", {\ar"12"},
"20", {\ar"12"},
"12", {\ar"21"},
"14", {\ar"13"},
"21", {\ar"13"},
"13", {\ar"22"},
"15", {\ar"14"},
"22", {\ar"14"},
"14", {\ar"23"},
"15", {\ar"17"},
"23", {\ar"15"},
"17", {\ar"16"},
"24", {\ar"16"},
"16", {\ar"25"},
"17", {\ar@/^2pc/"23"},
"25", {\ar"17"},
"17", {\ar"26"},
"19", {\ar"18"},
"28", {\ar"18"},
"18", {\ar"29"},
"20", {\ar"19"},
"29", {\ar"19"},
"19", {\ar"30"},
"21", {\ar"20"},
"30", {\ar"20"},
"20", {\ar"31"},
"22", {\ar"21"},
"23", {\ar"22"},
"31", {\ar"22"},
"32", {\ar"24"},
"26", {\ar"25"},
"25", {\ar"32"},
"32", {\ar"26"},
"34", {\ar"27"},
"27", {\ar"35"},
"29", {\ar"28"},
"35", {\ar"28"},
"28", {\ar"36"},
"30", {\ar"29"},
"36", {\ar"29"},
"29", {\ar"37"},
"31", {\ar"30"},
"37", {\ar"30"},
"30", {\ar"38"},
"38", {\ar"31"},
"33", {\ar"32"},
"35", {\ar"34"},
"36", {\ar"35"},
"37", {\ar"36"},
"38", {\ar"37"},
"39", {\ar"38"},
\end{xy}}
\end{minipage}}
\end{figure}
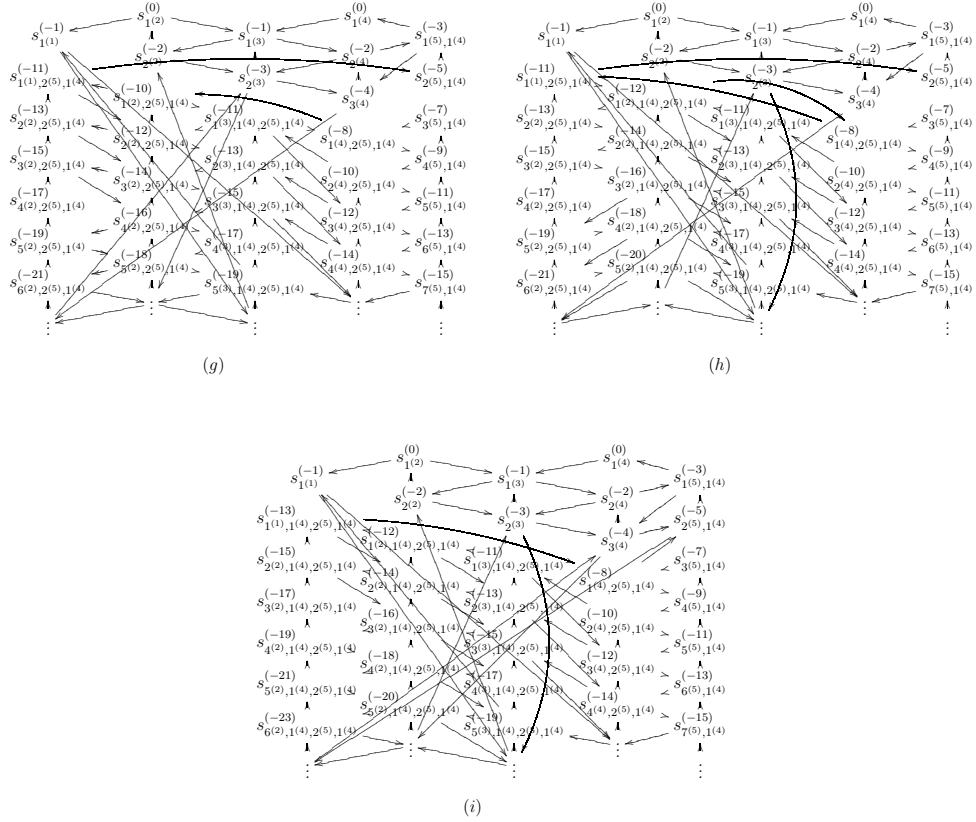
\begin{figure}[H]
\centering
\resizebox{.55\width}{.55\height}{
\begin{minipage}[b]{0.8\textwidth}
\centerline{
\begin{xy}
(0,70) *+{s^{(-1)}_{1^{(1)}}} ="0",
(0,60) *+[red]{s^{(-13)}_{1^{(1)},1^{(4)},2^{(5)},1^{(4)}}} ="1",
(0,50) *+{s^{(-15)}_{2^{(2)},1^{(4)},2^{(5)},1^{(4)}}} ="2",
(0,40) *+{s^{(-17)}_{3^{(2)},1^{(4)},2^{(5)},1^{(4)}}} ="3",
(0,30) *+{s^{(-19)}_{4^{(2)},1^{(4)},2^{(5)},1^{(4)}}} ="4",
(0,20) *+{s^{(-21)}_{5^{(2)},1^{(4)},2^{(5)},1^{(4)}}}="5",
(0,10) *+{s^{(-23)}_{6^{(2)},1^{(4)},2^{(5)},1^{(4)}}}="6",
(0,0) *+{\vdots} ="7",
(25,75) *+{s^{(0)}_{1^{(2)}}} ="8",
(25,65) *+{s^{(-2)}_{2^{(2)}}} ="9",
(25,55) *+[red]{s^{(-12)}_{1^{(2)},1^{(4)},2^{(5)},1^{(4)}}} ="10",
(25,45) *+{s^{(-14)}_{2^{(2)},1^{(4)},2^{(5)},1^{(4)}}} ="11",
(25,35) *+{s^{(-16)}_{3^{(2)},1^{(4)},2^{(5)},1^{(4)}}} ="12",
(25,25) *+{s^{(-18)}_{4^{(2)},1^{(4)},2^{(5)},1^{(4)}}}="13",
(25,15) *+{s^{(-20)}_{5^{(2)},1^{(4)},2^{(5)},1^{(4)}}}="14",
(25,5)  *+{\vdots} ="15",
(50,70) *+{s^{(-1)}_{1^{(3)}}} ="16",
(50,60) *+{s^{(-3)}_{2^{(3)}}} ="17",
(50,50) *+[red]{s^{(-11)}_{1^{(3)},1^{(4)},2^{(5)},1^{(4)}}}="18",
(50,40) *+{s^{(-13)}_{2^{(3)},1^{(4)},2^{(5)},1^{(4)}}}="19",
(50,30) *+{s^{(-15)}_{3^{(3)},1^{(4)},2^{(5)},1^{(4)}}}="20",
(50,20) *+{s^{(-17)}_{4^{(3)},1^{(4)},2^{(5)},1^{(4)}}}="21",
(50,10) *+{s^{(-19)}_{5^{(3)},1^{(4)},2^{(5)},1^{(4)}}}="22",
(50,0)  *+{\vdots} ="23",
(75,75) *+{s^{(0)}_{1^{(4)}}} ="24",
(75,65) *+{s^{(-2)}_{2^{(4)}}} ="25",
(75,55) *+{s^{(-4)}_{3^{(4)}}} ="26",
(75,45) *+[red]{s^{(-8)}_{1^{(4)},2^{(5)},1^{(4)}}} ="27",
(75,35) *+{s^{(-10)}_{2^{(4)},2^{(5)},1^{(4)}}} ="28",
(75,25) *+{s^{(-12)}_{3^{(4)},2^{(5)},1^{(4)}}} ="29",
(75,15) *+{s^{(-14)}_{4^{(4)},2^{(5)},1^{(4)}}} ="30",
(75,5)  *+{\vdots} ="31",
(95,70) *+{s^{(-3)}_{1^{(5)},1^{(4)}}} ="32",
(95,60) *+{s^{(-5)}_{2^{(5)},1^{(4)}}} ="33",
(95,50) *+{s^{(-7)}_{3^{(5)},1^{(4)}}} ="34",
(95,40) *+{s^{(-9)}_{4^{(5)},1^{(4)}}} ="35",
(95,30) *+{s^{(-11)}_{5^{(5)},1^{(4)}}} ="36",
(95,20) *+{s^{(-13)}_{6^{(5)},1^{(4)}}} ="37",
(95,10) *+{s^{(-15)}_{7^{(5)},1^{(4)}}} ="38",
(95,0)  *+{\vdots} ="39",
(40,-10) *+{(i)} ="40",
"8", {\ar"0"},
"22", {\ar"0"},
"0", {\ar"23"},
"0", {\ar"31"},
"2", {\ar"1"},
"10", {\ar"1"},
"1", {\ar"11"},
"1", {\ar@/^1pc/"27"},
"3", {\ar"2"},
"11", {\ar"2"},
"2", {\ar"12"},
"4", {\ar"3"},
"5", {\ar"4"},
"12", {\ar"4"},
"4", {\ar"13"},
"6", {\ar"5"},
"13", {\ar"5"},
"5", {\ar"14"},
"7", {\ar"6"},
"14", {\ar"6"},
"6", {\ar"15"},
"15", {\ar"7"},
"7", {\ar"26"},
"33", {\ar"7"},
"9", {\ar"8"},
"8", {\ar"16"},
"16", {\ar"9"},
"9", {\ar"17"},
"23", {\ar"9"},
"11", {\ar"10"},
"18", {\ar"10"},
"10", {\ar"19"},
"12", {\ar"11"},
"19", {\ar"11"},
"11", {\ar"20"},
"13", {\ar"12"},
"20", {\ar"12"},
"12", {\ar"21"},
"14", {\ar"13"},
"21", {\ar"13"},
"13", {\ar"22"},
"15", {\ar"14"},
"22", {\ar"14"},
"14", {\ar"23"},
"15", {\ar"17"},
"23", {\ar"15"},
"26", {\ar"15"},
"17", {\ar"16"},
"24", {\ar"16"},
"16", {\ar"25"},
"17", {\ar@/^2pc/"23"},
"25", {\ar"17"},
"17", {\ar"26"},
"19", {\ar"18"},
"28", {\ar"18"},
"18", {\ar"29"},
"20", {\ar"19"},
"29", {\ar"19"},
"19", {\ar"30"},
"21", {\ar"20"},
"30", {\ar"20"},
"20", {\ar"31"},
"22", {\ar"21"},
"23", {\ar"22"},
"31", {\ar"22"},
"32", {\ar"24"},
"26", {\ar"25"},
"25", {\ar"32"},
"32", {\ar"26"},
"26", {\ar"33"},
"27", {\ar"33"},
"34", {\ar"27"},
"27", {\ar"35"},
"29", {\ar"28"},
"35", {\ar"28"},
"28", {\ar"36"},
"30", {\ar"29"},
"36", {\ar"29"},
"29", {\ar"37"},
"31", {\ar"30"},
"37", {\ar"30"},
"30", {\ar"38"},
"38", {\ar"31"},
"33", {\ar"32"},
"35", {\ar"34"},
"36", {\ar"35"},
"37", {\ar"36"},
"38", {\ar"37"},
"39", {\ar"38"},
\end{xy}}
\end{minipage}}
\caption{In type $A_5$, $(a)$ is the original quiver of type $A_5$ and $(b)$, $(c)$, $(d)$, $(e)$, $(f)$, $(g)$, $(h)$, $(i)$ are the quivers after we mutate $\mathfrak{r}(\mathscr{L}_{4,0})$, $\mathfrak{l}(\mathscr{L}_{5,-5})$, $\mathfrak{l}(\tau_l(\mathscr{L}_{5,-5}))$, $\mathfrak{l}(\tau^{2}_l(\mathscr{L}_{5,-5}))$, $\mathfrak{l}(\tau^{3}_l(\mathscr{L}_{5,-5}))$, $C_{3,-9}$, $C_{2,-10}$, $C_{1,-11}$ respectively.}\label{Q_{1}}
\end{figure}

We assign numbers to the vertices in the initial quiver in this example, see Figure \ref{Q1}. Then the mutation sequence (\ref{mutation sequence in example 6.1}) is
\begin{align}\label{Q1 mutation sequence in example 6.1}
\begin{split}
& c_{51},c_{52}, \ldots;  c_{12},c_{13}, \ldots;  c_{23},c_{24}, \ldots; c_{33},c_{34}, \ldots;  c_{44},c_{45}, \ldots; \\
& c_{12},c_{13}, \ldots;  c_{23},c_{24}, \ldots; c_{33},c_{34}, \ldots; c_{12},c_{13}, \ldots;  c_{23},c_{24}, \ldots;\\
& c_{12},c_{13}, \ldots;c_{33},c_{34}, \ldots;c_{23},c_{24}, \ldots;c_{12},c_{13}, \ldots.
\end{split}
\end{align}
Here $c_{51},c_{52}, \ldots$ (the other sequences in (\ref{Q1 mutation sequence in example 6.1}) have similar meaning) means we mutate the corresponding column starting from $c_{51}$. When we mutate $c_{51},c_{52}, \ldots$, it is enough to mutate sufficient many vertices (not to infinity).

(2) Let $S=2_{-18}4_{-14}5_{-11}5_{-9}4_{-6}$. By Definition \ref{define fundamental segment},
\begin{align*}
\mathcal{FS}(S)=\{4_{-6}, \ 5_{-11}5_{-9}4_{-6}, \ 2_{-18}4_{-14}5_{-11} \}.
\end{align*}
The set of distinguished factors of $S$ is $\{4_{-6}, \ 5_{-11}\}$. The mutation sequence for $S$ is
\begin{align} \label{mutation sequence in example 6.2}
\begin{split}
&C_{2,0}, \ C_{4,0}, \ C_{1,-1}, \ C_{3,-1}, \ C_{5,-1}, \ C_{2,-2}, \ C_{4,-2}, \ C_{1,-3}, \ C_{3,-3}, \ C_{5,-3}, \\
&C_{2,-4}, \ C_{4-4}, \ C_{1,-5}, \ C_{3,-5}, \ C_{5,-5}, \ \mathfrak{r}(\mathscr{L}_{4,-6}), \ \mathfrak{l}(\mathscr{L}_{5,-11}),\ \mathfrak{l}(\tau_l(\mathscr{L}_{5,-11})), \\
&\mathfrak{l}(\tau^{2}_l(\mathscr{L}_{5,-11})), \ \mathfrak{l}(\tau^{3}_l(\mathscr{L}_{5,-11})), \ C_{3,-15}, \ C_{2,-16}, \ C_{1,-17}.
\end{split}
\end{align}
We obtain the prime snake module $\mathcal{S}=L(S)$ at the vertex which has the label $(2,-18)$.

We assign numbers to the vertices in the initial quiver in this example, see Figure \ref{Q1}. Then the mutation sequence (\ref{mutation sequence in example 6.2}) is
\begin{align*}
& c_{21},c_{22}, \ldots;  c_{41},c_{42}, \ldots;  c_{11},c_{12}, \ldots; c_{31},c_{32}, \ldots;  c_{51},c_{52}, \ldots; \\
& c_{21},c_{22}, \ldots;  c_{41},c_{42}, \ldots;  c_{11},c_{12}, \ldots; c_{31},c_{32}, \ldots;  c_{51},c_{52}, \ldots; \\
& c_{21},c_{22}, \ldots;  c_{41},c_{42}, \ldots;  c_{11},c_{12}, \ldots; c_{31},c_{32}, \ldots;  c_{51},c_{52}, \ldots; \\
& c_{51},c_{52}, \ldots; c_{12},  c_{13}, \ldots; c_{23}, c_{24}, \ldots; c_{33},c_{34}, \ldots;c_{44}, c_{45},\ldots; \\
& c_{12},c_{13}, \ldots; c_{23},c_{24}, \ldots; c_{33},c_{34}, \ldots; c_{12},c_{13}, \ldots; c_{23},c_{24}, \ldots;\\
& c_{12},c_{13}, \ldots; c_{33},c_{34},\ldots; c_{23},c_{24},\ldots; c_{12},c_{13},\ldots.
\end{align*}
\end{example}

\begin{example}
Suppose that $U_q \widehat{\mathfrak{g}}$ is of type $A_{4}$. Let $\scalemath{0.9}{S=3_{-25}3_{-21}2_{-16}2_{-12}3_{-9}2_{-6}2_{-4}1_{-1}}$. Then
\begin{gather}
\begin{align*}
\mathcal{FS}(3_{-25}3_{-21}2_{-16}2_{-12}3_{-9}2_{-6}2_{-4}1_{-1})=\{1_{-1}, \ 3_{-9}2_{-6}2_{-4}1_{-1},\ 2_{-12}3_{-9},\ 2_{-16}2_{-12},\ 3_{-21}2_{-16},\ 3_{-25}3_{-21}\}.
\end{align*}
\end{gather}
The set of distinguished factors of $S$ is $\{1_{-1},\ 3_{-9},\ 2_{-12},\ 2_{-16},\ 3_{-21} \}$. The mutation sequence for $S$ is
\begin{align}\label{mutation sequence in example 6.3}
\begin{split}
&\mathfrak{r}(\mathscr{L}_{1, -1}), \ \mathfrak{r}(\tau_r(\mathscr{L}_{1, -1})), \ \mathfrak{r}(\tau^2_r(\mathscr{L}_{1, -1})), \ C_{3,-5}, \ C_{4,-6}, \ C_{3,-7}, \ C_{4,-8}, \\
&\mathfrak{l}(\mathscr{L}_{3, -9}), \ \mathfrak{l}(\tau_l(\mathscr{L}_{3, -9})), \ \mathscr{L}_{2, -12}, \ \mathfrak{r}(\mathscr{L}_{2, -16}), \ \tau_r(\mathscr{L}_{2, -16}), \ \mathscr{L}_{3, -21}.
\end{split}
\end{align}
We obtain the prime snake module $\mathcal{S}=L(S)$ at the vertex which has the label $(3,-25)$.

We assign numbers to the vertices in the initial quiver in this example, see Figure \ref{Q3}. Then the mutation sequence (\ref{mutation sequence in example 6.3}) is
\begin{align*}
& c_{41},c_{42},\ldots; c_{31},c_{32},\ldots; c_{22},c_{23},\ldots;c_{41},c_{42},\ldots; c_{31},c_{32},\ldots; c_{41},c_{42},\ldots;\\
& c_{31},c_{32},\ldots;c_{41},c_{42},\ldots; c_{31},c_{32},\ldots;c_{41},c_{42},\ldots;  c_{15},c_{16},\ldots; c_{25},c_{26},\ldots;\\
&c_{15},c_{16},\ldots;c_{15},c_{16},\ldots; c_{42}, c_{43},\ldots;  c_{33},c_{34},\ldots; c_{26},c_{27},\ldots;c_{43},c_{44},\ldots;\\
& c_{34},c_{35},\ldots; c_{16},c_{17},\ldots; c_{27}, c_{28},\ldots;c_{43},c_{44},\ldots; c_{34},c_{35},\ldots; c_{17},c_{18},\ldots;\\
& c_{28},c_{29},\ldots;c_{44},c_{45},\ldots; c_{35},c_{36},\ldots.
\end{align*}
\end{example}

\begin{example}
Suppose that $U_q \widehat{\mathfrak{g}}$ is of type $B_3$. (1) Let $S=1_{-35}2_{-29}2_{-21}3_{-16}3_{-10}2_{-5}$. Then
\begin{gather}
\begin{align*}
\mathcal{FS}(1_{-35}2_{-29}2_{-21}3_{-16}3_{-10}2_{-5})=\{2_{-5}, 3_{-10}2_{-5}, 3_{-16}3_{-10}, 2_{-21}3_{-16}, 2_{-29}2_{-21}, 1_{-35}2_{-29}\}.
\end{align*}
\end{gather}
The set of distinguished factors of $S$ is $\{2_{-5},\ 3_{-10},\ 3_{-16},\ 2_{-21},\ 2_{-29} \}$. The mutation sequence for $S$ is
\begin{align}\label{mutation sequence in example 6.41}
\begin{split}
&C_{3,0}, \ C_{2,-1}, \ C_{1,-3}, \ C_{3,-2}, \ C_{1,-1}, \ C_{2,-3}, \ \mathfrak{r}(\mathscr{L}_{2, -5}), \\
&\mathscr{L}_{3, -10}, \ \mathfrak{l}(\mathscr{L}_{3, -16}), \ \mathfrak{l}(\tau_l(\mathscr{L}_{3, -16})), \ \mathscr{L}_{2, -21}, \ \mathfrak{l}(\mathscr{L}_{2, -29}).
\end{split}
\end{align}
We obtain the prime snake module $\mathcal{S}=L(S)$ at the vertex which has the label $(1,-35)$.

We assign numbers to the vertices in the initial quiver in this example, see Figure \ref{Q4}. Then the mutation sequence (\ref{mutation sequence in example 6.41}) is
\begin{align*}
& c_{31},c_{32},\ldots; c_{21},c_{22},\ldots; c_{11},c_{12},\ldots; c_{31},c_{32},\ldots; c_{51},c_{52},\ldots; c_{41},c_{42},\ldots; \\
&c_{51},c_{52},\ldots;c_{41},c_{42},\ldots;c_{33},c_{34},\ldots;c_{12},c_{13},\ldots;c_{23},c_{24},\ldots; c_{51},c_{52},\ldots;\\
& c_{41},c_{42},\ldots;c_{35},c_{36},\ldots; c_{12},c_{13},\ldots;c_{23},c_{24},\ldots; c_{12},c_{13},\ldots;c_{53},c_{54},\ldots;\\
& c_{43},c_{44},\ldots; c_{39},c_{40},\ldots; c_{12}, c_{13},\ldots; c_{24},c_{25},\ldots;c_{13},c_{14},\ldots.
\end{align*}

(2) Let $S=2_{-43}2_{-35}2_{-31}1_{-25}3_{-18}3_{-8}3_{-2}3_{0}$. Then
\begin{align*}
\mathcal{FS}(S)=\{3_{-2}3_{0}, \ 3_{-8}3_{-2},\ 3_{-18}3_{-8},\ 1_{-25}3_{-18},\ 2_{-35}2_{-31}1_{-25},\ 2_{-43}2_{-35} \}.
\end{align*}
The set of distinguished factors of $S$ is $\{ 3_{0}, \ 3_{-2},\ 3_{-8},\ 3_{-18},\ 1_{-25},\ 2_{-35}\}$. The mutation sequence for $S$ is
\begin{align}\label{mutation sequence in example 6.42}
&\mathscr{L}_{3, -2},\ \mathscr{L}_{3, -8},\ \tau(\mathscr{L}_{3, -8}),\ \mathfrak{r}(\mathscr{L}_{3, -18}),\ \mathfrak{r}(\tau_r(\mathscr{L}_{3, -18})), \ \mathfrak{l}(\mathscr{L}_{1, -25}),\ \mathfrak{l}(\tau_l(\mathscr{L}_{1, -25})),\ \mathscr{L}_{2, -35}.
\end{align}
We obtain the prime snake module $\mathcal{S}=L(S)$ at the vertex which has the label $(2,-43)$.

We assign numbers to the vertices in the initial quiver in this example, see Figure \ref{Q4}. Then the mutation sequence (\ref{mutation sequence in example 6.42}) is
\begin{align*}
& c_{11},c_{12},\ldots; c_{22},c_{23},\ldots; c_{51},c_{52},\ldots; c_{41},c_{42},\ldots; c_{34},c_{35},\ldots; c_{11},c_{12},\ldots; \\
&c_{22},c_{23},\ldots;c_{52},c_{53},\ldots;c_{42},c_{43},\ldots;c_{36},c_{37},\ldots;c_{11},c_{12},\ldots; c_{22},c_{23},\ldots;\\
& c_{52},c_{53},\ldots;c_{42},c_{43},\ldots; c_{37},c_{38},\ldots;c_{52},c_{53},\ldots;
c_{42},c_{43},\ldots; c_{52},c_{53},\ldots; \\
& c_{13},c_{14},\ldots; c_{24},c_{25},\ldots; c_{311}, c_{312},\ldots;c_{43},c_{44},\ldots; c_{13},c_{14},\ldots; c_{24}, c_{25},\ldots; \\
& c_{312},c_{313},\ldots;c_{55},c_{56},\ldots; c_{14},c_{15},\ldots; c_{25},c_{26},\ldots; c_{315},c_{316},\ldots;c_{45},c_{46},\ldots.
\end{align*}

\begin{figure}[H]
\centering
\resizebox{.7\width}{.7\height}{
\begin{minipage}[b]{0.45\textwidth}
\centerline{
\begin{xy}
(10,70) *+{c_{11}} ="0",
(10,60) *+{c_{12}} ="1",
(10,50) *+{c_{13}} ="2",
(10,40) *+{c_{14}} ="3",
(10,30) *+{c_{15}} ="4",
(10,20) *+{c_{16}} ="5",
(10,10) *+{c_{17}} ="6",
(10,0) *+{\vdots} ="7",
(25,75) *+{c_{21}} ="8",
(25,65) *+{c_{22}} ="9",
(25,55) *+{c_{23}} ="10",
(25,45) *+{c_{24}} ="11",
(25,35) *+{c_{25}} ="12",
(25,25) *+{c_{26}} ="13",
(25,15) *+{c_{27}} ="14",
(25,5)  *+{\vdots} ="15",
(40,70) *+{c_{31}} ="16",
(40,60) *+{c_{32}} ="17",
(40,50) *+{c_{33}} ="18",
(40,40) *+{c_{34}} ="19",
(40,30) *+{c_{35}} ="20",
(40,20) *+{c_{36}} ="21",
(40,10) *+{c_{37}}="22",
(40,0)  *+{\vdots} ="23",
(55,75) *+{c_{41}} ="24",
(55,65) *+{c_{42}} ="25",
(55,55) *+{c_{43}} ="26",
(55,45) *+{c_{44}} ="27",
(55,35) *+{c_{45}} ="28",
(55,25) *+{c_{46}} ="29",
(55,15) *+{c_{47}} ="30",
(55,5)  *+{\vdots} ="31",
(70,70) *+{c_{51}} ="32",
(70,60) *+{c_{52}} ="33",
(70,50) *+{c_{53}} ="34",
(70,40) *+{c_{54}} ="35",
(70,30) *+{c_{55}} ="36",
(70,20) *+{c_{56}} ="37",
(70,10) *+{c_{57}} ="38",
(70,0)  *+{\vdots} ="39",
"1", {\ar"0"},
"8", {\ar"0"},
"0", {\ar"9"},
"2", {\ar"1"},
"9", {\ar"1"},
"1", {\ar"10"},
"3", {\ar"2"},
"10", {\ar"2"},
"2", {\ar"11"},
"4", {\ar"3"},
"11", {\ar"3"},
"3", {\ar"12"},
"5", {\ar"4"},
"12", {\ar"4"},
"4", {\ar"13"},
"6", {\ar"5"},
"13", {\ar"5"},
"5", {\ar"14"},
"7", {\ar"6"},
"14", {\ar"6"},
"6", {\ar"15"},
"15", {\ar"7"},
"9", {\ar"8"},
"8", {\ar"16"},
"10", {\ar"9"},
"16", {\ar"9"},
"9", {\ar"17"},
"11", {\ar"10"},
"17", {\ar"10"},
"10", {\ar"18"},
"12", {\ar"11"},
"18", {\ar"11"},
"11", {\ar"19"},
"13", {\ar"12"},
"19", {\ar"12"},
"12", {\ar"20"},
"14", {\ar"13"},
"20", {\ar"13"},
"13", {\ar"21"},
"15", {\ar"14"},
"21", {\ar"14"},
"14", {\ar"22"},
"22", {\ar"15"},
"15", {\ar"23"},
"17", {\ar"16"},
"24", {\ar"16"},
"16", {\ar"25"},
"18", {\ar"17"},
"25", {\ar"17"},
"17", {\ar"26"},
"19", {\ar"18"},
"26", {\ar"18"},
"18", {\ar"27"},
"20", {\ar"19"},
"27", {\ar"19"},
"19", {\ar"28"},
"21", {\ar"20"},
"28", {\ar"20"},
"20", {\ar"29"},
"22", {\ar"21"},
"29", {\ar"21"},
"21", {\ar"30"},
"23", {\ar"22"},
"30", {\ar"22"},
"22", {\ar"31"},
"31", {\ar"23"},
"25", {\ar"24"},
"24", {\ar"32"},
"26", {\ar"25"},
"32", {\ar"25"},
"25", {\ar"33"},
"27", {\ar"26"},
"33", {\ar"26"},
"26", {\ar"34"},
"28", {\ar"27"},
"34", {\ar"27"},
"27", {\ar"35"},
"29", {\ar"28"},
"35", {\ar"28"},
"28", {\ar"36"},
"30", {\ar"29"},
"36", {\ar"29"},
"29", {\ar"37"},
"31", {\ar"30"},
"37", {\ar"30"},
"30", {\ar"38"},
"38", {\ar"31"},
"31", {\ar"39"},
"33", {\ar"32"},
"34", {\ar"33"},
"35", {\ar"34"},
"36", {\ar"35"},
"37", {\ar"36"},
"38", {\ar"37"},
"39", {\ar"38"},
\end{xy}}
\caption{The initial quiver in type $A_5$.}\label{Q1}
\end{minipage}
\begin{minipage}[b]{0.45\linewidth}
\centerline{
\begin{xy}
(10,70) *+{c_{11}} ="0",
(10,60) *+{c_{12}} ="1",
(10,50) *+{c_{13}} ="2",
(10,40) *+{c_{14}} ="3",
(10,30) *+{c_{15}} ="4",
(10,20) *+{c_{16}} ="5",
(10,10) *+{c_{17}} ="6",
(10,0) *+{\vdots} ="7",
(25,75) *+{c_{21}} ="8",
(25,65) *+{c_{22}} ="9",
(25,55) *+{c_{23}} ="10",
(25,45) *+{c_{24}} ="11",
(25,35) *+{c_{25}} ="12",
(25,25) *+{c_{26}} ="13",
(25,15) *+{c_{27}} ="14",
(25,5)  *+{\vdots} ="15",
(40,70) *+{c_{31}} ="16",
(40,60) *+{c_{32}} ="17",
(40,50) *+{c_{33}} ="18",
(40,40) *+{c_{34}} ="19",
(40,30) *+{c_{35}} ="20",
(40,20) *+{c_{36}} ="21",
(40,10) *+{c_{37}} ="22",
(40,0)  *+{\vdots} ="23",
(55,75) *+{c_{41}} ="24",
(55,65) *+{c_{42}} ="25",
(55,55) *+{c_{43}} ="26",
(55,45) *+{c_{44}} ="27",
(55,35) *+{c_{45}} ="28",
(55,25) *+{c_{46}} ="29",
(55,15) *+{c_{47}} ="30",
(55,5)  *+{\vdots} ="31",
"1", {\ar"0"},
"8", {\ar"0"},
"0", {\ar"9"},
"2", {\ar"1"},
"9", {\ar"1"},
"1", {\ar"10"},
"3", {\ar"2"},
"10", {\ar"2"},
"2", {\ar"11"},
"4", {\ar"3"},
"11", {\ar"3"},
"3", {\ar"12"},
"5", {\ar"4"},
"12", {\ar"4"},
"4", {\ar"13"},
"6", {\ar"5"},
"13", {\ar"5"},
"5", {\ar"14"},
"7", {\ar"6"},
"14", {\ar"6"},
"6", {\ar"15"},
"15", {\ar"7"},
"9", {\ar"8"},
"8", {\ar"16"},
"10", {\ar"9"},
"16", {\ar"9"},
"9", {\ar"17"},
"11", {\ar"10"},
"17", {\ar"10"},
"10", {\ar"18"},
"12", {\ar"11"},
"18", {\ar"11"},
"11", {\ar"19"},
"13", {\ar"12"},
"19", {\ar"12"},
"12", {\ar"20"},
"14", {\ar"13"},
"20", {\ar"13"},
"13", {\ar"21"},
"15", {\ar"14"},
"21", {\ar"14"},
"14", {\ar"22"},
"22", {\ar"15"},
"15", {\ar"23"},
"17", {\ar"16"},
"24", {\ar"16"},
"16", {\ar"25"},
"18", {\ar"17"},
"25", {\ar"17"},
"17", {\ar"26"},
"19", {\ar"18"},
"26", {\ar"18"},
"18", {\ar"27"},
"20", {\ar"19"},
"27", {\ar"19"},
"19", {\ar"28"},
"21", {\ar"20"},
"28", {\ar"20"},
"20", {\ar"29"},
"22", {\ar"21"},
"29", {\ar"21"},
"21", {\ar"30"},
"23", {\ar"22"},
"30", {\ar"22"},
"22", {\ar"31"},
"31", {\ar"23"},
"25", {\ar"24"},
"26", {\ar"25"},
"27", {\ar"26"},
"28", {\ar"27"},
"29", {\ar"28"},
"30", {\ar"29"},
"31", {\ar"30"},
\end{xy}}
\caption{The initial quiver in type $A_4$.}\label{Q3}
\end{minipage}
\begin{minipage}[b]{0.45\linewidth}
\centerline{
\begin{xy}
(0,100) *+{c_{11}} ="2",
(0,80) *+{c_{12}} ="3",
(0,60) *+{c_{13}} ="4",
(0,40) *+{c_{14}} ="5",
(0,20) *+{c_{15}} ="6",
(0,0) *+{\vdots} ="7",
(15,110) *+{c_{21}} ="9",
(15,90) *+{c_{22}} ="10",
(15,70) *+{c_{23}} ="11",
(15,50) *+{c_{24}} ="12",
(15,30) *+{c_{25}} ="13",
(15,10) *+{\vdots} ="14",
(30,120) *+{c_{31}} ="18",
(30,110) *+{c_{32}} ="19",
(30,100) *+{c_{33}} ="20",
(30,90) *+{c_{34}} ="21",
(30,80) *+{c_{35}} ="22",
(30,70) *+{c_{36}} ="23",
(30,60) *+{c_{37}} ="24",
(30,50) *+{c_{38}} ="25",
(30,40) *+{c_{39}} ="26",
(30,30) *+{c_{310}} ="27",
(30,20) *+{c_{311}} ="28",
(30,10) *+{c_{312}} ="29",
(30,0) *+{\vdots} ="30",
(45,100) *+{c_{41}} ="34",
(45,80) *+{c_{42}} ="35",
(45,60) *+{c_{43}} ="36",
(45,40) *+{c_{44}} ="37",
(45,20) *+{c_{45}} ="38",
(45,0) *+{\vdots} ="39",
(60,110) *+{c_{51}} ="41",
(60,90) *+{c_{52}} ="42",
(60,70) *+{c_{53}} ="43",
(60,50) *+{c_{54}} ="44",
(60,30) *+{c_{55}} ="45",
(60,10) *+{\vdots} ="46",
"3", {\ar"2"},
"9", {\ar"2"},
"2", {\ar"10"},
"4", {\ar"3"},
"10", {\ar"3"},
"3", {\ar"11"},
"5", {\ar"4"},
"11", {\ar"4"},
"4", {\ar"12"},
"6", {\ar"5"},
"7", {\ar"6"},
"12", {\ar"5"},
"5", {\ar"13"},
"13", {\ar"6"},
"6", {\ar"14"},
"10", {\ar"9"},
"18", {\ar"9"},
"9", {\ar"20"},
"11", {\ar"10"},
"20", {\ar"10"},
"10", {\ar"22"},
"12", {\ar"11"},
"22", {\ar"11"},
"11", {\ar"24"},
"13", {\ar"12"},
"24", {\ar"12"},
"12", {\ar"26"},
"14", {\ar"13"},
"14", {\ar"7"},
"26", {\ar"13"},
"13", {\ar"28"},
"28", {\ar"14"},
"14", {\ar"30"},
"19", {\ar"18"},
"20", {\ar"19"},
"19", {\ar"34"},
"21", {\ar"20"},
"22", {\ar"21"},
"34", {\ar"21"},
"21", {\ar"35"},
"23", {\ar"22"},
"24", {\ar"23"},
"35", {\ar"23"},
"23", {\ar"36"},
"25", {\ar"24"},
"26", {\ar"25"},
"36", {\ar"25"},
"25", {\ar"37"},
"27", {\ar"26"},
"28", {\ar"27"},
"37", {\ar"27"},
"27", {\ar"38"},
"29", {\ar"28"},
"30", {\ar"29"},
"38", {\ar"29"},
"29", {\ar"39"},
"35", {\ar"34"},
"41", {\ar"34"},
"34", {\ar"42"},
"36", {\ar"35"},
"42", {\ar"35"},
"35", {\ar"43"},
"37", {\ar"36"},
"43", {\ar"36"},
"36", {\ar"44"},
"38", {\ar"37"},
"44", {\ar"37"},
"37", {\ar"45"},
"39", {\ar"38"},
"45", {\ar"38"},
"38", {\ar"46"},
"46", {\ar"39"},
"42", {\ar"41"},
"43", {\ar"42"},
"44", {\ar"43"},
"45", {\ar"44"},
"46", {\ar"45"},
\end{xy}}
\caption{The initial quiver in type $B_3$.}\label{Q4}
\end{minipage}}
\end{figure}
\end{example}

\section{Proofs of Theorem \ref{real snake modules}}\label{prove theorem of Section 4}
In this section, we prove Theorem \ref{real snake modules}.

\subsection{Proof of Theorem \ref{real snake modules}}
Let $\mathcal{S}$ be a prime snake module and $S$ its highest $l$-weight monomial. Then $\mathcal{S}$ can be written as
\begin{align*}
\mathcal{S}=\mathcal{S}^{(t)}_{k_{1}^{(i_{1},j_{1})}, k_{2}^{(i_{2},j_{2})}, \ldots, k_{m-1}^{(i_{m-1},j_{m-1})}, k_{m}^{(i_{m})}},
\end{align*}
where $m \geq 1$, $j_\ell \geq 0$, $1 \leq \ell \leq m-1$, if $j_{\ell}=0$, then $i_{\ell} \neq i_{\ell+1}$, $k_1, k_2, \ldots, k_m \in \mathbb{Z}_{\geq 1}$, $t \in \mathbb{Z}$.

The theorem follows from the fact that $\chi_{q}(\mathcal{S})\chi_{q}(\mathcal{S})$ has only one dominant monomial $S^2$.

Let $L=\sum_{\ell=1}^{m} k_{\ell}$. Suppose that $\mathfrak{m}=\prod_{\ell=1}^{L}m(p_{\ell})$ (respectively, $\mathfrak{m}'=\prod_{\ell=1}^{L}m(p'_{\ell})$) is a monomial in the first (respectively, the second) $\chi_{q}(\mathcal{S})$ in $\chi_{q}(\mathcal{S})\chi_{q}(\mathcal{S})$, where
$(p_{1}, \ldots, p_{L}) \in \overline{\mathscr{P}}_{(c_{\ell},d_{\ell})_{1 \leq \ell \leq L}}$ (respectively, $(p'_{1}, \ldots, p'_{L}) \in \overline{\mathscr{P}}_{(c_{\ell},d_{\ell})_{1 \leq \ell \leq L}}$) is a tuple of non-overlapping paths and
\begin{align*}
c_{r+\sum_{\ell=1}^{j-1}k_{\ell}}=i_{j}, \quad d_{r+\sum_{\ell=1}^{j-1}k_{\ell}}=t+2d_{i_j}r-2d_{i_j}+\sum\nolimits_{\ell=1}^{j-1} n_\ell,
\end{align*}
where $1 \leq r \leq k_j$, $1 \leq j \leq m$, by convention $\sum_{\ell=1}^{0}k_{\ell}=0$.

Suppose that $\mathfrak{m}\mathfrak{m}'$ is dominant. If $p_{L}\neq p^{+}_{c_{L},d_{L}}$, then $\mathfrak{m}\mathfrak{m}'$ is right-negative and not dominant.
Therefore $p_{L} = p^{+}_{c_{L},d_{L}}$. Similarly, we have $p'_{L} = p'^{+}_{c'_{L},d'_{L}}$. By the non-overlapping property, we have $p'_{j}=p'^{+}_{c'_{j},d'_{j}}$, $p_{j}=p^{+}_{c_{j},d_{j}}$ for all $1+\sum_{\ell=1}^{m-1} k_{\ell} \leq j \leq L$.

Suppose that $p_{\sum_{\ell=1}^{m-1} k_{\ell}} \neq p^{+}_{c_{\sum_{\ell=1}^{m-1} k_{\ell}},d_{\sum_{\ell=1}^{m-1} k_{\ell}}}$. Then $m(p_{\sum_{\ell=1}^{m-1} k_{\ell}})$ has some negative factor $i_\ell$, where $(i,\ell) \in C^{-}_{p_{\sum_{\ell=1}^{m-1} k_{\ell}}}$. By Theorem \ref{path description of q-characters}, $\mathfrak{m}$ has the negative factor $i^{-1}_{\ell}$. Therefore, the negative factor $i^{-1}_{\ell}$ is canceled by $\mathfrak{m}'$. It follows that $\mathfrak{m}'\neq S$ since $i_{\ell}$ is not in $S$. But then $\mathfrak{m}\mathfrak{m}'$ has one of the following factors:
\begin{gather}
\begin{align*}
\begin{split}
\text{Type $A_n$: } & 1^{-1}_{\ell+i-1}, \ 2^{-1}_{\ell+i-2}, \ \ldots,  \ (i-1)^{-1}_{\ell+1}, (i+1)^{-1}_{\ell+1}, \ \ldots, \ (n-1)^{-1}_{\ell+n-i-1}, \ (n)^{-1}_{\ell+n-i}, \\
\text{Type $B_n$: } & i \neq n, \ \scalemath{0.92}{1^{-1}_{\ell+2i-2}, \ 2^{-1}_{\ell+2i-4}, \ \ldots, \ (i-1)^{-1}_{\ell+2}, \ (i+1)^{-1}_{\ell+2},\ \ldots, \ (n-1)^{-1}_{\ell+2n-2i-2}, \ (n)^{-1}_{\ell+2n-2i-1},} \\
&i= n, \ 1^{-1}_{\ell+2n-3}, \ 2^{-1}_{\ell+2n-5}, \ \ldots, \ (n-2)^{-1}_{\ell+3}, \ (n-1)^{-1}_{\ell+1}.
\end{split}
\end{align*}
\end{gather}
This contradicts the assumption that $\mathfrak{m}\mathfrak{m}'$ is dominant. Therefore, $p_{\sum_{\ell=1}^{m-1} k_{\ell}}= p^{+}_{c_{\sum_{\ell=1}^{m-1} k_{\ell}},d_{\sum_{\ell=1}^{m-1} k_{\ell}}}$.
By Theorem \ref{path description of q-characters}, we have $p_{j}=p^{+}_{c_{j},d_{j}}$ for all $1+\sum_{\ell=1}^{m-2} k_{\ell} \leq  j \leq \sum_{\ell=1}^{m-1} k_{\ell}$. By the same argument, we have $p_{j}=p^{+}_{c_{j},d_{j}}$, $1 \leq j \leq \sum_{\ell=1}^{m-2} k_{\ell}$. Therefore, $\mathfrak{m}=S$.

By a similar argument, we may show that $\mathfrak{m}'=S$. Therefore, the only dominant monomial in $\chi_{q}(\mathcal{S})\chi_{q}(\mathcal{S})$ is $S^2$.

\section{Proof of Theorem \ref{S-systems}} \label{proof of S-systems}
In this section, we prove Theorem \ref{S-systems}.

\subsection{Classification of dominant monomials}
First we classify all dominant monomials in each summand on the left- and right-hand sides of every equation in Theorem \ref{S-systems}. We have the following lemma.

\begin{lemma}\label{dominant monomials}
Let $[\mathcal{S}_{1}][\mathcal{S}_{2}]=[\mathcal{S}_{3}][\mathcal{S}_{4}]+[\mathcal{S}_{5}][\mathcal{S}_{6}]$ be any equation in the $S$-system of type $A_n$ (respectively, $B_n$) in Theorem \ref{S-systems}. Let $S_i$ be the highest $l$-weight monomial of $\mathcal{S}_i$, $i\in \{1,2,\ldots,6\}$. The dominant monomials in each summand on the left- and right-hand sides of $[\mathcal{S}_{1}][\mathcal{S}_{2}]=[\mathcal{S}_{3}][\mathcal{S}_{4}]+[\mathcal{S}_{5}][\mathcal{S}_{6}]$ are given in Table \ref{dominant monomials in the S-system of type A} (respectively, Table \ref{dominant monomials in the S-system of type B}).
\end{lemma}

In Table \ref{dominant monomials in the S-system of type A} and Table \ref{dominant monomials in the S-system of type B}, $M\prod_{0\leq j\leq r}A^{-1}_{i,s} = M$ for $r = -1$, $s \in \mathbb{Z}$; ``DMs'' means ``Dominant monomials''.

\subsection{Proof of Theorem \ref{S-systems}}
By Tables \ref{dominant monomials in the S-system of type A} and \ref{dominant monomials in the S-system of type B}, the dominant monomials of the $q$-characters of the left-hand side and of the right-hand side of every equation in Theorem \ref{S-systems} are the same and have the same multiplicities (the monomials occurring in the $q$-character of snake modules have multiplicity one, see {\cite[Theorem 6.1]{MY12a}}, {\cite[Theorem 6.5]{MY12b}} or Theorem \ref{path description of q-characters}). Therefore by Proposition \ref{dominant monomials determine q-characters}, the theorem is true.
\begin{table}[H]
\resizebox{.8\width}{.8\height}{
\begin{tabular}{c c c c}
\hline %
$\mathcal{S}_1$ & DMs of $\chi_q(\mathcal{S}_1)\chi_q(\mathcal{S}_2)$ & DMs of $\chi_q(\mathcal{S}_3)\chi_q(\mathcal{S}_4)$ & DMs of $\chi_q(\mathcal{S}_5)\chi_q(\mathcal{S}_6)$ \\
\hline %
(1) in Table \ref{definition S_1 in type A_{n}} & $\substack{S_1S_2 \prod_{0 \leq j \leq r} A^{-1}_{i_{1}, t+2k_{1}-2j-1}, \\ -1 \leq r \leq k_{1}-1}$ & $\substack{S_3S_4 \prod_{0 \leq j \leq r} A^{-1}_{i_{1}, t+2k_1-2j-1}, \\-1 \leq r \leq k_{1}-2}$ & $\substack{S_5 S_6}$  \\
\hline %
(2) in Table \ref{definition S_1 in type A_{n}} & $\substack{S_1S_2 \prod_{0 \leq j \leq r} A^{-1}_{i_{1}, t+2k_{1}-2j-1}, \\ -1 \leq r \leq k_{1}-1}$ & $\substack{S_3S_4 \prod_{0 \leq j \leq r} A^{-1}_{i_{1}, t+2k_1-2j-1}, \\-1 \leq r \leq k_{1}-2}$ & $\substack{S_5 S_6}$  \\
\hline %
(3) in Table \ref{definition S_1 in type A_{n}} & $\substack{S_1S_2 \prod_{0 \leq j \leq r} A^{-1}_{i_{1}, t+2k_{1}-2j-1}, \\ -1 \leq r \leq k_{1}-1}$ & $\substack{S_3S_4 \prod_{0 \leq j \leq r} A^{-1}_{i_{1}, t+2k_1-2j-1}, \\-1 \leq r \leq k_{1}-2}$ & $\substack{S_5 S_6}$  \\
\hline %
(4) in Table \ref{definition S_1 in type A_{n}} & $\substack{S_1S_2 \prod_{0 \leq j \leq r} A^{-1}_{i_{1}, t+2k_{1}-2j-1}, \\ -1 \leq r \leq k_{1}-1}$ & $\substack{S_3S_4 \prod_{0 \leq j \leq r} A^{-1}_{i_{1}, t+2k_1-2j-1}, \\-1 \leq r \leq k_{1}-2}$ & $\substack{S_5 S_6}$  \\
\hline %
(5) in Table \ref{definition S_1 in type A_{n}} & $\substack{S_1S_2 \prod_{0 \leq j \leq r} A^{-1}_{i_{1}, t+2k_{1}-2j-1}, \\ -1 \leq r \leq k_{1}-1}$ & $\substack{S_3S_4 \prod_{0 \leq j \leq r} A^{-1}_{i_{1}, t+2k_1-2j-1}, \\-1 \leq r \leq k_{1}-2}$ & $\substack{S_5 S_6}$  \\
\hline %
(6) in Table \ref{definition S_1 in type A_{n}} & $\substack{S_1S_2 \prod_{0 \leq j \leq r} A^{-1}_{i_{1}, t+2k_{1}-2j-1}, \\ -1 \leq r \leq k_{1}-1}$ & $\substack{S_3S_4 \prod_{0 \leq j \leq r} A^{-1}_{i_{1}, t+2k_1-2j-1}, \\-1 \leq r \leq k_{1}-2}$ & $\substack{S_5 S_6}$  \\
\hline %
(7) in Table \ref{definition S_1 in type A_{n}} & $\substack{S_1S_2 \prod_{0 \leq j \leq r} A^{-1}_{i_{1}, t+2k_{1}-2j-1}, \\ -1 \leq r \leq k_{1}-1}$ & $\substack{S_3S_4 \prod_{0 \leq j \leq r} A^{-1}_{i_{1}, t+2k_1-2j-1}, \\-1 \leq r \leq k_{1}-2}$ & $\substack{S_5 S_6}$  \\
\hline %
\end{tabular}}
\caption{Classification of dominant monomials in the $S$-system of type $A_n$.}\label{dominant monomials in the S-system of type A}
\end{table}

\begin{table}[H]
\resizebox{.7\width}{.7\height}{
\begin{tabular}{c c c c}
\hline %
$\mathcal{S}_1$ & DMs of $\chi_q(\mathcal{S}_1)\chi_q(\mathcal{S}_2)$ & DMs of $\chi_q(\mathcal{S}_3)\chi_q(\mathcal{S}_4)$ & DMs of $\chi_q(\mathcal{S}_5)\chi_q(\mathcal{S}_6)$ \\
\hline %
(1) in Table \ref{definition S_1 in type B_{n}} & $\substack{S_1S_2 \prod_{0 \leq j \leq r} A^{-1}_{i_{1}, t+2d_{i_1}k_{1}-2d_{i_1}j-d_{i_1}},\\ -1 \leq r \leq k_{1}-1}$ & $\substack{S_3S_4 \prod_{0 \leq j \leq r} A^{-1}_{i_{1}, t+2d_{i_1}k_{1}-2d_{i_1}j-d_{i_1}}, \\ -1 \leq r \leq k_{1}-2}$ & $\substack{S_5 S_6}$  \\
\hline %
(2) in Table \ref{definition S_1 in type B_{n}} & $\substack{S_1S_2 \prod_{0 \leq j \leq r} A^{-1}_{i_{1}, t+4k_{1}-4j-2},\\ -1 \leq r \leq k_{1}-1}$ & $\substack{S_3S_4 \prod_{0 \leq j \leq r} A^{-1}_{i_{1}, t+4k_{1}-4j-2}, \\ -1 \leq r \leq k_{1}-2}$ & $\substack{S_5 S_6}$  \\
\hline %
(3) in Table \ref{definition S_1 in type B_{n}} & $\substack{S_1S_2 \prod_{0 \leq j \leq r} A^{-1}_{i_{1}, t+4k_{1}-4j-2},\\ -1 \leq r \leq k_{1}-1}$ & $\substack{S_3S_4 \prod_{0 \leq j \leq r} A^{-1}_{i_{1}, t+4k_{1}-4j-2}, \\ -1 \leq r \leq k_{1}-2}$ & $\substack{S_5 S_6}$  \\
\hline %
(4) in Table \ref{definition S_1 in type B_{n}} & $\substack{S_1S_2 \prod_{0 \leq j \leq r} A^{-1}_{i_{1}, t+4k_{1}-4j-2},\\ -1 \leq r \leq k_{1}-1}$ & $\substack{S_3S_4 \prod_{0 \leq j \leq r} A^{-1}_{i_{1}, t+4k_{1}-4j-2}, \\ -1 \leq r \leq k_{1}-2}$ & $\substack{S_5 S_6}$ \\
\hline %
(5) in Table \ref{definition S_1 in type B_{n}} & $\substack{S_1S_2 \prod_{0 \leq j \leq r} A^{-1}_{i_{1}, t+4k_{1}-4j-2},\\ -1 \leq r \leq k_{1}-1}$ & $\substack{S_3S_4 \prod_{0 \leq j \leq r} A^{-1}_{i_{1}, t+4k_{1}-4j-2}, \\ -1 \leq r \leq k_{1}-2}$ & $\substack{S_5 S_6}$  \\
\hline %
(6) in Table \ref{definition S_1 in type B_{n}} & $\substack{S_1S_2 \prod_{0 \leq j \leq r} A^{-1}_{i_{1}, t+4k_{1}-4j-2},\\ -1 \leq r \leq k_{1}-1}$ & $\substack{S_3S_4 \prod_{0 \leq j \leq r} A^{-1}_{i_{1}, t+4k_{1}-4j-2}, \\ -1 \leq r \leq k_{1}-2}$ & $\substack{S_5 S_6}$ \\
\hline %
(7) in Table \ref{definition S_1 in type B_{n}} & $\substack{S_1S_2 \prod_{0 \leq j \leq r} A^{-1}_{i_{1}, t+4k_{1}-4j-2},\\ -1 \leq r \leq k_{1}-1}$ & $\substack{S_3S_4 \prod_{0 \leq j \leq r} A^{-1}_{i_{1}, t+4k_{1}-4j-2}, \\ -1 \leq r \leq k_{1}-2}$ & $\substack{S_5 S_6}$\\
\hline %
(8) in Table \ref{definition S_1 in type B_{n}} & $\substack{S_1S_2 \prod_{0 \leq j \leq r} A^{-1}_{i_{1}, t+4k_{1}-4j-2},\\ -1 \leq r \leq k_{1}-1}$ & $\substack{S_3S_4 \prod_{0 \leq j \leq r} A^{-1}_{i_{1}, t+4k_{1}-4j-2}, \\ -1 \leq r \leq k_{1}-2}$ & $\substack{S_5 S_6}$  \\
\hline %
(9) in Table \ref{definition S_1 in type B_{n}} & $\substack{S_1S_2 \prod_{0 \leq j \leq r} A^{-1}_{i_{1}, t+2k_{1}-2j-1},\\ -1 \leq r \leq k_{1}-1}$ & $\substack{S_3S_4 \prod_{0 \leq j \leq r} A^{-1}_{i_{1}, t+2k_{1}-2j-1}, \\ -1 \leq r \leq k_{1}-2}$ & $\substack{S_5 S_6}$  \\
\hline %
(10) in Table \ref{definition S_1 in type B_{n}} & $\substack{S_1S_2 \prod_{0 \leq j \leq r} A^{-1}_{i_{1}, t+2k_{1}-2j-1},\\ -1 \leq r \leq k_{1}-1}$ & $\substack{S_3S_4 \prod_{0 \leq j \leq r} A^{-1}_{i_{1}, t+2k_{1}-2j-1}, \\ -1 \leq r \leq k_{1}-2}$ & $\substack{S_5 S_6}$ \\
\hline %
(11) in Table \ref{definition S_1 in type B_{n}} & $\substack{S_1S_2 \prod_{0 \leq j \leq r} A^{-1}_{i_{1}, t+2k_{1}-2j-1},\\ -1 \leq r \leq k_{1}-1}$ & $\substack{S_3S_4 \prod_{0 \leq j \leq r} A^{-1}_{i_{1}, t+2k_{1}-2j-1}, \\ -1 \leq r \leq k_{1}-2}$ & $\substack{S_5 S_6}$  \\
\hline %
(12) in Table \ref{definition S_1 in type B_{n}} & $\substack{S_1S_2 \prod_{0 \leq j \leq r} A^{-1}_{i_{1}, t+2k_{1}-2j-1},\\ -1 \leq r \leq k_{1}-1}$ & $\substack{S_3S_4 \prod_{0 \leq j \leq r} A^{-1}_{i_{1}, t+2k_{1}-2j-1}, \\ -1 \leq r \leq k_{1}-2}$ & $\substack{S_5 S_6}$  \\
\hline %
(13) in Table \ref{definition S_1 in type B_{n}} & $\substack{S_1S_2 \prod_{0 \leq j \leq r} A^{-1}_{i_{1}, t+4k_{1}-4j-2},\\ -1 \leq r \leq k_{1}-1}$ & $\substack{S_3S_4 \prod_{0 \leq j \leq r} A^{-1}_{i_{1}, t+4k_{1}-4j-2}, \\ -1 \leq r \leq k_{1}-2}$ & $\substack{S_5 S_6}$  \\
\hline %
(14) in Table \ref{definition S_1 in type B_{n}}  & $\substack{S_1S_2 \prod_{0 \leq j \leq r} A^{-1}_{i_{1}, t+4k_{1}-4j-2},\\ -1 \leq r \leq k_{1}-1}$ & $\substack{S_3S_4 \prod_{0 \leq j \leq r} A^{-1}_{i_{1}, t+4k_{1}-4j-2}, \\ -1 \leq r \leq k_{1}-2}$ & $\substack{S_5 S_6}$    \\
\hline %
(15) in Table \ref{definition S_1 in type B_{n}} & $\substack{S_1S_2 \prod_{0 \leq j \leq r} A^{-1}_{i_{1}, t+4k_{1}-4j-2},\\ -1 \leq r \leq k_{1}-1}$ & $\substack{S_3S_4 \prod_{0 \leq j \leq r} A^{-1}_{i_{1}, t+4k_{1}-4j-2}, \\ -1 \leq r \leq k_{1}-2}$ & $\substack{S_5 S_6}$  \\
\hline %
(16) in Table \ref{definition S_1 in type B_{n}}  & $\substack{S_1S_2 \prod_{0 \leq j \leq r} A^{-1}_{i_{1}, t+4k_{1}-4j-2},\\ -1 \leq r \leq k_{1}-1}$ & $\substack{S_3S_4 \prod_{0 \leq j \leq r} A^{-1}_{i_{1}, t+4k_{1}-4j-2}, \\ -1 \leq r \leq k_{1}-2}$ & $\substack{S_5 S_6}$  \\
\hline %
(17) in Table \ref{definition S_1 in type B_{n}}  & $\substack{S_1S_2 \prod_{0 \leq j \leq r} A^{-1}_{i_{1}, t+4k_{1}-4j-2},\\ -1 \leq r \leq k_{1}-1}$ & $\substack{S_3S_4 \prod_{0 \leq j \leq r} A^{-1}_{i_{1}, t+4k_{1}-4j-2}, \\ -1 \leq r \leq k_{1}-2}$ & $\substack{S_5 S_6}$  \\
\hline %
(18) in Table \ref{definition S_1 in type B_{n}} & $\substack{S_1S_2 \prod_{0 \leq j \leq r} A^{-1}_{i_{1}, t+2k_{1}-2j-1},\\ -1 \leq r \leq k_{1}-1}$ & $\substack{S_3S_4 \prod_{0 \leq j \leq r} A^{-1}_{i_{1}, t+2k_{1}-2j-1}, \\ -1 \leq r \leq k_{1}-2}$ & $\substack{S_5 S_6}$  \\
\hline %
\end{tabular}}
\caption{Classification of dominant monomials in the $S$-system of type $B_n$.}\label{dominant monomials in the S-system of type B}
\end{table}

\subsection{Proof of Lemma \ref{dominant monomials}}
We will prove the lemma for the case corresponding to (5) in Table \ref{definition S_1 in type A_{n}} in type $A_n$ and the cases corresponding to (6), (11), (14), (18) in Table \ref{definition S_1 in type B_{n}} in type $B_n$. The other cases are similar.

\begin{proof}[\bf Proof of the case corresponding to (5) in Table \ref{definition S_1 in type A_{n}}.] \label{proof classification of dominant monomials $A_{n}$}
Let
\begin{align*}
\mathcal{S}_{1}= \mathcal{S}^{(t+2)}_{k_{1}^{(i_{1})}, {k_{2}}^{(i_{2}+1)}, k_{3}^{(i_{3},j_{3})}, k_{4}^{(i_{4},j_{4})}, \ldots, k_{m}^{(i_{m})}},\quad
\mathcal{S}_{2}= \mathcal{S}^{(t)}_{k_{1}^{(i_{1})}, k_{2}^{(i_{2})}, k_{3}^{(i_{3},j_{3})}, k_{4}^{(i_{4},j_{4})}, \ldots, k_{m}^{(i_{m})}},
\end{align*}
where $m\geq 3$, $i_1 > i_2$, $i_3 > i_2$, $j_\ell \geq 0$, $3\leq \ell \leq m-1$. We will firstly prove that the dominant monomials in $\chi_q(\mathcal{S}_1)\chi_q(\mathcal{S}_2)$ must be of the form $S_1 \chi_q(\mathcal{S}_2)$.

Let $L=\sum_{\ell=1}^{m} k_{\ell}$. Suppose that $\mathfrak{m}=\prod_{j=1}^{L}m(p_{j})$ is a monomial in $\chi_{q}(\mathcal{S}_{2})$, where
$(p_{1}, \ldots, p_{L}) \in \overline{\mathscr{P}}_{(c_{j},d_{j})_{1 \leq j \leq L}}$ is a tuple of non-overlapping paths and
\begin{align*}
c_{r+\sum_{\ell=1}^{j-1}k_{\ell}}=i_{j}, \quad d_{r+\sum_{\ell=1}^{j-1}k_{\ell}}=t+2r-2+\sum\nolimits_{\ell=1}^{j-1}n_\ell,
\end{align*}
where $1 \leq r \leq k_j$, $1 \leq j \leq m$, by convention $\sum_{\ell=1}^{0}k_{\ell}=0$.

Let $\mathfrak{m}'=\prod^{L}_{j=1}m(p'_{j})$ be a monomial in $\chi_{q}(\mathcal{S}_{1})$, where $(p'_{1}, \ldots, p'_{L}) \in \overline{\mathscr{P}}_{(c'_{j},d'_{j})_{1\leq j \leq L}}$ is a tuple of non-overlapping paths and
\begin{align*}
&c'_{1}=c'_{2}=\cdots=c'_{k_{1}}=i_{1}, \ c'_{k_{1}+1}=c'_{k_{1}+2}=\cdots=c'_{k_{1}+k_{2}}=i_{2}+1, \\
&d'_{1}=t+2,\ d'_{2}=t+4,\ \ldots,\ d'_{k_{1}}=t+2k_{1}, \\
&d'_{k_{1}+1}=t+n_1+1,\ d'_{k_{1}+2}=t+n_1+3,\ \ldots,\ d'_{k_{1}+k_{2}}=t+n_1+2k_{2}-1,\\
&c'_{r+\sum_{\ell=1}^{j-1}k_{\ell}}=i_{j},\quad d'_{r+\sum_{\ell=1}^{j-1}k_{\ell}}=t+2r-2+\sum\nolimits_{\ell=1}^{j-1}n_\ell,
\end{align*}
where $1 \leq r \leq k_j$, $3 \leq j \leq m$. Thus we have $c_{\ell}=c_{\ell'}$, $d_{\ell}=d_{\ell'}$ for $k_1+k_2+1\leq \ell \leq L$, $c'_{\ell}=c_{\ell}+1$, $d'_{\ell}=d_{\ell}+1$ for $k_1+1\leq \ell \leq k_1+k_2$, and $c'_{\ell}=c_{\ell}$, $d'_{\ell}=d_{\ell}+2$ for $1\leq \ell \leq k_1$.

Suppose that $\mathfrak{m}\mathfrak{m}'$ is dominant. By the same arguments as the arguments in the proof of Theorem \ref{real snake modules}, we have $p_{j}=p^{+}_{c_{j},d_{j}}$, $k_{1}+1 \leq j \leq L$ and $p'_{j}=p'^{+}_{c'_{j},d'_{j}}$, $1 \leq j \leq L$. Therefore, the dominant monomials in $\chi_q(\mathcal{S}_1)\chi_q(\mathcal{S}_2)$ must be of the form $S_1 \chi_q(\mathcal{S}_2)$.

If $p_{k_{1}}=p^{+}_{c_{k_{1}},d_{k_{1}}}$, then $p_{j}=p^{+}_{c_{j},d_{j}}$ for all $1 \leq j \leq k_{1}-1$. Therefore, $\mathfrak{m}\mathfrak{m}'$=$S_{1}S_{2}$. If $p_{k_{1}}=p^{+}_{c_{k_{1}},d_{k_{1}}}A^{-1}_{i_{1}, t+2k_{1}-1}$, then $p_{j} \in \{ p^{+}_{c_{j},d_{j}}, \ p^{+}_{c_{j},d_{j}}A^{-1}_{i_{1}, t+2j-1 } \}$, $1 \leq j \leq k_{1}-1$. Therefore, $\mathfrak{m}\mathfrak{m}'$ is one of the dominant monomials $S_1S_2\prod_{j=0}^{r}A^{-1}_{i_{1}, t+2k_{1}-2j-1}$, $0 \leq r \leq k_{1}-1$. If $p_{k_{1}} \not\in \{ p^{+}_{c_{k_{1}},d_{k_{1}}}, \ p^{+}_{c_{k_{1}},d_{k_{1}}}A^{-1}_{i_{1}, t+2k_{1}-1} \}$, then by the same arguments as the arguments in the proof of Theorem \ref{real snake modules}, it follows that $\mathfrak{m}\mathfrak{m}'$ is not dominant which contradicts our assumption.

By the same arguments as the arguments of dominant monomials in $\chi_q(\mathcal{S}_1)\chi_q(\mathcal{S}_2)$, the dominant monomials in $\chi_q(\mathcal{S}_3)\chi_q(\mathcal{S}_4)$ must be of the form $\scalemath{0.92}{S_3S_4\prod_{j=0}^{r}A^{-1}_{i_{1}, t+2k_{1}-2j-1}}$, $-1 \leq r \leq k_{1}-2$, and the dominant monomial in $\chi_q(\mathcal{S}_5)\chi_q(\mathcal{S}_6)$ is $S_5S_6$.
\end{proof}

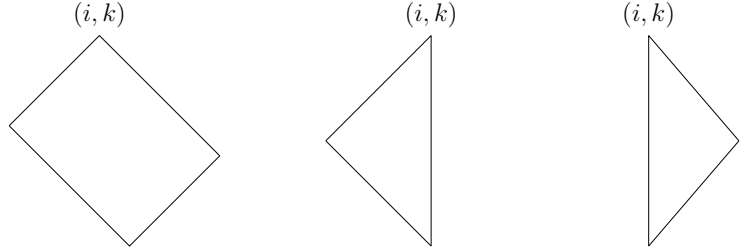
\begin{figure}[H]
\centering
\resizebox{0.8\width}{0.8\height}{
\begin{minipage}[b]{0.3\textwidth}
\centerline{
\begin{tikzpicture}
\draw (0,0)--(2,-2) (0,0)--(-1.5,-1.5);
\draw (-1.5,-1.5)--(0.5,-3.5)--(2,-2);
\node [above] at (0,0) {$(i,k)$};
\end{tikzpicture}}
\end{minipage}
\begin{minipage}[b]{0.3\textwidth}
\centerline{
\begin{tikzpicture}
\draw  (0.25,-0.75)--(-1.5,-2.5) (0.25,-0.75)--(0.25,-4.25);
\draw (-1.5,-2.5)--(0.25,-4.25);
\node[above] at (0.25,-0.75) {$(i,k)$};
\end{tikzpicture}}
\end{minipage}
\begin{minipage}[b]{0.3\textwidth}
\centerline{
\begin{tikzpicture}
\draw  (0.25,-0.75)--(1.75,-2.5) (0.25,-0.75)--(0.25,-4.25);
\draw (1.75,-2.5)--(0.25,-4.25);
\node[above] at (0.25,-0.75) {$(i,k)$};
\end{tikzpicture}}
\end{minipage}}
\caption{The rectangular box or triangle corresponding to $\mathscr{P}_{(i,k)}$.}\label{rectangular box or a triangle}
\end{figure}

Recall that each $\mathscr{P}_{(i,k)}$, $i\in I$, $k\in \mathbb{Z}$, defined in Section \ref{Path description of q-characters}, corresponds to a rectangular box or a triangle, see Figure \ref{rectangular box or a triangle}. The monomials appearing in $\chi_q(i_k)$ correspond to paths in the rectangular box or the triangle. In the following proofs, we will use frequently the arguments in the proof of Theorem \ref{real snake modules} and the fact that snake modules are special (Theorem \ref{path description of q-characters}).

\begin{proof}[\bf Proof of the case corresponding to (6) in Table \ref{definition S_1 in type B_{n}}.] Suppose that $\mathcal{S}_1$ is of the form (6) in Table \ref{definition S_1 in type B_{n}}.
Let
\begin{align*}
&\mathcal{S}_{1}=\mathcal{S}^{(t+4)}_{k_{1}^{(i_{1})}, {k_{2}}^{(i_{2}+\sgn(i_1-i_2),j_2-1)}, k_{3}^{(i_{3},j_{3})}, k_{4}^{(i_{4},j_{4})}, \ldots, k_{m}^{(i_{m})}},
&\mathcal{S}_{2}= \mathcal{S}^{(t)}_{{k_{1}}^{(i_1)}, k_{2}^{(i_{2},j_{2})}, k_{3}^{(i_{3},j_{3})}, \ldots, k_{m}^{(i_{m})}},
\end{align*}
where $i_1, i_2 \neq n$, $j_1=0$, $j_2 \geq 1$, $m\geq 3$, $(i_2-i_1)(i_3-i_2)\geq 0$, and $j_\ell \geq 0$, $3 \leq \ell \leq m-1$.

Let $L=\sum_{\ell=1}^{m} k_{\ell}$. Suppose that $\mathfrak{m}=\prod_{j=1}^{L}m(p_{j})$ (respectively, $\mathfrak{m}'=\prod^{L}_{\ell=1}m(p'_{j})$) is a monomial in $\chi_{q}(\mathcal{S}_{2})$ (respectively, $\chi_{q}(\mathcal{S}_{1})$), where
$(p_{1}, \ldots, p_{L}) \in \overline{\mathscr{P}}_{(c_{j},d_{j})_{1 \leq j \leq L}}$ (respectively, $(p'_{1}, \ldots, p'_{L}) \in \overline{\mathscr{P}}_{(c'_{j},d'_{j})_{1\leq j \leq L}}$) is a tuple of non-overlapping paths, where $c_{\ell}=c_{\ell'}$, $d_{\ell}=d_{\ell'}$ for $k_1+k_2+1\leq \ell \leq L$, $c'_{\ell}=c_{\ell}+\sgn(i_1-i_2)$, $d'_{\ell}=d_{\ell}+2$ for $k_1+1\leq \ell \leq k_1+k_2$, and $c'_{\ell}=c_{\ell}$, $d'_{\ell}=d_{\ell}+4$ for $1\leq \ell \leq k_1$. That is, $\mathscr{P}_{(c_{\ell},d_{\ell})}=\mathscr{P}_{(c'_{\ell},d'_{\ell})}$ for $k_1+k_2+1\leq \ell \leq L$, see (a)--(c) in Figure \ref{proof 1}, relations between $\mathscr{P}_{(c_{\ell},d_{\ell})}$ and $\mathscr{P}_{(c'_{\ell},d'_{\ell})}$ for $k_1+1\leq \ell \leq k_1+k_2$ are shown in (d), (e) in Figure \ref{proof 1}, and relations between $\mathscr{P}_{(c_{\ell},d_{\ell})}$ and $\mathscr{P}_{(c'_{\ell},d'_{\ell})}$ for $1\leq \ell \leq k_1$ are shown in (f) in Figure \ref{proof 1}.
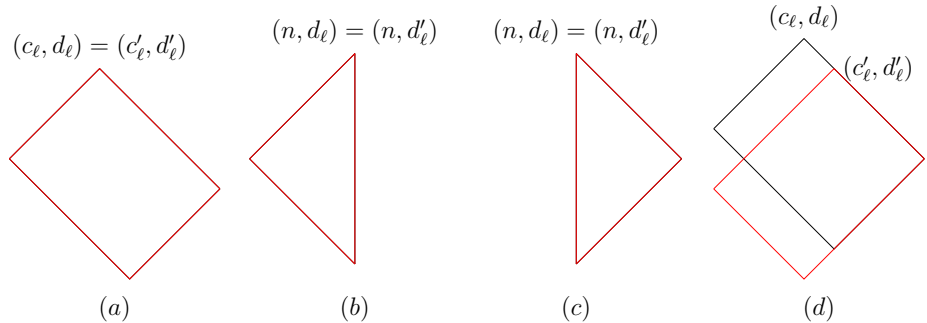
\begin{figure}[H]
\centering
\resizebox{0.8\width}{0.8\height}{
\begin{minipage}[b]{0.25\textwidth}
\centerline{
\begin{tikzpicture}
\draw (0,0)--(2,-2) (0,0)--(-1.5,-1.5);
\draw (-1.5,-1.5)--(0.5,-3.5)--(2,-2);
\draw[red] (0,0)--(2,-2) (0,0)--(-1.5,-1.5);
\draw[red] (-1.5,-1.5)--(0.5,-3.5)--(2,-2);
\node [above] at (0,0) {$(c_\ell,d_\ell)=(c'_\ell,d'_\ell)$};
\node at (0.25,-4) {$(a)$};
\end{tikzpicture}}
\end{minipage}
\begin{minipage}[b]{0.25\textwidth}
\centerline{
\begin{tikzpicture}
\draw (0.25,-0.25)--(-1.5,-2) (0.25,-0.25)--(0.25,-3.75);
\draw[red] (0.25,-0.25)--(-1.5,-2) (0.25,-0.25)--(0.25,-3.75);
\draw (-1.5,-2)--(0.25,-3.75);
\draw[red] (-1.5,-2)--(0.25,-3.75);
\node [above] at (0.25,-0.25) {$(n,d_\ell)=(n,d'_\ell)$};
\node at (0.25,-4.5) {$(b)$};
\end{tikzpicture}}
\end{minipage}
\begin{minipage}[b]{0.25\textwidth}
\centerline{
\begin{tikzpicture}
\draw (0.25,-0.25)--(2,-2) (0.25,-0.25)--(0.25,-3.75);
\draw[red]  (0.25,-0.25)--(2,-2) (0.25,-0.25)--(0.25,-3.75);
\draw (2,-2)--(0.25,-3.75);
\draw[red] (2,-2)--(0.25,-3.75);
\node [above] at (0.25,-0.25) {$(n,d_\ell)=(n,d'_\ell)$};
\node at (0.25,-4.5) {$(c)$};
\end{tikzpicture}}
\end{minipage}
\begin{minipage}[b]{0.25\textwidth}
\centerline{
\begin{tikzpicture}
\draw (0,0)--(2,-2) (0,0)--(-1.5,-1.5);
\draw (-1.5,-1.5)--(0.5,-3.5) -- (2,-2);
\draw [red] (0.5,-0.5)--(2,-2) (0.5,-0.5)--(-1.5,-2.5);
\draw[red] (-1.5,-2.5)--(0,-4) -- (2,-2);
\node [above] at (0,0) {$(c_\ell,d_\ell)$};
\node [right] at (0.5,-0.5) {$(c'_\ell,d'_\ell)$};
\node at (0.25,-4.5) {$(d)$};
\end{tikzpicture}}
\end{minipage}}
\end{figure}

\begin{figure}[H]
\centering
\resizebox{0.8\width}{0.8\height}{
\begin{minipage}[b]{0.3\textwidth}
\centerline{
\begin{tikzpicture}
\draw (0,0)--(1.5,-1.5) (0,0)--(-2,-2);
\draw (1.5,-1.5)--(-0.5,-3.5) -- (-2,-2);
\draw[red] (-0.5,-0.5)--(1.5,-2.5) (-0.5,-0.5)--(-2,-2);
\draw[red] (1.5,-2.5)--(0,-4) -- (-2,-2);
\node [above] at (0,0) {$(c_\ell,d_\ell)$};
\node [left] at (-0.5,-0.5) {$(c'_\ell,d'_\ell)$};
\node at (-0.25,-4.5) {$(e)$};
\end{tikzpicture}}
\end{minipage}
\begin{minipage}[b]{0.3\textwidth}
\centerline{
\begin{tikzpicture}
\draw (0,0)--(1.5,-1.5) (0,0)--(-1.5,-1.5);
\draw (-1.5,-1.5)--(0,-3) -- (1.5,-1.5);
\draw[red] (0,-0.8)--(1.5,-2.3) (0,-0.8)--(-1.5,-2.3);
\draw[red] (1.5,-2.3)--(0,-3.8) -- (-1.5,-2.3);
\node [above] at (0,0) {$(c_\ell,d_\ell)$};
\node [above] at (0,-1) {$(c'_\ell,d'_\ell)$};
\node at (0,-4.5) {$(f)$};
\end{tikzpicture}}
\end{minipage}}
\caption{Relations between $\mathscr{P}_{(c_{\ell},d_{\ell})}$ (black) and $\mathscr{P}_{(c'_{\ell},d'_{\ell})}$ (red).}\label{proof 1}
\end{figure}

Suppose that $\mathfrak{m}\mathfrak{m}'$ is dominant. By the same arguments as the arguments in the proof of Theorem \ref{real snake modules} and relations between $\mathscr{P}_{(c_{\ell},d_{\ell})}$ and $\mathscr{P}_{(c'_{\ell},d'_{\ell})}$, we have $p_{j}=p^{+}_{c_{j},d_{j}}$, $k_{1}+1 \leq j \leq L$ and $p'_{j}=p'^{+}_{c'_{j},d'_{j}}$, $1 \leq j \leq L$. Therefore, the dominant monomials in $\chi_q(\mathcal{S}_1)\chi_q(\mathcal{S}_2)$ must be of the form $S_1 \chi_q(\mathcal{S}_2)$.

If $p_{k_{1}}=p^{+}_{c_{k_{1}},d_{k_{1}}}$, then $p_{j}=p^{+}_{c_{j},d_{j}}$ for all $1 \leq j \leq k_{1}-1$. Therefore, $\mathfrak{m}\mathfrak{m}' = S_{1}S_{2}$. If $p_{k_{1}}=p^{+}_{c_{k_{1}},d_{k_{1}}}A^{-1}_{i_{1}, t+4k_{1}-2}$, then $p_{j} \in \{ p^{+}_{c_{j},d_{j}}, \ p^{+}_{c_{j},d_{j}}A^{-1}_{i_{1}, t+4j-2} \}$, $1 \leq j \leq k_{1}-1$. Therefore, $\mathfrak{m}\mathfrak{m}'$ is one of the dominant monomials $S_1S_2\prod_{j=0}^{r}A^{-1}_{i_{1}, t+4k_{1}-4j-2}$, $0 \leq r \leq k_{1}-1$. If $p_{k_{1}} \not\in \{ p^{+}_{c_{k_{1}},d_{k_{1}}}, \ p^{+}_{c_{k_{1}},d_{k_{1}}}A^{-1}_{i_{1}, t+4k_{1}-2} \}$, then by the same arguments as the arguments in the proof of Theorem \ref{real snake modules}, it follows that $\mathfrak{m}\mathfrak{m}'$ is not dominant which contradicts our assumption.

By the same arguments as the arguments of dominant monomials in $\chi_q(\mathcal{S}_1)\chi_q(\mathcal{S}_2)$, the dominant monomials in $\chi_q(\mathcal{S}_3)\chi_q(\mathcal{S}_4)$ must be of the form $\scalemath{0.92}{S_3S_4\prod_{j=0}^{r}A^{-1}_{i_{1}, t+4k_{1}-4j-2}}$, $-1 \leq r \leq k_{1}-2$, and the dominant monomial in $\chi_q(\mathcal{S}_5)\chi_q(\mathcal{S}_6)$ is $S_5S_6$.
\end{proof}

\begin{proof}[\bf Proof of the case corresponding to (11) in Table \ref{definition S_1 in type B_{n}}.] Suppose that $\mathcal{S}_1$ is of the form (11) in Table \ref{definition S_1 in type B_{n}}. Let
\begin{align*}
&\mathcal{S}_{1}=\mathcal{S}^{(t+2)}_{(k_{1}+1)^{(n)}, (\frac{d_{i_{2}}}{d_{i_{2}+1}}k_{2})^{(i_{2}+1,j_{2}-1)}, k_{3}^{(i_{3},j_{3})}, k_{4}^{(i_{4},j_{4})}, \ldots, k_{m}^{(i_{m})}},
&\mathcal{S}_{2}= \mathcal{S}^{(t)}_{{k_{1}}^{(n)}, k_{2}^{(i_{2},j_{2})}, k_{3}^{(i_{3},j_{3})}, \ldots, k_{m}^{(i_{m})}},
\end{align*}
where $j_1=0$, $j_2 \geq 1$, $m\geq 3$, $i_2 \geq i_3$, and $j_\ell \geq 0$ for $3 \leq \ell \leq m-1$.

Let $L=\sum_{\ell=1}^{m} k_{\ell}$ and $L'=k_{1}+1+(\frac{d_{i_{2}}}{d_{i_{2}+1}}k_{2})+\sum_{\ell=3}^{m} k_{\ell}$.  Suppose that $\mathfrak{m}=\prod_{j=1}^{L}m(p_{j})$ (respectively, $\mathfrak{m}'=\prod^{L'}_{\ell=1}m(p'_{j})$) is a monomial in $\chi_{q}(\mathcal{S}_{2})$ (respectively, $\chi_{q}(\mathcal{S}_{1})$), where
$(p_{1}, \ldots, p_{L}) \in \overline{\mathscr{P}}_{(c_{j},d_{j})_{1 \leq j \leq L}}$ (respectively, $(p'_{1}, \ldots, p'_{L'}) \in \overline{\mathscr{P}}_{(c'_{j},d'_{j})_{1\leq j \leq L'}}$) is a tuple of non-overlapping paths, where $c_{k_1+k_2+\ell}=c'_{k_{1}+1+(\frac{d_{i_{2}}}{d_{i_{2}+1}}k_{2})+\ell}$, $d_{k_1+k_2+\ell}=d'_{k_{1}+1+(\frac{d_{i_{2}}}{d_{i_{2}+1}}k_{2})+\ell}$ for $1 \leq \ell \leq \sum_{\ell=3}^{m} k_{\ell}$, $c'_{k_1+1+\ell}=i_2+1$ for $1 \leq \ell \leq (\frac{d_{i_{2}}}{d_{i_{2}+1}}k_{2})$, $c_{k_1+\ell}=i_2$ for $1 \leq \ell \leq k_{2}$, and $c'_{\ell}=c_{\ell}=c'_{k_1+1}=n$, $d'_{\ell}=d_{\ell}+2$ for $1\leq \ell \leq k_1$. That is, $\mathscr{P}_{(c_{k_1+k_2+\ell},d_{k_1+k_2+\ell})}=\mathscr{P}_{(c'_{k_{1}+1+(\frac{d_{i_{2}}}
{d_{i_{2}+1}}k_{2})+\ell},d'_{k_{1}+1+(\frac{d_{i_{2}}}{d_{i_{2}+1}}k_{2})+\ell})}$ for $1 \leq \ell \leq \sum_{\ell=3}^{m} k_{\ell}$, see (a)--(c) in Figure \ref{proof 2}, relations between $\mathscr{P}_{(i_2,d_{k_1+\ell})}$ for $1 \leq \ell \leq k_{2}$ and $\mathscr{P}_{(i_2+1,d'_{k_1+1+\ell})}$ for $1 \leq \ell \leq (\frac{d_{i_{2}}}{d_{i_{2}+1}}k_{2})$ are shown in (d)--(f) in Figure \ref{proof 2}, and relations between  $\mathscr{P}_{(n,d_{\ell})}$ and $\mathscr{P}_{(n,d'_{\ell})}$ for $1\leq \ell \leq k_1$ are shown in (g), (h) in Figure \ref{proof 2}.
\begin{figure}[H]
\centering
\resizebox{0.8\width}{0.8\height}{
\begin{minipage}[b]{0.25\textwidth}
\centerline{
\begin{tikzpicture}
\draw (0,0)--(2,-2) (0,0)--(-1.5,-1.5);
\draw (-1.5,-1.5)--(0.5,-3.5)--(2,-2);
\draw[red] (0,0)--(2,-2) (0,0)--(-1.5,-1.5);
\draw[red] (-1.5,-1.5)--(0.5,-3.5)--(2,-2);
\node [above] at (0,0) {$(c_\ell,d_\ell)=(c'_\ell,d'_\ell)$};
\node at (0.25,-4) {$(a)$};
\end{tikzpicture}}
\end{minipage}
\begin{minipage}[b]{0.25\textwidth}
\centerline{
\begin{tikzpicture}
\draw (0.25,-0.25)--(-1.5,-2) (0.25,-0.25)--(0.25,-3.75);
\draw[red] (0.25,-0.25)--(-1.5,-2) (0.25,-0.25)--(0.25,-3.75);
\draw (-1.5,-2)--(0.25,-3.75);
\draw[red] (-1.5,-2)--(0.25,-3.75);
\node [above] at (0.25,-0.25) {$(n,d_\ell)=(n,d'_\ell)$};
\node at (0.25,-4.5) {$(b)$};
\end{tikzpicture}}
\end{minipage}
\begin{minipage}[b]{0.25\textwidth}
\centerline{
\begin{tikzpicture}
\draw (0.25,-0.25)--(2,-2) (0.25,-0.25)--(0.25,-3.75);
\draw[red]  (0.25,-0.25)--(2,-2) (0.25,-0.25)--(0.25,-3.75);
\draw (2,-2)--(0.25,-3.75);
\draw[red] (2,-2)--(0.25,-3.75);
\node [above] at (0.25,-0.25) {$(n,d_\ell)=(n,d'_\ell)$};
\node at (0.25,-4.5) {$(c)$};
\end{tikzpicture}}
\end{minipage}
\begin{minipage}[b]{0.25\textwidth}
\centerline{
\begin{tikzpicture}
\draw (0,0)--(2,-2) (0,0)--(-1.5,-1.5);
\draw (-1.5,-1.5)--(0.5,-3.5)--(2,-2);
\draw[red] (0.25,-0.25)--(2,-2) (0.25,0.25)--(-1.5,-1.5) (0.25,0.25)--(0.25,-3.25);
\draw[red] (2,-2)--(0.25,-3.75) (-1.5,-1.5)--(0.25,-3.25) (0.25,-0.25)--(0.25,-3.75);
\node [left] at (0,0) {$(n-1,d_\ell)$};
\node [above] at (0.25,0.25) {$(n,d_\ell-1)$};
\node [right] at (0.25,-0.25) {$(n,d_\ell+1)$};
\node at (0.25,-4.5) {$(d)$};
\end{tikzpicture}}
\end{minipage}}
\end{figure}

\begin{figure}[H]
\centering
\resizebox{0.8\width}{0.8\height}{
\begin{minipage}[b]{0.25\textwidth}
\centerline{
\begin{tikzpicture}
\draw (0,0)--(-2,-2) (0,0)--(1.5,-1.5);
\draw (1.5,-1.5)--(-0.5,-3.5)--(-2,-2);
\draw[red] (-0.25,-0.25)--(-2,-2) (-0.25,0.25)--(1.5,-1.5) (-0.25,0.25)--(-0.25,-3.25);
\draw[red] (-2,-2)--(-0.25,-3.75) (1.5,-1.5)--(-0.25,-3.25) (-0.25,-0.25)--(-0.25,-3.75);
\node [right] at (0,0) {$(n-1,d_\ell)$};
\node [above] at (0.25,0.25) {$(n,d_\ell-1)$};
\node [left] at (-0.25,-0.25) {$(n,d_\ell+1)$};
\node at (0.25,-4.5) {$(e)$};
\end{tikzpicture}}
\end{minipage}
\begin{minipage}[b]{0.25\textwidth}
\centerline{
\begin{tikzpicture}
\draw (0,0)--(2,-2) (0,0)--(-1.5,-1.5);
\draw (-1.5,-1.5)--(0.5,-3.5) -- (2,-2);
\draw [red] (0.5,-0.5)--(2,-2) (0.5,-0.5)--(-1.5,-2.5);
\draw[red] (-1.5,-2.5)--(0,-4) -- (2,-2);
\node [above] at (0,0) {$(c_\ell,d_\ell)$};
\node [right] at (0.5,-0.5) {$(c'_\ell,d'_\ell)$};
\node at (0.25,-4.5) {$(f)$};
\end{tikzpicture}}
\end{minipage}
\begin{minipage}[b]{0.25\textwidth}
\centerline{
\begin{tikzpicture}
\draw (0.25,-0.25)--(-1.5,-2) (0.25,-0.25)--(0.25,-3.75);
\draw[red] (0.25,-0.75)--(2,-2.5) (0.25,-0.75)--(0.25,-4.25);
\draw (-1.5,-2)--(0.25,-3.75);
\draw[red] (2,-2.5)--(0.25,-4.25);
\node [above] at (0.25,-0.25) {$(n,d_\ell)$};
\node [right] at (0.25,-0.75) {$(n,d'_\ell)$};
\node at (0.25,-4.5) {$(g)$};
\end{tikzpicture}}
\end{minipage}
\begin{minipage}[b]{0.25\textwidth}
\centerline{
\begin{tikzpicture}
\draw (0.25,-0.25)--(2,-2) (0.25,-0.25)--(0.25,-3.75);
\draw[red]  (0.25,-0.75)--(-1.5,-2.5) (0.25,-0.75)--(0.25,-4.25);
\draw (2,-2)--(0.25,-3.75);
\draw[red] (-1.5,-2.5)--(0.25,-4.25);
\node [above] at (0.25,-0.25) {$(n,d_\ell)$};
\node [left] at (0.25,-0.75) {$(n,d'_\ell)$};
\node at (0.25,-4.5) {$(h)$};
\end{tikzpicture}}
\end{minipage}}
\caption{Relations between $\mathscr{P}_{(c_{\ell},d_{\ell})}$ (black) and $\mathscr{P}_{(c'_{\ell},d'_{\ell})}$ (red).}\label{proof 2}
\end{figure}

Suppose that $\mathfrak{m}\mathfrak{m}'$ is dominant. By the same arguments as the arguments in the proof of Theorem \ref{real snake modules} and relations between $\mathscr{P}_{(c_{\ell},d_{\ell})}$ and $\mathscr{P}_{(c'_{\ell},d'_{\ell})}$, we have $p_{j}=p^{+}_{c_{j},d_{j}}$, $k_{1}+1 \leq j \leq L$ and $p'_{j}=p'^{+}_{c'_{j},d'_{j}}$, $1 \leq j \leq L'$. Therefore, the dominant monomials in $\chi_q(\mathcal{S}_1)\chi_q(\mathcal{S}_2)$ must be of the form $S_1 \chi_q(\mathcal{S}_2)$.

If $p_{k_{1}}=p^{+}_{c_{k_{1}},d_{k_{1}}}$, then $p_{j}=p^{+}_{c_{j},d_{j}}$ for all $1 \leq j \leq k_{1}-1$. Therefore, $\mathfrak{m}\mathfrak{m}' = S_{1}S_{2}$. If $p_{k_{1}}=p^{+}_{c_{k_{1}},d_{k_{1}}}A^{-1}_{i_{1}, t+2k_{1}-1}$, then $p_{j} \in \{ p^{+}_{c_{j},d_{j}}, \ p^{+}_{c_{j},d_{j}}A^{-1}_{i_{1}, t+2j-1} \}$, $1 \leq j \leq k_{1}-1$. Therefore, $\mathfrak{m}\mathfrak{m}'$ is one of the dominant monomials $S_1S_2\prod_{j=0}^{r}A^{-1}_{i_{1}, t+2k_{1}-2j-1}$, $0 \leq r \leq k_{1}-1$. If $p_{k_{1}} \not\in \{ p^{+}_{c_{k_{1}},d_{k_{1}}}, \ p^{+}_{c_{k_{1}},d_{k_{1}}}A^{-1}_{i_{1}, t+2k_{1}-1} \}$, then by the same arguments as the arguments in the proof of Theorem \ref{real snake modules}, it follows that $\mathfrak{m}\mathfrak{m}'$ is not dominant which contradicts our assumption.

By the same arguments as the arguments of dominant monomials in $\chi_q(\mathcal{S}_1)\chi_q(\mathcal{S}_2)$, the dominant monomials in $\chi_q(\mathcal{S}_3)\chi_q(\mathcal{S}_4)$ must be of the form $\scalemath{0.92}{S_3S_4\prod_{j=0}^{r}A^{-1}_{i_{1}, t+2k_{1}-2j-1}}$, $-1 \leq r \leq k_{1}-2$, and the dominant monomial in $\chi_q(\mathcal{S}_5)\chi_q(\mathcal{S}_6)$ is $S_5S_6$.
\end{proof}

\begin{proof}[\bf Proof of the case corresponding to (14) in Table \ref{definition S_1 in type B_{n}}.] Suppose that $\mathcal{S}_1$ is of the form (14) in Table \ref{definition S_1 in type B_{n}}. Let
\begin{gather}
\begin{align*}
&\mathcal{S}_{1}= \mathcal{S}^{(t+4)}_{k_{1}^{(i_{1})}, (\frac{k_{2}-1}{2})^{(n-1)}, 1^{(n)}, (\frac{d_{i_{3}}}{d_{i_{3}+1}}k_{3})^{(i_{3}+1,j_{3}-\delta_{i_3i_4})}, k_{4}^{(i_{4},j_{4})},\ldots, k_{m}^{(i_{m})}},
&\mathcal{S}_{2}= \mathcal{S}^{(t)}_{{k_{1}}^{(i_1)}, k_{2}^{(n)}, k_{3}^{(i_{3},j_{3})}, k_{4}^{(i_{4},j_{4})},\ldots, k_{m}^{(i_{m})}},
\end{align*}
\end{gather}
where $k_2$ is odd, $m\geq 4$, $i_3\leq i_4$, and $j_\ell \geq 0$ for $3\leq \ell \leq m-1$.

Let $L=\sum_{\ell=1}^{m} k_{\ell}$ and $L'=k_1+(\frac{k_{2}-1}{2})+(\frac{d_{i_{3}}}{d_{i_{3}+1}}k_{3})+1+\sum_{\ell=4}^{m} k_{\ell}$. Suppose that $\mathfrak{m}=\prod_{j=1}^{L}m(p_{j})$ (respectively, $\mathfrak{m}'=\prod^{L'}_{\ell=1}m(p'_{j})$) is a monomial in $\chi_{q}(\mathcal{S}_{2})$ (respectively, $\chi_{q}(\mathcal{S}_{1})$), where
$(p_{1}, \ldots, p_{L}) \in \overline{\mathscr{P}}_{(c_{j},d_{j})_{1 \leq j \leq L}}$ (respectively, $(p'_{1}, \ldots, p'_{L'}) \in \overline{\mathscr{P}}_{(c'_{j},d'_{j})_{1\leq j \leq L'}}$) is a tuple of non-overlapping paths, where $c'_{k_1+(\frac{k_{2}-1}{2})+(\frac{d_{i_{3}}}{d_{i_{3}+1}}k_{3})+1+\ell}=c_{k_1+k_2+k_3+\ell}$, $d'_{k_1+(\frac{k_{2}-1}{2})+(\frac{d_{i_{3}}}{d_{i_{3}+1}}k_{3})+1+\ell}=d_{k_1+k_2+k_3+\ell}$ for $1\leq \ell \leq \sum_{\ell=4}^{m} k_{\ell}$ and $c'_{\ell}=c_{\ell}$, $d'_{\ell}=d_{\ell}+4$ for $1\leq \ell \leq k_1$. That is, $$\mathscr{P}_{(c_{k_1+k_2+k_3+\ell},d_{k_1+k_2+k_3+\ell})}=\mathscr{P}_{(c'_{k_1+(\frac{k_{2}-1}{2})+(\frac{d_{i_{3}}}{d_{i_{3}+1}}k_{3})+1+\ell},
d'_{k_1+(\frac{k_{2}-1}{2})+(\frac{d_{i_{3}}}{d_{i_{3}+1}}k_{3})+1+\ell})}$$ for $1\leq \ell \leq \sum_{\ell=4}^{m} k_{\ell}$, see (a)--(c) in Figure \ref{proof 3}, relations between $\mathscr{P}_{(c_{\ell},d_{\ell})}$ for $k_1+1 \leq \ell \leq k_1+k_2+k_3$ and $\mathscr{P}_{(c'_{\ell},d'_{\ell})}$ for $k_1+1 \leq \ell \leq k_1+(\frac{k_{2}-1}{2})+(\frac{d_{i_{3}}}{d_{i_{3}+1}}k_{3})+1$ are shown in (d)--(f) in Figure \ref{proof 3}, and relations between $\mathscr{P}_{(c_{\ell},d_{\ell})}$ and $\mathscr{P}_{(c'_{\ell},d'_{\ell})}$ for $1\leq \ell \leq k_1$ are shown in (g) in Figure \ref{proof 4}.
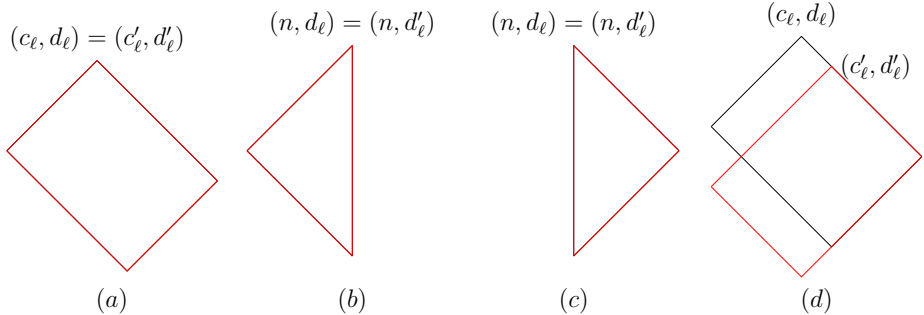
\begin{figure}[H]
\centering
\resizebox{0.8\width}{0.8\height}{
\begin{minipage}[b]{0.25\textwidth}
\centerline{
\begin{tikzpicture}
\draw (0,0)--(2,-2) (0,0)--(-1.5,-1.5);
\draw (-1.5,-1.5)--(0.5,-3.5)--(2,-2);
\draw[red] (0,0)--(2,-2) (0,0)--(-1.5,-1.5);
\draw[red] (-1.5,-1.5)--(0.5,-3.5)--(2,-2);
\node [above] at (0,0) {$(c_\ell,d_\ell)=(c'_\ell,d'_\ell)$};
\node at (0.25,-4) {$(a)$};
\end{tikzpicture}}
\end{minipage}
\begin{minipage}[b]{0.25\textwidth}
\centerline{
\begin{tikzpicture}
\draw (0.25,-0.25)--(-1.5,-2) (0.25,-0.25)--(0.25,-3.75);
\draw[red] (0.25,-0.25)--(-1.5,-2) (0.25,-0.25)--(0.25,-3.75);
\draw (-1.5,-2)--(0.25,-3.75);
\draw[red] (-1.5,-2)--(0.25,-3.75);
\node [above] at (0.25,-0.25) {$(n,d_\ell)=(n,d'_\ell)$};
\node at (0.25,-4.5) {$(b)$};
\end{tikzpicture}}
\end{minipage}
\begin{minipage}[b]{0.25\textwidth}
\centerline{
\begin{tikzpicture}
\draw (0.25,-0.25)--(2,-2) (0.25,-0.25)--(0.25,-3.75);
\draw[red]  (0.25,-0.25)--(2,-2) (0.25,-0.25)--(0.25,-3.75);
\draw (2,-2)--(0.25,-3.75);
\draw[red] (2,-2)--(0.25,-3.75);
\node [above] at (0.25,-0.25) {$(n,d_\ell)=(n,d'_\ell)$};
\node at (0.25,-4.5) {$(c)$};
\end{tikzpicture}}
\end{minipage}
\begin{minipage}[b]{0.25\textwidth}
\centerline{
\begin{tikzpicture}
\draw (0,0)--(2,-2) (0,0)--(-1.5,-1.5);
\draw (-1.5,-1.5)--(0.5,-3.5) -- (2,-2);
\draw [red] (0.5,-0.5)--(2,-2) (0.5,-0.5)--(-1.5,-2.5);
\draw[red] (-1.5,-2.5)--(0,-4) -- (2,-2);
\node [above] at (0,0) {$(c_\ell,d_\ell)$};
\node [right] at (0.5,-0.5) {$(c'_\ell,d'_\ell)$};
\node at (0.25,-4.4) {$(d)$};
\end{tikzpicture}}
\end{minipage}}
\end{figure}
\begin{figure}[H]
\centering
\resizebox{0.8\width}{0.8\height}{
\begin{minipage}[b]{0.3\textwidth}
\centerline{
\begin{tikzpicture}
\draw (0,0)--(2,-2) (0,0)--(-1.5,-1.5);
\draw (-1.5,-1.5)--(0.5,-3.5)--(2,-2);
\draw[red] (0.25,-0.25)--(2,-2)  (0.25,0.25)--(-1.5,-1.5) (0.25,0.25)--(0.25,-3.25);
\draw[red] (2,-2)--(0.25,-3.75)  (-1.5,-1.5)--(0.25,-3.25) (0.25,0.25)--(0.25,-3.75);
\node [left] at (0,0) {$(n-1,d_\ell)$};
\node [above] at (0.25,0.25) {$(n,d_\ell-1)$};
\node [right] at (0.25,-0.25) {$(n,d_\ell+1)$};
\node at (0.25,-4.2) {$(e)$};
\end{tikzpicture}}
\end{minipage}
\begin{minipage}[b]{0.3\textwidth}
\centerline{
\begin{tikzpicture}
\draw (0,0)--(-2,-2) (0,0)--(1.5,-1.5);
\draw (1.5,-1.5)--(-0.5,-3.5)--(-2,-2);
\draw[red] (-0.25,-0.25)--(-2,-2) (-0.25,0.25)--(1.5,-1.5) (-0.25,0.25)--(-0.25,-3.25);
\draw[red] (-2,-2)--(-0.25,-3.75) (1.5,-1.5)--(-0.25,-3.25) (-0.25,-0.25)--(-0.25,-3.75);
\node [right] at (0,0) {$(n-1,d_\ell)$};
\node [above] at (0.25,0.25) {$(n,d_\ell-1)$};
\node [left] at (-0.25,-0.25) {$(n,d_\ell+1)$};
\node at (0.25,-4.5) {$(f)$};
\end{tikzpicture}}
\end{minipage}
\begin{minipage}[b]{0.3\textwidth}
\centerline{
\begin{tikzpicture}
\draw (0,0)--(1.5,-1.5) (0,0)--(-1.5,-1.5);
\draw (-1.5,-1.5)--(0,-3) -- (1.5,-1.5);
\draw[red] (0,-0.8)--(1.5,-2.3) (0,-0.8)--(-1.5,-2.3);
\draw[red] (1.5,-2.3)--(0,-3.8) -- (-1.5,-2.3);
\node [above] at (0,0) {$(c_\ell,d_\ell)$};
\node [above] at (0,-1) {$(c'_\ell,d'_\ell)$};
\node at (0,-4.2) {$(g)$};
\end{tikzpicture}}
\end{minipage}}
\caption{Relations between $\mathscr{P}_{(c_{\ell},d_{\ell})}$ (black) and $\mathscr{P}_{(c'_{\ell},d'_{\ell})}$ (red).}\label{proof 3}
\end{figure}
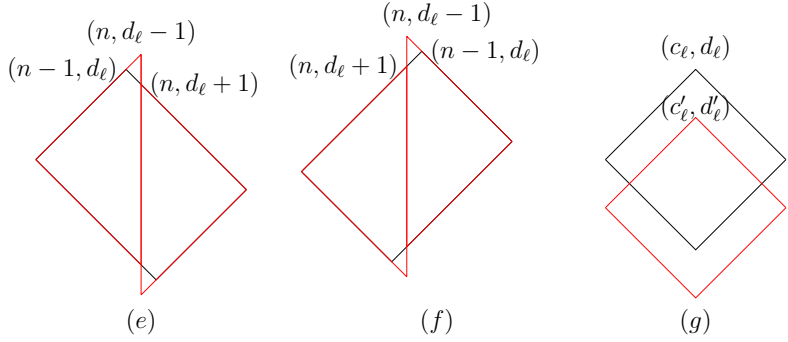

Suppose that $\mathfrak{m}\mathfrak{m}'$ is dominant. By the same arguments as the arguments in the proof of Theorem \ref{real snake modules} and relations between $\mathscr{P}_{(c_{\ell},d_{\ell})}$ and $\mathscr{P}_{(c'_{\ell},d'_{\ell})}$, $\mathfrak{m}\mathfrak{m}'$ is one of the dominant monomials $S_1S_2\prod_{j=0}^{r}A^{-1}_{i_{1}, t+4k_{1}-4j-2}$, $-1 \leq r \leq k_{1}-1$.

By the same arguments as the arguments of dominant monomials in $\chi_q(\mathcal{S}_1)\chi_q(\mathcal{S}_2)$, the dominant monomials in $\chi_q(\mathcal{S}_3)\chi_q(\mathcal{S}_4)$ must be of the form $\scalemath{0.92}{S_3S_4\prod_{j=0}^{r}A^{-1}_{i_{1}, t+4k_{1}-4j-2}}$, $-1 \leq r \leq k_{1}-2$, and the dominant monomial in $\chi_q(\mathcal{S}_5)\chi_q(\mathcal{S}_6)$ is $S_5S_6$.
\end{proof}

\begin{proof}[\bf Proof of the case corresponding to (18) in Table \ref{definition S_1 in type B_{n}}.] Suppose that $\mathcal{S}_1$ is of the form (18) in Table \ref{definition S_1 in type B_{n}}. Let
\begin{align*}
&\mathcal{S}_{1}= \mathcal{S}^{(t+2)}_{(k_{1}+1)^{(n,j_{1}-1)}, k_{2}^{(i_{2},j_{2})}, k_{3}^{(i_{3},j_{3})}, \ldots, k_{m}^{(i_{m})}},
&\mathcal{S}_{2}= \mathcal{S}^{(t)}_{{k_{1}}^{(n,j_{1})}, k_{2}^{(i_{2},j_{2})}, k_{3}^{(i_{3},j_{3})}, \ldots, k_{m}^{(i_{m})}},
\end{align*}
where $j_1 \geq 1$, $m \geq 2$, and $j_\ell \geq 0$ for $2 \leq \ell \leq m-1$.

Let $L=\sum_{\ell=1}^{m} k_{\ell}$. Suppose that $\mathfrak{m}=\prod_{j=1}^{L}m(p_{j})$ (respectively, $\mathfrak{m}'=\prod^{L+1}_{\ell=1}m(p'_{j})$) is a monomial in $\chi_{q}(\mathcal{S}_{2})$ (respectively, $\chi_{q}(\mathcal{S}_{1})$), where
$(p_{1}, \ldots, p_{L}) \in \overline{\mathscr{P}}_{(c_{j},d_{j})_{1 \leq j \leq L}}$ (respectively, $(p'_{1}, \ldots, p'_{L+1}) \in \overline{\mathscr{P}}_{(c'_{j},d'_{j})_{1\leq j \leq L+1}}$) is a tuple of non-overlapping paths, where $c_{k_1+\ell}=c'_{k_1+1+\ell}$, $d_{k_1+\ell}=d'_{k_1+1+\ell}$ for $1\leq \ell \leq L-k_1$ and $c_{\ell}=c'_{\ell}=c'_{k_1+1}=n$, $d_{\ell}+2=d'_{\ell}$ for $1\leq \ell \leq k_1$. That is, $\mathscr{P}_{(c_{k_1+\ell},d_{k_1+\ell})}=\mathscr{P}_{(c'_{k_1+1+\ell},d'_{k_1+1+\ell})}$ for $1\leq \ell \leq L-k_1$, see (a)--(c) in Figure \ref{proof 4} and relations between $\mathscr{P}_{(n,d_{\ell})}$ and $\mathscr{P}_{(n,d'_{\ell})}$ for $1\leq \ell \leq k_1$ are shown in (d), (e) in Figure \ref{proof 4}.

\begin{figure}[H]
\centering
\resizebox{0.8\width}{0.8\height}{
\begin{minipage}[b]{0.25\textwidth}
\centerline{
\begin{tikzpicture}
\draw (0,0)--(2,-2) (0,0)--(-1.5,-1.5);
\draw (-1.5,-1.5)--(0.5,-3.5)--(2,-2);
\draw[red] (0,0)--(2,-2) (0,0)--(-1.5,-1.5);
\draw[red] (-1.5,-1.5)--(0.5,-3.5)--(2,-2);
\node [above] at (0,0) {$(c_\ell,d_\ell)=(c'_\ell,d'_\ell)$};
\node at (0.25,-4) {$(a)$};
\end{tikzpicture}}
\end{minipage}
\begin{minipage}[b]{0.25\textwidth}
\centerline{
\begin{tikzpicture}
\draw (0.25,-0.25)--(-1.5,-2) (0.25,-0.25)--(0.25,-3.75);
\draw[red] (0.25,-0.25)--(-1.5,-2) (0.25,-0.25)--(0.25,-3.75);
\draw (-1.5,-2)--(0.25,-3.75);
\draw[red] (-1.5,-2)--(0.25,-3.75);
\node [above] at (0.25,-0.25) {$(n,d_\ell)=(n,d'_\ell)$};
\node at (0.25,-4.5) {$(b)$};
\end{tikzpicture}}
\end{minipage}
\begin{minipage}[b]{0.25\textwidth}
\centerline{
\begin{tikzpicture}
\draw (0.25,-0.25)--(2,-2) (0.25,-0.25)--(0.25,-3.75);
\draw[red]  (0.25,-0.25)--(2,-2) (0.25,-0.25)--(0.25,-3.75);
\draw (2,-2)--(0.25,-3.75);
\draw[red] (2,-2)--(0.25,-3.75);
\node [above] at (0.25,-0.25) {$(n,d_\ell)=(n,d'_\ell)$};
\node at (0.25,-4.5) {$(c)$};
\end{tikzpicture}}
\end{minipage}
\begin{minipage}[b]{0.25\textwidth}
\centerline{
\begin{tikzpicture}
\draw (0.25,-0.25)--(-1.5,-2) (0.25,-0.25)--(0.25,-3.75);
\draw[red] (0.25,-0.75)--(2,-2.5) (0.25,-0.75)--(0.25,-4.25);
\draw (-1.5,-2)--(0.25,-3.75);
\draw[red] (2,-2.5)--(0.25,-4.25);
\node [above] at (0.25,-0.25) {$(n,d_\ell)$};
\node [right] at (0.25,-0.75) {$(n,d'_\ell)$};
\node at (0.25,-4.5) {$(d)$};
\end{tikzpicture}}
\end{minipage}}
\end{figure}
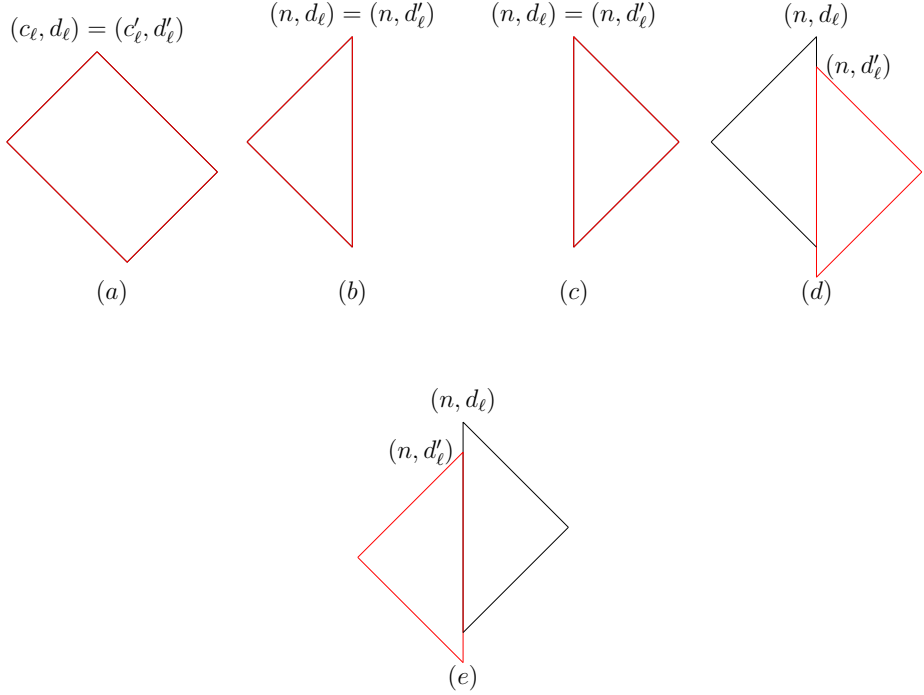
\begin{figure}[H]
\centering
\resizebox{0.8\width}{0.8\height}{
\begin{minipage}[b]{0.3\textwidth}
\centerline{
\begin{tikzpicture}
\draw (0.25,-0.25)--(2,-2) (0.25,-0.25)--(0.25,-3.75);
\draw[red]  (0.25,-0.75)--(-1.5,-2.5) (0.25,-0.75)--(0.25,-4.25);
\draw (2,-2)--(0.25,-3.75);
\draw[red] (-1.5,-2.5)--(0.25,-4.25);
\node [above] at (0.25,-0.25) {$(n,d_\ell)$};
\node [left] at (0.25,-0.75) {$(n,d'_\ell)$};
\node at (0.25,-4.5) {$(e)$};
\end{tikzpicture}}
\end{minipage}}
\caption{Relations between $\mathscr{P}_{(c_{\ell},d_{\ell})}$ (black) and $\mathscr{P}_{(c'_{\ell},d'_{\ell})}$ (red).}\label{proof 4}
\end{figure}

Suppose that $\mathfrak{m}\mathfrak{m}'$ is dominant. By the same arguments as the arguments in the proof of Theorem \ref{real snake modules} and relations between $\mathscr{P}_{(c_{\ell},d_{\ell})}$ and $\mathscr{P}_{(c'_{\ell},d'_{\ell})}$, $\mathfrak{m}\mathfrak{m}'$ is one of the dominant monomials $S_1S_2\prod_{j=0}^{r}A^{-1}_{i_{1}, t+2k_{1}-2j-1}$, $-1 \leq r \leq k_{1}-1$.

By the same arguments as the arguments of dominant monomials in $\chi_q(\mathcal{S}_1)\chi_q(\mathcal{S}_2)$, the dominant monomials in $\chi_q(\mathcal{S}_3)\chi_q(\mathcal{S}_4)$ must be of the form $\scalemath{0.95}{S_3S_4\prod_{j=0}^{r}A^{-1}_{i_{1}, t+2k_{1}-2j-1}}$, $-1 \leq r \leq k_{1}-2$, and the dominant monomial in $\chi_q(\mathcal{S}_5)\chi_q(\mathcal{S}_6)$ is $S_5S_6$.
\end{proof}

\section{Proof of Theorem \ref{equation simple}}  \label{proof irreducible}
In this section, we prove Theorem \ref{equation simple}.

By Lemma \ref{dominant monomials}, we have the following corollary.
\begin{corollary}
The modules in the second summand on the right-hand side of every equation of the $S$-systems for types $A_n$ and $B_n$ are special. In particular, they are simple.
\end{corollary}

Therefore, in order to prove Theorem \ref{equation simple}, we only need to show that the modules in the first summand on the right-hand side of every equation of the $S$-systems for types $A_n$ and $B_n$ are simple. Note that if $k_1=1$, then the modules in the first summand on the right-hand side of every equation of the $S$-systems for types $A_n$ and $B_n$ are special (has only a dominant monomial, see Lemma \ref{dominant monomials}) and hence they are simple. So, we consider only the case where $k_1\geq 2$.

Let $\mathcal{S}$ be a module corresponding to the first summand on the right-hand side of every equation of the $S$-systems for types $A_n$ and $B_n$. It suffices to prove that for each non-highest dominant monomial $M$ in $\mathcal{S}$, we have $\chi_q(L(M))\not \subseteq \chi_q(\mathcal{S})$, see \cite{Her06}, \cite{MY12a}.

We will prove that in the case of type $A_n$,
\begin{align}\label{one case of type An proof irreducible}
\mathcal{S}^{(t)}_{(k_{1}+1)^{(i_{1})}, {k_{2}}^{(i_{2}+1)}, k_{3}^{(i_{3},j_{3})}, k_{4}^{(i_{4},j_{4})}, \ldots, k_{m}^{(i_{m})}} \otimes \mathcal{S}^{(t+2)}_{(k_{1}-1)^{(i_{1})}, k_{2}^{(i_{2})}, k_{3}^{(i_{3},j_{3})}, k_{4}^{(i_{4},j_{4})}, \ldots, k_{m}^{(i_{m})}},
\end{align}
where $m\geq 3$, $i_1 > i_2$, $i_3 > i_2$, $j_\ell \geq 0$, $3\leq \ell \leq m-1$, is simple (this is $\mathcal{S}_3 \otimes \mathcal{S}_4$ which corresponds to the $5$-th line of the third column of Table \ref{dominant monomials in the S-system of type A}), and in the case of type $B_n$,
\begin{align} \label{one case of type Bn proof irreducible}
\begin{split}
\scalemath{0.96}{\mathcal{S}^{(t)}_{(k_{1}+1)^{(i_{1})}, (\frac{k_{2}-1}{2})^{(n-1)}, 1^{(n)}, (\frac{d_{i_{3}}}{d_{i_{3}+1}}k_{3})^{(i_{3}+1,j_{3}-\delta_{i_3i_4})}, k_{4}^{(i_{4},j_{4})},\ldots, k_{m}^{(i_{m})}}\otimes \mathcal{S}^{(t+4)}_{(k_{1}-1)^{(i_1)}, k_{2}^{(n)}, k_{3}^{(i_{3},j_{3})}, k_{4}^{(i_{4},j_{4})},\ldots, k_{m}^{(i_{m})}}}
\end{split}
\end{align}
is simple (this is $\mathcal{S}_3 \otimes \mathcal{S}_4$ which corresponds to the $14$-rd line of the third column of Table \ref{dominant monomials in the S-system of type B}). The other cases are similar.

\begin{proof}[\bf Proof of the fact that (\ref{one case of type An proof irreducible}) is simple.] By Lemma \ref{dominant monomials}, the dominant monomials appearing in $\chi_{q}(\mathcal{S}_{3}) \chi_{q}(\mathcal{S}_{4})$ are
\begin{align*}
M_{r} = S_3 S_4 \prod_{j=0}^{r} A^{-1}_{i_{1},t+2k_{1}-2j-1}, \quad  -1 \leq r \leq k_{1}-2,
\end{align*}
where
\begin{gather}
\begin{align*}
S_3 & =(i_1)_{t}(i_1)_{t+2}\cdots (i_1)_{t+2k_1} \prod_{r=0}^{k_{2}-1} (i_{2}+1)_{t+2k_1+i_1-i_2+2r+1} \prod_{j=3}^{m} \left( \prod_{r=0}^{k_{j}-1} (i_{j})_{t+2r+\sum_{\ell=1}^{j-1}n_{\ell}} \right), \\
S_4 & = (i_1)_{t+2}(i_1)_{t+4}\cdots (i_1)_{t+2k_1-2} \prod_{r=0}^{k_{2}-1} (i_{2})_{t+2k_1+i_1-i_2+2r} \prod_{j=3}^{m} \left( \prod_{r=0}^{k_{j}-1} (i_{j})_{t+2r+\sum_{\ell=1}^{j-1}n_{\ell}} \right),
\end{align*}
\end{gather}
and $n_{\ell}$ is defined in (\ref{n_1 in type A_n}). We need to show that $\chi_{q}(M_{r}) \nsubseteq \chi_{q}(S_3)\chi_{q}(S_4)$ for $0 \leq r \leq k_{1}-2$.

We have
\begin{align*}
\beta_{i_1}(M_r)& = {\bf i_1}^{(1)} {\bf i_1}^{(2)} {\bf i_1}^{(3)} \cdots {\bf i_1}^{(p)},
\end{align*}
where $p$ is some integer,
\begin{align*}
& {\bf i_1}^{(1)}=(i_1)_{t} (i_1)_{t+2} (i_1)_{t+4} \cdots (i_1)_{t+2k_1-2r-4} (i_1)_{t+2k_1-2r-2}, \\
& {\bf i_1}^{(2)} = (i_1)_{t+2} (i_1)_{t+4} \cdots (i_1)_{t+2k_1-2r-4},
\end{align*}
and ${\bf i_1}^{(1)}, \ldots, {\bf i_1}^{(p)}$ are $q_{i_1}$-strings which are pairwise in general position.

Let $n_{r}= M_{r} A^{-1}_{i_{1},t+2k_{1}-2r-1}$. By Corollary \ref{Uqsl_2 arguments}, the monomial $n_r \in \chi_{q}(M_r)$. We have
\begin{align*}
n_{r} & = M_{r} A^{-1}_{i_{1},t+2k_{1}-2r-1}  \\
& = \left(S_{3} S_{4}\prod_{j=0}^{r} A^{-1}_{i_{1},t+2k_{1}-2j-1}\right) A^{-1}_{i_{1},t+2k_{1}-2r-1} \\
& = \left(S_{3}  A^{-1}_{i_{1},t+2k_{1}-2r-1}\right) \left(S_{4}\prod_{j=0}^{r} A^{-1}_{i_{1},t+2k_{1}-2j-1}\right).
\end{align*}

By Corollary \ref{Uqsl_2 arguments}, the monomial $S_{4}\prod_{j=0}^{r} A^{-1}_{i_{1},t+2k_{1}-2j-1}$ is in $\chi_q(\mathcal{S}_{4})$, since
$$\beta_{i_1}(S_{4})={\bf i_1}^{(1)} \cdots {\bf i_1}^{(p)},$$
where $p$ is some integer, ${\bf i_1}^{(1)}=(i_1)_{t+2} (i_1)_{t+4} \cdots (i_1)_{t+2k_1-2}$ and ${\bf i_1}^{(1)}, \ldots, {\bf i_1}^{(p)}$ are $q_{i_1}$-strings which are pairwise in general position. Since $(i_1)_{t+2k_{1}-2r-2}$ is not a factor of $S_{4}\prod_{j=0}^{r} A^{-1}_{i_{1},t+2k_{1}-2j-1}$, the monomial
\begin{align}\label{9.3}
\left(S_{4}\prod_{j=0}^{r} A^{-1}_{i_{1},t+2k_{1}-2j-1}\right)A^{-1}_{i_{1},t+2k_{1}-2r-1}
\end{align}
is not in $\varphi_{i_1}\left(S_{4}\prod_{j=0}^{r} A^{-1}_{i_{1},t+2k_{1}-2j-1}\right)$. Therefore (\ref{9.3}) is not in $\chi_q(\mathcal{S}_{4})$ by Lemma \ref{need to expand from right most factors}.

Suppose that $n_{r} \in \chi_{q}(S_3)\chi_{q}(S_4)$. Then $S_{3}  A^{-1}_{i_{1},t+2k_{1}-2r-1}$ would be in $\chi_q(\mathcal{S}_{3})$.
We have
\[
\beta_{i_1}(S_{3})={\bf i_1}^{(1)} \cdots {\bf i_1}^{(p)},
\]
where $p$ is some integer, ${\bf i_1}^{(1)}= (i_1)_{t}(i_1)_{t+2}\cdots (i_1)_{t+2k_1}$, ${\bf i_1}^{(1)}, \ldots, {\bf i_1}^{(p)}$ are $q_{i_1}$-strings which are pairwise in general position, and ${\bf i}^{(j')}_1 \not\subset {\bf i}^{(j'')}_1$, $j' \neq j''$, $j', j'' \in \{1, \ldots, p\}$ ($S_3$ is a snake module). Since $(i_1)_{t+2k_1}$ is a factor of $S_{3}  A^{-1}_{i_{1},t+2k_{1}-2r-1}$, by Lemma \ref{need to expand from right most factors} we have $S_{3}  A^{-1}_{i_{1},t+2k_{1}-2r-1}\not\in \chi_q(S_{3})$. This is a contradiction.

Therefore $n_r \not\in \chi_{q}(S_3)\chi_{q}(S_4)$ and hence $\chi_{q}(M_r)\nsubseteq \chi_q(S_3) \chi_q(S_4)$.
\end{proof}

\begin{proof}[\bf Proof of the fact that (\ref{one case of type Bn proof irreducible}) is simple.] By Lemma \ref{dominant monomials}, the dominant monomials appearing in $\chi_q(\mathcal{S}_3)$ $\chi_q(\mathcal{S}_4)$ are
\begin{align*}
M_{r} = S_3 S_4 \prod_{j=0}^{r} A^{-1}_{i_{1},t+4k_{1}-4j-2}, \quad  -1 \leq r \leq k_{1}-2.
\end{align*}

We need to show that $\chi_{q}(M_{r})\nsubseteq \chi_q(S_3) \chi_q(S_4)$ for $0 \leq r \leq k_1-2$. By the expression $M_r$ and Corollary \ref{Uqsl_2 arguments}, the monomial $n_{r} = M_{r} A^{-1}_{i_{1},t+4k_{1}-4r-2}$ is in $\chi_{q}(M_r)$, since
\begin{align*}
\beta_{i_1}(M_r)& = {\bf i_1}^{(1)} {\bf i_1}^{(2)} {\bf i_1}^{(3)} \cdots {\bf i_1}^{(p)},
\end{align*}
where $p$ is some integer,
\begin{align*}
& {\bf i_1}^{(1)}=(i_1)_{t} (i_1)_{t+4} (i_1)_{t+8} \cdots (i_1)_{t+4k_1-4r-8} (i_1)_{t+4k_1-4r-4}, \\
& {\bf i_1}^{(2)} = (i_1)_{t+4} (i_1)_{t+8} \cdots (i_1)_{t+4k_1-4r-8},
\end{align*}
and ${\bf i_1}^{(1)}, \ldots, {\bf i_1}^{(p)}$ are $q_{i_1}$-strings which are pairwise in general position.

By Corollary \ref{Uqsl_2 arguments}, the monomial $S_4 \prod_{j=0}^{r} A^{-1}_{i_{1},t+4k_{1}-4j-2}$ is in $\chi_q(\mathcal{S}_{4})$, since
$$\beta_{i_1}(S_{4})={\bf i_1}^{(1)} \cdots {\bf i_1}^{(p)},$$ where $p$ is some integer, ${\bf i_1}^{(1)}=(i_1)_{t+4} (i_1)_{t+8} \cdots (i_1)_{t+4k_1-4}$ and ${\bf i_1}^{(1)}, \ldots, {\bf i_1}^{(p)}$ are $q_{i_1}$-strings which are pairwise in general position.

Suppose that $n_{r}\in \chi_{q}(S_3)\chi_{q}(S_4)$. Then $n_r=m_{1}m_{2}$, where $m_1 \in \chi_{q}(S_3)$, $m_2 \in \chi_{q}(S_4)$. Since $n_r=\left(S_3 A^{-1}_{i_{1},t+4k_{1}-4r-2}\right) \left(S_4 \prod_{j=0}^{r} A^{-1}_{i_{1},t+4k_{1}-4j-2}\right)$ and $(i_{1})_{t+4k_{1}-4r-4}$ is not a factor of $S_4 \prod_{j=0}^{r}A^{-1}_{i_{1},t+4k_{1}-4j-2}$, by the expressions $S_3$ and $S_4$ we must have
\begin{align*}
m_1 =S_3 A^{-1}_{i_{1},t+4k_{1}-4r-2} \in \chi_{q}(S_3), \quad m_{2} = S_4 \prod_{j=0}^{r} A^{-1}_{i_{1},t+4k_{1}-4j-2}\in \chi_{q}(S_4),
\end{align*}
which contradicts Lemma \ref{need to expand from right most factors}:
$$\beta_{i_1}(S_{3})={\bf i_1}^{(1)} \cdots {\bf i_1}^{(p)},$$
where $p$ is some integer, ${\bf i_1}^{(1)}=(i_{1})_{t}(i_{1})_{t+4}\ldots(i_{1})_{t+4k_{1}}$, ${\bf i_1}^{(1)}, \ldots, {\bf i_1}^{(p)}$ are $q_{i_1}$-strings which are pairwise in general position and ${\bf i}^{(j')}_1 \not\subset {\bf i}^{(j'')}_1$, $j' \neq j''$, $j', j'' \in \{1, \ldots, p\}$ ($S_3$ is a snake module). Since $(i_{1})_{t+4k_{1}}$ is a factor of $S_3 A^{-1}_{i_{1},t+4k_{1}-4r-2}$, by Lemma \ref{need to expand from right most factors} we have $S_3 A^{-1}_{i_{1},t+4k_{1}-4r-2}\not\in \chi_q(S_{3})$. Therefore $n_r \not\in \chi_{q}(S_3)\chi_{q}(S_4)$ and hence $\chi_{q}(M_r)\nsubseteq \chi_q(S_3) \chi_q(S_4)$.
\end{proof}

\section*{Acknowledgement}
The authors are very grateful to the anonymous referee for the comments and suggestions that have been very helpful to improve the quality of this paper. The authors would like to express their gratitude to Professor Vyjayanthi Chari for helpful discussions about prime modules. This work was partially supported by the National Natural Science Foundation of China (no. 11771191, 11371177, 11501267, 11401275), and the Fundamental Research Funds for the Central Universities of China (no. lzujbky-2015-78). The research of J.-R. Li on this project is supported by the Minerva foundation with funding from the Federal German Ministry for Education and Research, ERC AdG Grant 247049, and the PBC Fellowship Program of Israel for Outstanding Post-Doctoral Researchers from China and India.


\begin{bibdiv}
\begin{biblist}

\bib{Car05}{book}{
author= {Carter, R. W.},
title={Lie algebras of finite and affine type},
series={Cambridge Studies in Advanced Mathematics},
volume={96},
publisher={Cambridge University Press},
address={Cambridge},
year={2005}}

\bib{C95}{article}{
author={Chari, V.},
title={Minimal affinizations of representations of quantum groups: the rank $2$ case}, journal={Publ. Res. Inst. Math. Sci.},
date={1995},
volume={31},
number={5},
pages={873--911}}

\bib{CMY13}{article}{
author={Chari, V.},author={Moura, A.},author={Young, C. A.},
title={Prime representations from a homological perspective},
journal={Math. Z.},
date={2013},
volume={274},
number={1--2},
pages={613--645}}

\bib{CP91}{article}{
author={Chari, V.},author={Pressley, A.},
title={Quantum affine algebras},
journal={Comm. Math. Phys.},
date={1991},
volume={142},
number={2},
pages={261--283}}

\bib{CP94}{book}{
author={Chari, V.},author= {Pressley, A.},
title={A guide to quantum groups},
publisher={Cambridge University Press},
address={Cambridge},
year={1994}}

\bib{CP95a}{article}{
author={Chari, V.},author= {Pressley, A.},
title={Quantum affine algebras and their representations},
conference={
title={Representations of groups},
address={Banff, AB},
date={1994}},
book={
series={CMS Conf. Proc.},
volume={16},
publisher={Amer. Math. Soc.},
address={Providence, RI},
date={1995}},
pages={59--78}}

\bib{CP96}{article}{
author={Chari, V.},author={Pressley, A.},
title={Quantum affine algebras and affine Hecke algebras},
journal={Pacific J. Math.},
date={1996},
volume={174},
number={2},
pages={295--326}}

\bib{CP97}{article}{
author={Chari, V.},author= {Pressley, A.},
title={Factorization of representations of quantum affine algebras},
conference={
title={Modular interfaces},
address={Riverside, CA},
date={1995}},
book={series={AMS/IP Stud. Adv. Math.},
volume={4},
publisher={Amer. Math. Soc.},
address={Providence, RI},
date={1997}},
pages={33--40}}

\bib{Dri87}{article}{
author={Drinfeld, V. G.},
title={Quantum groups},
conference={
title={Proceedings of the International Congress of Mathematicians},
address={Berkeley, Calif.},
date={1986}},
book={volume={1,2},
publisher={Amer. Math. Soc.},
address={Providence, RI},
date={1987}},
pages={798--820}}

\bib{Dri88}{article}{
author={Drinfeld, V. G.},
title={A new realization of Yangians and of quantum affine algebras},
journal={Dokl. Akad. Nauk SSSR},
date={1987},
volume={296},
number={1},
pages={13--17 [Russian]; translation in Soviet Math. Dokl. \textbf{36} (1988), no. 2, 212--216}}

\bib{FM01}{article}{
author={Frenkel, E.},author={Mukhin, E.},
title={Combinatorics of $q$-characters of finite-dimensional representations of quantum affine algebras},
journal={Comm. Math. Phys.},
date={2001},
volume={216},
number={1},
pages={23--57}}

\bib{FR98}{article}{
author={Frenkel, E.},author={Reshetikin, N.},
title={The $q$-characters of representations of quantum affine algebras and deformations of $W$-algebras},
conference={title={Recent developments in quantum affine algebras and related topics},
address={Raleigh, NC},
date={1998}},
book={series={Contemp. Math.},
volume={248},
publisher={Amer. Math. Soc.},
address={Providence, RI},
date={1999}},
pages={163--205}}


\bib{FZ02}{article}{
author={Fomin, S.},author={Zelevinsky, A.},
title={Cluster algebras I: Foundations},
journal={J. Amer. Math. Soc.},
date={2002},
volume={15},
number={2},
pages={497--529}}

\bib{GG14}{article}{
author={Grabowski, J.},author={Gratz,S.},
title={Cluster algebras of infinite rank, with an appendix by Michael Groechenig},
journal={J. Lond. Math. Soc.},
date={2014},
volume={89},
number={2},
pages={337--363}}

\bib{Her05}{article}{
author={Hernandez, D.},
title={Monomials of $q$ and $q,t$-characters for non simply-laced quantum affinizations},
journal={Math. Z.},
date={2005},
volume={250},
number={2},
pages={443--473}}

\bib{Her06}{article}{
author={Hernandez, D.},
title={The Kirillov--Reshetikhin conjecture and solutions of T-systems},
journal={J. Reine Angew. Math.},
date={2006},
volume={596},
number={2},
pages={63--87}}

\bib{HL10}{article}{
author={Hernandez, D.},author={Leclerc, B.},
title={Cluster algebras and quantum affine algebras},
journal={Duke Math. J.},
date={2010},
volume={154},
number={2},
pages={265--341}}

\bib{HL13}{article}{
author={Hernandez, D.},author={Leclerc, B.},
title={A cluster algebra approach to $q$-characters of Kirillov-Reshetikhin modules},
journal={J. Eur. Math. Soc. (JEMS)},
date={2016},
volume={18},
number={5},
pages={1113--1159}}

\bib{Jim85}{article}{
author={Jimbo, M.},
title={A $q$-difference analogue of $U(\mathfrak{g})$ and the Yang-Baxter equation},
journal={Lett. Math. Phys.},
date={1985},
volume={10},
number={1},
pages={63--69}}

\bib{K12}{article}{
author={Keller, B.},
title={Quiver mutation in Java},
journal={http://people.math.jussieu.fr/keller/quivermutation/},
date={2006, accessed April 2012}}

\bib{KQ14}{article}{
author={Kimura, Y.},author={Qin, F.},
title={Graded quiver varieties, quantum cluster algebras and dual canonical basis},
journal={Adv. Math.},
date={2014},
volume={262},
pages={261--312}}

\bib{KNS94}{article}{
author={Kuniba, A.},author={Nakanishi, T.},author={Suzuki, J.},
title={Functional relations in solvable lattice models: I. Functional relations and representation theory},
journal={Internat. J. Modern Phys. A},
date={1994},
volume={9},
number={30},
pages={5215--5266}}

\bib{LM16}{article}{
author={Lapid, E.},author={Minguez, A.},
title={Geometric conditions for $\square$-irreducibility of certain representations of the general linear group over a non-archimedean local field},
journal={Adv. Math.},
date={2018},
volume={339},
pages={113--190}} 

\bib{Le03}{article}{
author={Leclerc, B.},
title={Imaginary vectors in the dual canonical basis of $U_{q}(n)$},journal={Transform Groups},
date={2003},
volume={8},
number={1},
pages={95--104}}

\bib{Le10}{article}{
author={Leclerc, B.},
title={Cluster algebras and representation theory},
conference={
title={Proceedings of the International Congress of Mathematicians}},
book={series={Hindustan Book Agency, New Delhi},
volume={IV},
date={2010}},
pages={2471--2488}}

\bib{Lee13}{article}{
author={Lee, K.},
title={Every finite acyclic quiver is a full subquiver of a quiver mutation equivalent to a bipartite quiver},
journal={arXiv:1311.0711},
date={2013},
pages={1--2}}

\bib{LQ17}{article}{
author={Li, J.-R.},author={Qiao, L.},
title={Three-term recurrence relations of minimal affinizations of type $G_2$},
journal={Journal of Lie Theory},
date={2017},
volume={27},
number={4},
pages={1119--1140}}

\bib{MY12a}{article}{
author={Mukhin, E.},author={Young, C. A. S.},
title={Path description of type $B$ $q$-characters},
journal={Adv. Math.},
date={2012},
volume={231},
number={2},
pages={1119--1150}}

\bib{MY12b}{article}{
author={Mukhin, E.},author={Young, C. A. S.},
title={Extended $T$-systems},
journal={Selecta Math. (N.S.)},
date={2012},
volume={18},
number={3},
pages={591--631}}

\bib{Nak11}{article}{
author={Nakajima, H.},
title={Quiver varieties and cluster algebras},
journal={Kyoto J. Math.},
date={2011},
volume={51},
number={1},
pages={71--126}}

\bib{NN11}{article}{
author={Nakai, W.},author={Nakanishi, T.},
title={On Frenkel-Mukhin algorithm for $q$-character of quantum affine algebras},
conference={
title={Exploring new structures and natural constructions in mathematical physics}},
book={series={Adv. Stud. Pure Math.},
volume={61},
publisher={Math. Soc. Japan},
address={Tokyo},
date={2011}},
pages={327--347}}

\bib{Q15}{article}{
author={Qin, F.},
title={Triangular bases in quantum cluster algebras and monoidal categorification conjectures},
journal={Duke Math. J.},
date={2017},
volume={166},
number={12},
pages={2337--2442}}

\bib{YMLZ15}{article}{
author={Yang, Y.-M.},author={Ma, H.-T.},author={Lin, B.-S.},author={Zheng, Z.-J.},
title={Cluster algebra structure on the finite dimensional representations of affine quantum group $U_{q}(\widehat{A}_3)$},
journal={Chinese Phys. B},
date={2015},
volume={24},
number={1},
pages={010201}}

\bib{ZDLL16}{article}{
author={Zhang, Q.-Q.},author={Duan, B.},author={Li, J.-R.},author={Luo, Y.-F.},
title={M-systems and cluster algebras},
journal={Int. Math. Res. Not. IMRN},
volume={2016},
number={14},
pages={4449--4486}}

\end{biblist}
\end{bibdiv}
\end{document}